\documentclass[preprint,3p,12pt]{elsarticle}

\usepackage{graphicx}
\usepackage{amsmath}
\usepackage{amsfonts}
\usepackage{enumitem} 	\usepackage{xspace}						% Resume emurate 		
\usepackage[labelfont=bf,tableposition=top]{caption} 	% Make captions blod (Figure & Table)
\usepackage{subfig}									% Subfigure package
\usepackage{epstopdf}								% Eps to pdf, now it is possible to use both PNG and EPS images
\usepackage{floatrow}
\newfloatcommand{capbtabbox}{table}[][\FBwidth]
\newsavebox{\bigleftbox}
\floatstyle{plaintop}
\restylefloat{table}
\usepackage{color}
\usepackage{xcolor}
\usepackage{mathrsfs}  
\usepackage{multirow}
\usepackage{amsmath}
\usepackage[capitalise]{cleveref}
\usepackage{amsthm}
\usepackage{titlesec}
\usepackage{etoolbox}
\usepackage{lipsum}
\usepackage{chngpage}
\usepackage{wrapfig}
\usepackage{tikz}
\usetikzlibrary{arrows.meta, positioning, calc}
\usepackage{pgfplots}

\usepackage{float}
\usepackage[section]{placeins}
\theoremstyle{remark}
\newtheorem{remark}{Remark}
\newcommand{\bair}{e^{q}}

\newcommand{\qair}{q}
\newcommand{\cair}{\phi}   %\phi reminder

\newcommand{\capacitor}{C}
\newcommand{\cotwo}{CO\textsubscript{2}\xspace}

\newcommand{\bth}{e^{R}}

\begin{document}

\begin{frontmatter}

\title{\large Data-driven moving-window Bayesian inference for transient \cotwo--temperature network models of buildings}

\author[address1]{Zhijian Wang\corref{cor1}}
\ead{z.wang7@tue.nl}
\author[address1]{Stein K.F. Stoter}
\ead{k.f.s.stoter@tue.nl}
\author[address1]{Clemens V. Verhoosel}
\ead{c.v.verhoose@tue.nl}
\author[address1]{Idoia Cortes Garcia}
\ead{i.cortes.garcia@tue.nl}

\cortext[cor1]{Corresponding author}

\address[address1]{Department of Mechanical Engineering, Eindhoven University of Technology, The Netherlands}
%\address[address2]{Adress for author 2 }

\begin{abstract}
%Indoor thermal dynamics in buildings are strongly affected by ventilation and occupancy, yet these drivers are not easily accessible by standard building instrumentation. This limit challenges the building thermal operation systems, especially when those states change regularly such as in offices.

In this work, we proposes a \cotwo--temperature network model that links multi-zone mass transport and thermal dynamics through shared latent drivers, airflow and occupancy. The thermal component is formulated as a resistance-capacitance (RC) network augmented with airflow-driven convective exchange, while the \cotwo component is governed by inter-zonal convective transport. To calibrate the model and track time-varying operating conditions based on sparse sensing, we introduce a moving-window Bayesian inference procedure that jointly estimates thermal parameters, airflow and occupancy trajectories. The estimation also provides room-specific sensor noise levels, yielding posterior predictive forecasts with credible intervals.

The framework is assessed using a controlled synthetic benchmark, and a scaled physical validation experiment using \cotwo and temperature sensing. In both settings, the posterior accurately reconstructs trajectories within windows and delivers low forecast errors. When inference windows overlap abrupt regime transitions, the widened uncertainty bands and increased inferred noise levels provide an interpretable diagnostic of model--data mismatch, followed by rapid recovery once the new regime is observed. 
%Overall, coupling \cotwo-informed airflow with thermal dynamics improves robustness of temperature prediction, supporting practical monitoring and energy-performance assessment under limited sensing.
Overall, coupling \cotwo-informed airflow with thermal dynamics yields a robust approach for conductive and advective temperature prediction, supporting practical monitoring and energy-performance assessment under limited sensing.
\end{abstract}


\end{frontmatter}
\newpage
\tableofcontents

%\linenumbers

\newpage

%=====================================================
\section{Introduction}
\label{sec:introduction}

% building consumes energy
Buildings account for approximately 40\% of energy use in Europe~\cite{EU_EPBD_2024}, 70\% of which is attributed to heating, ventilation, and air conditioning (HVAC) systems that maintain indoor thermal comfort and indoor air quality~\cite{muncan2024state}. Consequently, HVAC control plays a critical role in optimizing building energy consumption~\cite{Wang2025DynamicFacade}. While low-cost sensors for, e.g., temperature, humidity, \cotwo, and particulate matter, are nowadays commonly used for monitoring indoor air quality and thermal conditions, the full potential of these sensor streams for real-time control remains largely unutilized. Most building condition regulation systems adhere to coarse, long-term schedules rather than adaptive, data-driven strategies~\cite{esrafilian2021occupancy}. Integrated use of sensor information in closed-loop HVAC optimization can lead to substantial energy saving and is therefore warranted. Our aim is to create a framework that combines the robustness of models based on physical laws with real-time sensor data that can ultimately be used for improved thermal control of office buildings.

In typical office buildings, conductive heat transfer through walls, windows, and other boundaries, together with convective exchange via inter-room airflow, jointly determine both temperature and indoor air quality, with occupancy patterns driving these dynamics. {Incorporating} natural ventilation, inter-zone airflow, and occupancy into the control system remains challenging and is often often handled indirectly through static schedules and simplified assumptions, which can lead to suboptimal operation \cite{Hobson2021_OCC,Hahn2022_InformationGap}. A practical alternative to direct airflow or occupancy measurement is to infer these {latent} drivers from the dynamics observed in the temperature and \cotwo concentration measurements. This offers a non-intrusive indication that can be embedded in supervisory control. As typical sources of airflow in office buildings are ventilation% for occupancy
, adding \cotwo as an inference proxy within the HVAC loop thus has the potential for large efficiency gains without additional invasive instrumentation. 

Given the above motivation, we aim to develop a model that captures coupled temperature and airflow dynamics while remaining lightweight enough for real-time use in office buildings. Resistance--capacitance (RC) network models are widely used for thermal identification and control because they balance physical interpretability and computational speed. Comprehensive reviews on modeling choices, parameterization, and applications across simulation and control pipelines can be found in \citep{Zhang2021RSER,Hou2021RSER,Chong2018EAB_Guidelines}. Recent case studies demonstrate practical workflows for parameter identification and controller tuning in real buildings \citep{Rivera2022Frontiers}, and newer approaches emphasize physical coherence so that calibrated parameters remain consistent with building physics \citep{Lecocq2025EnBuild}. However, these RC-focused works primarily address conductive heat transfer and typically do not explicitly represent inter-zonal airflow and its impact on heat transport \citep{Dols2016CONTAM}, which becomes critical in multi-zone settings with ventilation and air exchange.

\cotwo mass-balance and transport models are widely used to infer indoor occupancy and ventilation conditions when direct airflow measurement is impractical. Gray-box methods, which refers to low-order, physics-based dynamic models, can recover time-varying occupancy and ventilation rates directly from \cotwo trajectories \citep{Wolf2019ApEn,Wolf2019CLIMA}, and simulation-backed studies suggest that combining an airflow-network structure with \cotwo sensing enables multi-zone occupancy estimation at the office scale \citep{Jorissen2017BS}. Field deployments using low-cost sensors further report robust estimation performance across diverse scenarios \citep{Liang2024JBE}.

Bayesian inference has also been widely used to calibrate building energy and indoor-environment models while quantifying parameter uncertainty \citep{Hou2021RSER,Chong2018EAB_Guidelines}. Beyond one-time calibration, several studies propose sequential, continuous, or moving-window updates to track parameter drift and regime changes over time \citep{Chong2019IBPSA,Chong2019EAB_Continuous}. In the context of air-quality sensing, Bayesian formulations have also been applied to infer latent occupancy and ventilation drivers from \cotwo measurements, including filtering-based approaches \citep{Jiang2020EnB} and moving-window methods that estimate time-varying occupancy and supply flows using \cotwo and ventilation data \citep{Zavala2013ANL}.

Although RC thermal models, \cotwo-based occupancy/ventilation inference, and Bayesian calibration have each been studied, they are usually developed in isolation. Occupancy or ventilation is typically inferred from \cotwo data alone, while RC parameters are calibrated separately, and the two are rarely coupled. Moreover, online multi-zone methods often assume a fixed measurement-noise level and return only point estimates, which can be unreliable when data are sparse or sensor quality varies. As a result, such estimators are more vulnerable to model mismatch, sensor drift, and regime changes.

In the current work, we develop a model to forecast temperature, based on historical data while inferring occupancy, airflow, and heat transfer parameters. Rather than aiming to represent the full operation of HVAC systems, including active heating and cooling control, the framework focuses on the coupled evolution of indoor temperature and \cotwo concentration, where occupancy driven \cotwo production and airflow provide additional information for the thermal state inference. A Bayesian inference framework is proposed for a coupled \cotwo--temperature network, that makes use of \cotwo and temperature measurements to update heat transfer, airflow, and occupancy parameters online using a moving window approach. Hereby, the model is able to capture both conductive and convective heat transfer through zones of the building. The key contributions of this work are: (i) a joint likelihood formulation that integrates both airflow transport and heat transfer, (ii) a windowed Bayesian updating scheme that captures gradual system changes while limiting parameter drift, and (iii) a robust inference design that remains effective under model mismatch and sparse sensing, while explicitly estimating noise and model uncertainty to reveal prediction reliability. 

This paper is outlined as follows: Section \ref{sec:Model} defines the coupled \cotwo and temperature networks and discusses the model coupling via occupancy and airflow conditions.
Section \ref{sec:Bayesian} introduces the setup of the Bayesian inference framework with moving window parameter updating. 
Section \ref{sec:Synthetic benchmark problem} presents a synthetic benchmark problem based on a characteristic building layout, which we use to quantify the performance of our model and the Bayesian inference scheme in a controlled environment. 
In Section \ref{sec:Experiment validation} we present a physical experiment of a (scaled) setup with the same layout as the synthetic benchmark to validate our model. %, and reports the results.
Section \ref{sec:conclusion} draws the conclusions of the paper.

%===============================================
%===============================================
\section{\cotwo--temperature network model}
\label{sec:Model}

\begin{figure}[H]
    \centering
    \includegraphics[width=0.7\linewidth]{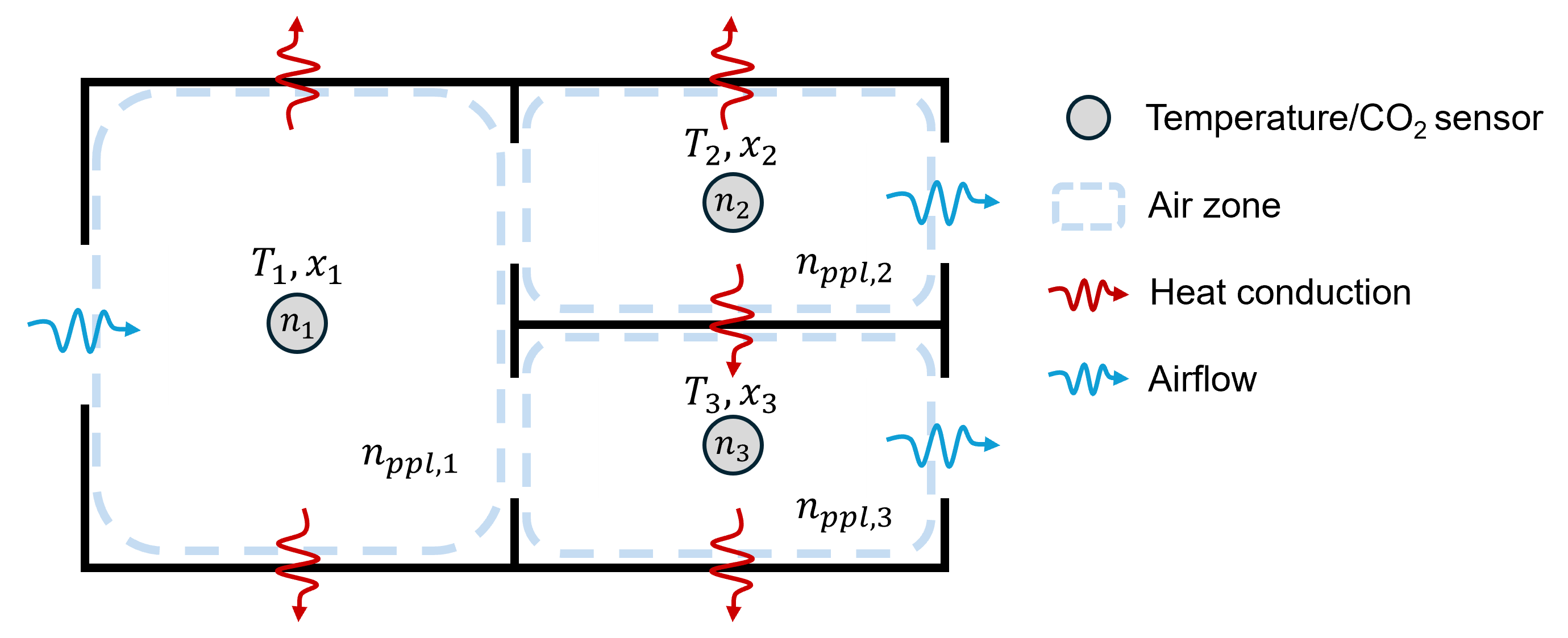}
    \caption{Schematic of the room layout.}
    \label{fig:Room schematic}
\end{figure}

In buildings, air motion simultaneously transports heat and gaseous pollutants such as \cotwo generated by occupants, while occupancy numbers largely determine heat and \cotwo production. To model the dynamics of both the temperature and the \cotwo concentration, we consider a building divided into zones, as shown in Figure~\ref{fig:Room schematic}. One zone typically, though not necessarily, represents one room. Each zone \(n_i\) assumes uniform conditions, and employs a sensor that reports the temperature $T_i$ and \cotwo concentration \(\cair_i\). The complete system is then modeled as two lumped networks, a \cotwo network and a temperature network, containing the same set of zones (nodes). 
Between the zones, there is advective exchange in the \cotwo network and both convective and conductive heat transfer in the temperature network. Which zones exchange \cotwo or heat either advectively or conductively is represented by connections between nodes (edges), which do not necessarily coincide between the two networks. %\icg{Whereas all connections with advective exchange included in the \cotwo network also have conductive heat transfer, some connections (e.g. between walls) only have conductive heat transfer and no advective one. Therefore, all edges of the \cotwo network are also part of the temperature network but not the other way around.} 
The latent, time-varying drivers are the inter-zonal airflows and zonal occupancy numbers. They couple the two networks through airflow-driven \cotwo advection and airflow-driven heat advection. Additionally, the temperature network involves conductive heat transfer parameters.

To formalize the model, we define the sets of nodal states as $\boldsymbol{\cair}=\{\cair_i\}$ and $\boldsymbol{T}=\{T_i\}$. %, the signed edge flows as $\boldsymbol{\qair}=\{\qair_{i}\}$ (one for each edge in the convection network), the occupancies as $\boldsymbol{n}_{\mathrm{ppl}}=\{n_{\mathrm{ppl},i}\}$. 
Denoting the ``unknown'' (to be inferred) parameters that determine the \cotwo and temperatuer models as $\boldsymbol{\theta}^{\cair}$ and $\boldsymbol{\theta}^{T}$, and denoting the ``known'' model parameters as $\boldsymbol{\kappa}^{\cair}$ and $\boldsymbol{\kappa}^{T}$, the forward models are expressed as the residuals
\begin{subequations}\label{eq:residuals}
\begin{align}
\boldsymbol{r}^{\cair}\!\left(\boldsymbol{\cair};\,\boldsymbol{\theta}^{\cair},\boldsymbol{\kappa}^{\cair}\right) &= \boldsymbol{0},
\label{eq:airresidual}\\[2mm]
\boldsymbol{r}^{T}\!\left(\boldsymbol{T};\,\boldsymbol{\theta}^{ T },\boldsymbol{\kappa}^{ T }\right) &= \boldsymbol{0},
\label{eq:thermalresidual}
\end{align}
\end{subequations}
where \(\boldsymbol{r}^{\boldsymbol{\cair}}\) encodes the (mass) balance equation dictating the evolution of \cotwo, and $\boldsymbol{r}^{ T }$ encodes the (thermal energy) balance equation that describes the temperature evolution.

Given the information above, the forward model integrates the networks to predict $\boldsymbol{\cair}$ and $\boldsymbol{T}$ under airflow, occupancy, and thermal parameters. The remainder of this section specifies the concrete forms of the residuals $\boldsymbol{r}^{\boldsymbol{\cair}}$ and $\boldsymbol{r}^{ T }$ and the parameter sets $\boldsymbol{\theta}^{\cair}$, $\boldsymbol{\theta}^{T}$, $\boldsymbol{\kappa}^{\cair}$ and $\boldsymbol{\kappa}^{T}$. % and the coupling terms between them.

%%%=============================================
\subsection{The \cotwo network model}
\label{subsec:Airflow model}

%%%%%%%%%%%%%%----------------%%%%%%%%%%%%%%
% \begin{figure}[H]
% \centering
% \begin{tikzpicture}[
%     node/.style={draw, circle, minimum size=1.2cm},
%     edgelabel/.style={midway, fill=none, inner sep=1pt, sloped},
%     every edge/.style={draw, -{Stealth}, thick}
% ]
% % Nodes
% \node[node] (nodei)  at (0,0)   {$n_{i}$};
% \node[node] (nodei1) at (0,3.5) {$n_{i,1}$};
% \node[node] (nodei2) at (2.1,2.7) {$n_{i,2}$};
% \node[node] (nodeij) at (3.5,0) {$n_{i,j}$};
% % Edges with rotated labels (sloped)
% \draw[-{Stealth[scale=1.3]}, line width=0.8pt] (nodei)  -- node[edgelabel, above] {$\bair_{i,1}$, $\qair_{i,1}$} (nodei1);

% \draw[-{Stealth[scale=1.3]}, line width=0.8pt] (nodei) -- node[edgelabel, above] {$\bair_{i,2}$, $\qair_{i,2}$} (nodei2);

% \draw[-{Stealth[scale=1.3]}, line width=0.8pt] (nodei) -- node[edgelabel, above] {$\bair_{i,j}$, $\qair_{i,j}$} (nodeij);
% % Floating dots
% \node[rotate=-60] at (2.4,1) {$\cdots$};
% \end{tikzpicture}
% \caption{Schematic \cotwo network centered at node $n_{i}$.}
% \label{fig:airflow_network}
% \end{figure}
%%%%%%%%%%%%%%----------------%%%%%%%%%%%%%%
\begin{figure}[H]
\centering
\begin{tikzpicture}[
    node/.style={draw, circle, minimum size=1.cm},
    edgelabel/.style={midway, fill=none, inner sep=1pt, sloped},
    every edge/.style={draw, -{Stealth}, thick}
]
\def \elen{3}
% Nodes
\node[node] (node1)  at (0,0)   {$n_{1}$};
\node[node] (node2) at (\elen,0) {$n_{2}$};
\node[node] (node3) at (\elen,-\elen) {$n_{3}$};
\node[node] (node4) at (0,-\elen) {$n_{4}$};
\node[node] (node5) at (0,-2*\elen) {$n_{5}$};
% Edges with rotated labels (sloped)
\draw[-{Stealth[scale=1.3]}, line width=0.8pt] (node1)  -- node[edgelabel, above] {$e_1,\,\qair_1$} (node2);

\draw[-{Stealth[scale=1.3]}, line width=0.8pt] (node1) -- node[edgelabel, above] {$e_2,\,\qair_2$} (node3);

\draw[-{Stealth[scale=1.3]}, line width=0.8pt] (node3) -- node[edgelabel, above] {$e_3,\,\qair_3$} (node4);

\draw[-{Stealth[scale=1.3]}, line width=0.8pt] (node4) -- node[edgelabel, above] {$e_4,\,\qair_4$} (node1);

\draw[-{Stealth[scale=1.3]}, line width=0.8pt] (node3) -- node[edgelabel, above] {$e_5,\,\qair_5$} (node5);

\draw[-{Stealth[scale=1.3]}, line width=0.8pt] (node5) -- node[edgelabel, above] {$e_6,\,\qair_6$} (node4);

%\node[align=left,] at (\elen+4,0) {$\mathcal{N}_1=\{n_2,\,n_3,\,n_4\}\Rightarrow$ \\ $n_{1,1}=n_2, \; n_{1,2}=n_3, \; n_{1,3} = n_4$\\ $\mathcal{B}_1=\{e_1,\,e_2,\,e_4\} \Rightarrow$\\ $b_{1,1}=e_1, \; b_{1,2}=e_2, \; b_{1,3} = e_4$\\ $\mathcal{A}_i=\{1,\,1,\,-1\} \Rightarrow$\\ $\nu_{1,1}=1,\; \nu_{1,2}=1,\; \nu_{1,3}=-1$};
\node[align=left] at (\elen+5,-\elen) {$\mathcal{N}^q_1=\{n_2,\,n_3,\,n_4\}\Rightarrow\begin{cases}n_{1,1}^q=n_2\\n_{1,2}^q=n_3\\ n_{1,3}^q = n_4\end{cases}$\\ $\mathcal{E}^q_1=\{e_1,\,e_2,\,e_4\}\hspace{0.5em}\Rightarrow\begin{cases}e_{1,1}^q=e_1\\e_{1,2}^q=e_2\\ e_{1,3}^q = e_4\end{cases} $\\ $\mathcal{A}_1=\{1,\,1,\,-1\}\hspace{0.7em}\Rightarrow \begin{cases} a_{1,1}=1\\ a_{1,2}=1\\ a_{1,3}=-1 \end{cases}$\\ $\mathcal{Q}_1=\{q_1,\,q_2,\,q_4\}\hspace{0.7em}\Rightarrow \begin{cases} q_{1,1}=q_1\\ q_{1,2}=q_2\\ q_{1,3}=q_4 \end{cases}$};
\end{tikzpicture}
\caption{Schematic of local-to-global correspondence in the network. Here, $e_i$ denotes the $i$th oriented edge of the graph, where $q_i$ flows in the specified direction. %\icg{To be placed where fits!}
}
\label{fig:airflow_network}
\end{figure}
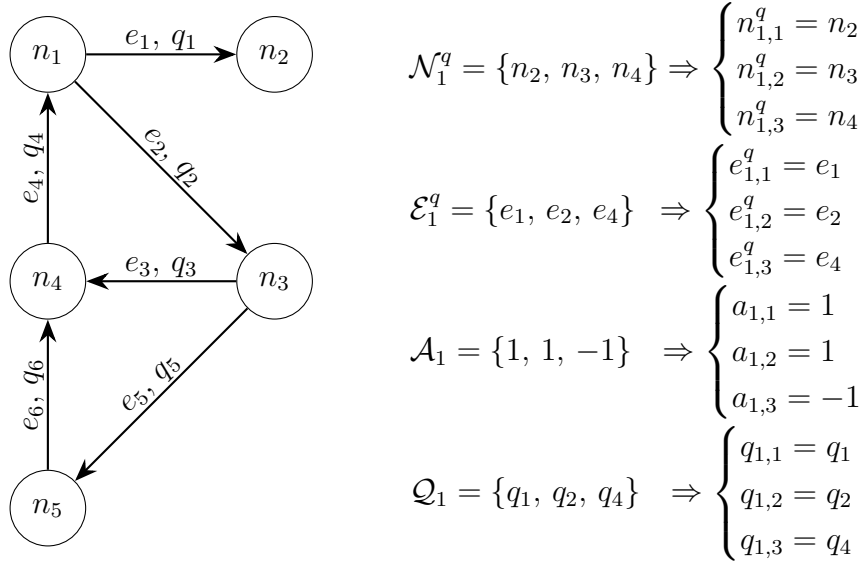
%%%%%%%%%%%%%%%%%--------%%%%%%%%%%%%%%%%%%%%%%%%

The \cotwo evolution is modeled as purely convective, and hence the network is described by an oriented graph whose nodes represent zones that store \cotwo mass and whose edges represent inter-zonal convective exchange. For each (globally counted) node \(n_i\), the set of neighboring nodes that are connected through convective transport are denoted as:
\[
\mathcal{N}^q_i=\{n^q_{i,1},n^q_{i,2},\ldots,n^q_{i,j},\ldots\},\qquad
\text{ and }\qquad
\mathcal{E}^{q}_i=\{\bair_{i,1},\bair_{i,2},\ldots,\bair_{i,j},\ldots\},
\]
where \(\bair_{i,j}\) links node \(n_i\) to node \(n^q_{i,j}\). See Figure~\ref{fig:airflow_network} for an illustration of these sets and the relation between global and local numbering of nodes and edges.

The edges in this graph are oriented; pointing from one node to the other. The volumetric air current along edge \(e_{i}\) is represented by the signed quantity \(\qair_{i}\), whose magnitude corresponds to the flow rate. For each node $n_i$, the orientations of the connected edges and the volumetric flow rates that live on these edges are encoded in the sets
\[
\mathcal{A}_i=\{a_{i,1},a_{i,2},\ldots,a_{i,j},\ldots\}\qquad
\text{ and }\qquad
\mathcal{Q}_i=\{q_{i,1},q_{i,2},\ldots,q_{i,j},\ldots\},
\]
where $a_{i,j} = 1$ means the edge points out of node $n_i$ and into node $n^q_{i,j}$, and $a_{i,j} = -1$ means the converse.

Assuming air incompressibility, the nodal mass balance reads
\begin{equation}
\sum_{j=1}^{|\mathcal{E}^q_i|} a_{i,j} \qair_{i,j}=0
\quad \forall\, i,
\label{eq:massb}
\end{equation}
which represents a constraint on the (to-be-inferred) parameters $q_i$.

\begin{remark}
The mass-balance equation is underdetermined, as the constraints provide less independent equations than the number of edge airflows. Therefore, a subset of edge flows can be chosen freely, while the remaining flows are determined by \eqref{eq:massb}. To chose the linearly independent quantities, we adopt a tree-cotree decomposition as done in electric network theory \cite{Nagel1975}. Specifically, we select a spanning tree of the airflow graph and partition the edge flows into (i) independent flows on cotree edges and (ii) dependent flows on tree edges. Given values of the independent flows, the dependent ones are obtained by solving the mass-balance equations, yielding a unique, constraint-consistent airflow field. In the remainder of the paper, the vector of airflow inference parameters $\boldsymbol{\qair}$ refers only to these independent (cotree) flows.

\end{remark}

In our \cotwo transport model, we exclusively consider convective transport between nodes, as this is expected to be the dominant driver across open pathways. We adopt an upwind scheme, whereby the \cotwo balance in zone \(i\) of volume \(V_i\) is described as:
\begin{equation}
\begin{aligned}
V_i\,\frac{d}{dt}\cair_i
&=\sum_{j=1}^{|\mathcal{E}^q_i|} \Big(
\underbrace{\langle - a_{i,j} \qair_{i,j}\rangle\,\cair_{n_{i,j}}}_{\text{convective inflow}}
- \!\!\!
\underbrace{\langle a_{i,j}\qair_{i,j}\rangle\,\cair_i}_{\text{convective outflow}} \!\!
\Big)
+\underbrace{n_{\mathrm{ppl},i}\,q_{\mathrm{exh}}\,\cair_{\mathrm{exh}}}_{\text{occupancy source}}
\quad \forall \, i.
\end{aligned}
\label{eq:co2b}
\end{equation}
Here, $\cair_i$ and $\cair_{n_{i,j}}$ denote the \cotwo concentrations in node $n_i$ and its neighbor $n_{i,j}$, respectively, $n_{ppl,i}$ is the occupancy in zone $i$, $q_{\mathrm{exh}}$ is the exhalation flow rate per-person, $\cair_{\mathrm{exh}}$ is the \cotwo volume fraction of exhaled air, and % To distinguish the \cotwo contribution associated with each airflow direction, the inflow and outflow terms must be separated using the
$\langle \bullet \rangle$ denote the Macaulay brackets:
\begin{equation}
\langle x\rangle = 
\begin{cases}
x, & \text{if } x \geq 0,\\
0, & \text{if } x < 0.
\end{cases}\label{eq:Rx}
\end{equation}

Under the assumption of the nodal mass balance \eqref{eq:massb}, the \cotwo balance \eqref{eq:co2b} over all nodes defines the residual operator \eqref{eq:airresidual}, i.e.,
\(
\boldsymbol{r}^{\boldsymbol{\cair}}\!\big(\boldsymbol{\cair};\,
\boldsymbol{\theta}^{\boldsymbol{\cair}},\boldsymbol{\kappa}^{\boldsymbol{\cair}}\big),
\)
with
\begin{subequations}\label{eq:air-params}
\begin{align}
\boldsymbol{\theta}^{\boldsymbol{\cair}} &= \{\boldsymbol n_{\mathrm{ppl}},\,\boldsymbol{\qair}\}
\quad\text{(unknown/inference parameters)}, \label{eq:air-params-a}\\[2mm]
\boldsymbol{\kappa}^{\boldsymbol{\cair}} &= \{\,\{V_i\},\, \qair_{\mathrm{exh}},\, \cair_{\mathrm{exh}},\, \text{BC}^{\boldsymbol{\cair}}\,\}
\quad\text{(known inputs)}, \label{eq:air-params-b}
\end{align}
\end{subequations}
where $\text{BC}^{\boldsymbol{\cair}}$ denotes the known boundary/ambient conditions. 
Given an initial state $\boldsymbol{\cair}_0$ and the above sets of parameters% and constants $(\boldsymbol{\theta}^{\boldsymbol{\cair}},\boldsymbol{\kappa}^{\boldsymbol{\cair}})$
, the forward solution $\boldsymbol{\cair}$ follows from satisfying $\boldsymbol{r}^{\boldsymbol{\cair}}=\boldsymbol{0}$ for $t>0$. This defines the forward map
\[
\mathcal F^{\boldsymbol{\cair}}:\; \big(\boldsymbol{\theta}^{\boldsymbol{\cair}},\boldsymbol{\kappa}^{\boldsymbol{\cair}},\boldsymbol{\cair}_0\big)
\longmapsto \boldsymbol{\cair}.
\]

%%%==============================================
\subsection{The temperature network model}
% \label{subsec:Thermal network}
% \begin{figure}[H]
% \centering
% \begin{tikzpicture}[
%     node/.style={draw, circle, minimum size=1.2cm},
%     edgelabel/.style={midway, fill=none, inner sep=1pt, sloped},
%     every edge/.style={draw, thick} % no arrows for conduction
% ]
% % Nodes
% \node[node] (nodei)  at (0,0)   {$n_{i}$};
% \node[node] (nodei1) at (0,3.5) {$n_{i,1}$};
% \node[node] (nodei2) at (2.1,2.7) {$n_{i,2}$};
% \node[node] (nodeij) at (3.5,0) {$n_{i,j}$};
% % Edges without arrows (conductive)
% \draw (nodei)  -- node[edgelabel, above] {$\bth_{i,1}$} (nodei1);
% \draw (nodei2) -- node[edgelabel, above] {$\bth_{i,2}$} (nodei);
% \draw (nodeij) -- node[edgelabel, above] {$\bth_{i,j}$} (nodei);
% % Floating dots
% \node[rotate=-60] at (2.4,1) {$\cdots$};
% \end{tikzpicture}
% \caption{Schematic temperature network centered at node $n_{i}$ with sloped labels.}
% \label{fig:thermal_network}
% \end{figure}

The temperature evolution is not only governed by convection, but also by conduction. Hence, we introduce a second network that connects the same set of zones as the \cotwo network, $n_i$, but with undirected edges, $e^R_i$. For each node \(n_i\), the node-local sets of neighboring nodes and incident edges of this graph are denoted
\[
\mathcal{N}^R_i=\{n^R_{i,1},n^R_{i,2},\ldots,n^R_{i,j},\ldots\}\qquad
\text{ and }\qquad
\mathcal{E}^{R}_i=\{\bth_{i,1},\bth_{i,2},\ldots,\bth_{i,j},\ldots\},
\]
where $\bth_{i,j}$ links node $n_{i,j}$ to node $n_i$.

The lumped resistance to conductive transport through edge $e^R_i$ is described by the resistance parameter $R_{i}$. For each node $n_i$, the resistance parameters for all its connected edges may be collected in the set
\[
\mathcal{R}_i=\{R_{i,1},R_{i,2},\ldots,R_{i,j},\ldots\}\,.
\]

While the set of nodes of the convective and conductive networks coincide, their connections, i.e., the edges, do not. Some connections (such as those through walls) only permit conductive heat transfer. Therefore, all edges of the convective network are also part of the conductive network but not vice versa.
%This means a one-to-one mapping between $\mathcal{N}^{\boldsymbol{\cair}}_i$ and $\mathcal{N}^{\mathrm{T}}_i$ or between $\mathcal{B}^{\boldsymbol{\cair}}_i$ and $\mathcal{B}^{\mathrm{T}}_i$ is not required. 

To define the evolution equation of the temperatures in our thermal network, we assign a lumped heat capacitance $\capacitor_i>0$ to each interior node $n_i$. Nodes that represent Dirichlet boundary conditions, which are, among others, exterior nodes, can be handled by fixing  their temperature $T_i$ and omitting $\capacitor_i$. %Each edge of the temperature network $\bth_{i,j}$ is assigned an effective thermal resistance $R_{i,j}>0$.
%With these definitions, 
The nodal energy balance for zone $n_i$ then becomes
\begin{equation}
\begin{aligned}
\capacitor_{i}\,\frac{d}{dt} T_i
&=\sum_{j}^{|\mathcal{E}^R_i|}\underbrace{\frac{T_{n_{i,j}}-T_i}{R_{i,j}} }_{\text{conduction}}
+\sum_{j}^{|\mathcal{E}^q_i|}\Big(
\underbrace{\langle - a_{i,j} \qair_{i,j}\rangle\,T_{n_{i,j}}}_{\text{convective inflow}}
-\!\!\!
\underbrace{\langle  a_{i,j} \qair_{i,j}\rangle\,T_{i}}_{\text{convective outflow}}
\!\!\!\Big)\,c_{p,\mathrm{air}}\,\rho_\mathrm{air}+\!\!\!\! \underbrace{n_{\mathrm{ppl},i}\,Q_{\mathrm{ppl}}}_{\text{occupancy source}}\!\! ,
\end{aligned}
\label{eq:tlnode}
\end{equation}
where %$T_i$ and $T_{n_{i,j}}$ denote the temperatures at node $n_i$ and its neighbor $n_{i,j}$, respectively, 
$Q_\mathrm{ppl}$ is the sensible heat per person, and $c_{p,\mathrm{air}}$ and $\rho_{\mathrm{\cair}}$ are the specific heat and density of air. The convective transport terms again involve the signed edge flows $\qair_{i,j}$, and we again adopt an upwind scheme. %from the \cotwo network account for advective heat transfer, and the Macaulay brackets $\langle \bullet\rangle$ ensure that the correct inflow and outflow contributions are considered in each term.

Collecting the nodal energy balances \eqref{eq:tlnode} over all zones defines the residual operator (\ref{eq:thermalresidual}), 
\(
\boldsymbol{r}^{T}\!\big(\boldsymbol{T};\,
\boldsymbol{\theta}^{T},\boldsymbol{\kappa}^{T}\big),
\)
where
\begin{subequations}\label{eq:thetakappa_tl}
\begin{align}
\boldsymbol{\theta}^{T} &= \{\,\boldsymbol{R},\,\boldsymbol{\capacitor},\boldsymbol n_{\mathrm{ppl}},\,\boldsymbol{\qair} \}
\quad\text{(unknown/inference parameters)}, \label{eq:thetatl}\\[2mm]
\boldsymbol{\kappa}^{T} &= \{\,c_{p,\mathrm{air}},\,\rho_\mathrm{air},\, Q_{\mathrm{ppl}},\,\text{BC}^{T}\,\}
\quad\text{(known inputs)}. \label{eq:kappatl}
\end{align}
\end{subequations}
Here $\boldsymbol{R}=\{R_{i}\}$ and $\boldsymbol{\capacitor}=\{\capacitor_i\}$ are the thermal resistances and capacitances, respectively, and 
$\text{BC}^{T}$ denotes known thermal boundary/ambient specifications (e.g., outdoor temperature). %The coupling to the airflow side enters through $\boldsymbol{\theta}^{\boldsymbol{\cair}}=\{\boldsymbol n_{\mathrm{ppl}},\boldsymbol{\qair}\}$ via sensible heat gains $\,\boldsymbol n_{\mathrm{ppl}}Q_{\mathrm{ppl}}$, and advective heat transport terms proportional to $\boldsymbol{\qair}$, $c_{p,\mathrm{air}}$, and $\rho_\mathrm{air}$.

Given an initial state $\boldsymbol{T}_0$ and parameter sets $(\boldsymbol{\theta}^{T},\boldsymbol{\kappa}^{T})$, the forward solution $\boldsymbol T$ follows from satisfying $\boldsymbol{r}^{T}=\boldsymbol{0}$ for $t>0$, defining the forward map
\[
\mathcal F^{ T }:\; \big(\boldsymbol{\theta}^{ T },\boldsymbol{\kappa}^{ T },\boldsymbol{T}_0\big)
\longmapsto \boldsymbol{T}.
\]

%%%%%%%%%%%%%%%%%%%%%%%%%%% bayesian part %%%%%%%%%%%%%%%%
\section{A Bayesian inference framework}
\label{sec:Bayesian}
\begin{figure}[H]
    \centering
    \includegraphics[width=0.7\linewidth]{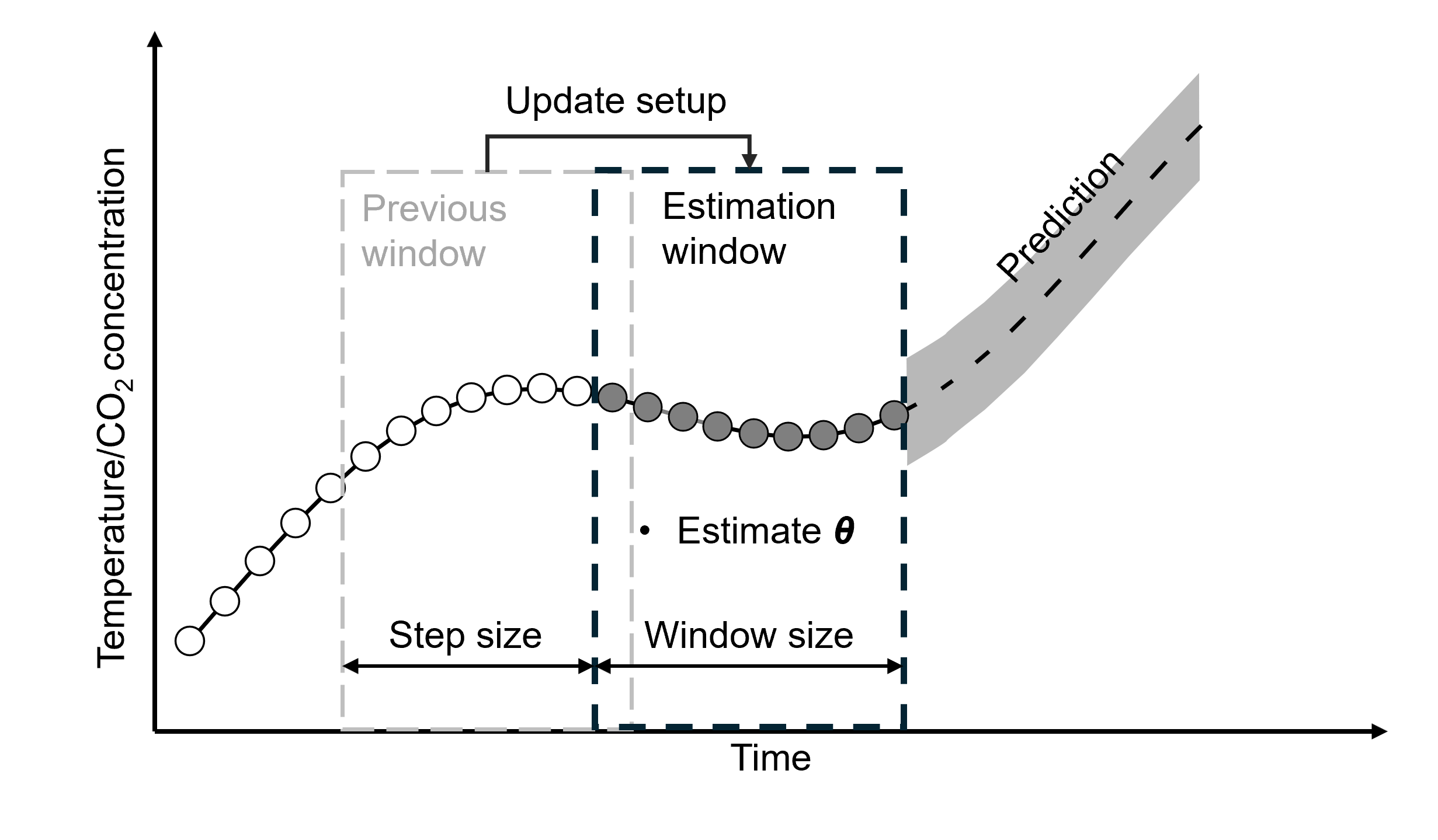}
    \caption{Schematic of the moving-window Bayesian framework.}
    \label{fig:bayesian_framework}
\end{figure}

In this section, we propose a moving-window Bayesian inference framework to enable data-driven, model-informed prediction under time-varying operating conditions. Rather than conditioning on the full measurement history, the inference is performed over a finite recent window so that parameter estimates remain representative of the current regime, while keeping the computational burden bounded for real-time application \cite{Valluru2017MovingWindow}. As illustrated in Figure~\ref{fig:bayesian_framework}, we use data $\mathcal{D} = \mathcal{D}^{\cair} \cup \mathcal{D}^{T}$ contained in the current inference window to infer the parameter vector $\boldsymbol{\theta} = \boldsymbol{\theta}^{\cair} \cup \boldsymbol{\theta}^{T}$ of the coupled \cotwo-temperature model.

The inference step produces not only a single parameter estimate, but a range of plausible parameter values consistent with the noisy measurements.
We use the resulting parameter samples to predict how \cotwo and temperature behave beyond the inference window, and report the variability as uncertainty bands around the prediction.
%We then propagate these parameter samples through the forward model to predict how \cotwo and temperature may behave, and report the resulting variability as uncertainty bands around the prediction. 
The window is then advanced by a fixed step, and the parameters are updated using the %posterior 
results from the previous window as initial guess.

To outline our inference framework, we commence with the introduction of the constituents of Bayes' rule: the prior, the likelihood, and the posterior. We then discuss the posterior evaluation through Markov-Chain Monte-Carlo sampling. In what follows, $p(\cdot)$ denotes a probability density function, and $p(a \mid b)$ denotes the conditional probability density of $a$ given that $b$ holds. 
%, i.e., the probability density distribution of $a$ after accounting for information $b$.

%=====================================================
\subsection{Bayes' rule}
\label{subsec:bayes_rule}
Eventually, Bayesian inference updates our initial beliefs about unknown model parameters using measurement data. The initial beliefs are encoded in the so-called ``prior'' distribution, whose probability density is denoted by $p(\boldsymbol{\theta})$. %, where $\boldsymbol{\theta}$ collects the parameters to be identified. 
The prior is a user-specified, assumed, distribution of the model parameters. Presumably, this probability distribution does not match reality, as expressed by the data $\mathcal{D}$. After observing the window data $\mathcal{D}$, we thus wish to obtain an updated probability distribution for $\boldsymbol{\theta}$, called the ``posterior'', denoted $p(\boldsymbol{\theta}\mid\mathcal{D})$. The posterior encodes the probability density of the model parameters $\boldsymbol{\theta}$ given the observed data $\mathcal{D}$, and thus quantifies which parameter values remain plausible after accounting for the data. 

Bayes' rule describes the relation between the prior and posterior given the data as 
\begin{equation}
p(\boldsymbol{\theta}\mid \mathcal{D})
= \frac{p(\mathcal{D}\mid \boldsymbol{\theta})\,
p(\boldsymbol{\theta})}{p(\mathcal{D})}.
\end{equation}
Here, $p(\mathcal{D}\mid \boldsymbol{\theta})$ denotes the ``likelihood'', i.e., the probability density distribution of the observed data given a set of model parameters. It quantifies how well the model predictions generated by $\boldsymbol{\theta}$ agree with the observed \cotwo and temperature data under an assumed measurement-noise model. The denominator, $p(\mathcal{D})$, is called the evidence (or marginal likelihood). In practice, it is a normalizing constant that ensures the posterior integrates to one and depends on the data and the chosen model/prior, but not on $\boldsymbol{\theta}$. 

The posterior distribution can be sampled without the need to evaluate the evidence, allowing us to work with the proportional form of Bayes' rule~\citep{gelman2013bayesian}:
\begin{equation}\label{eq:BayesProp}
p(\boldsymbol{\theta}\mid \mathcal{D})
\;\propto\; p(\mathcal{D}\mid \boldsymbol{\theta})\,p(\boldsymbol{\theta}).
\end{equation}
This form highlights the practical interpretation of Bayes' rule: the posterior favors parameter values that are both physically plausible under the prior and produce model outputs that match the measurements well.

The likelihood $p(\mathcal{D}\mid \boldsymbol{\theta})$ in equation \eqref{eq:BayesProp} is constructed from the windowed measurements $\mathcal{D}$ as follows.
Let $U^m(\boldsymbol{\theta}^m)\in\mathbb{R}^{J\times N}$ denote the solutions of our forward models \eqref{eq:airresidual} and  \eqref{eq:thermalresidual} with parameter values $\boldsymbol{\theta}^m$, for $m\in\{\cair,T\}$, respectively, evaluated at $J$ nodes corresponding to sensor positions
within $N$ time instants. Let $\mathcal{D}^m\in\mathbb{R}^{J\times N}$ be the corresponding sensor measurements.
%Let $U^m(\boldsymbol{\theta})\in\mathbb{R}^{J\times N}$ denote the model response for $m\in\{\cair,T\}$ evaluated at $J$ sensor nodes
%and $N$ time instants, and let $\mathcal{D}^m\in\mathbb{R}^{J\times N}$ be the corresponding measurements.
We assume additive, independent Gaussian measurement noise with node-specific standard deviations that are constant over the window as
\[
\mathcal{D}^m_{j,k} \;=\; U^m_{j,k}(\boldsymbol{\theta}^m) \;+\; \varepsilon^m_{j,k},
\qquad
\varepsilon^m_{j,k}\sim\mathcal N\!\big(0,({\sigma^m_j})^2\big),
\]
with $\mathcal N(\mu,\sigma^2)$ a univariate normal distribution with mean $\mu$ and variance $\sigma^2$.
The parameters $\sigma^m_j$ therefore quantify the measurement-noise scale at node $j$ for $m\in\{\cair,T\}.$

The measurement noise is assumed independent across all sensor nodes and time instants. Under these independence assumptions, the joint likelihood factorizes as
\[
p(\mathcal{D} \mid \boldsymbol{\theta}^m)
=p(\mathcal{D}^{\cair} \mid \boldsymbol{\theta}^{\cair})\;
 p(\mathcal{D}^T \mid \boldsymbol{\theta}^T),
\]
where
\begin{equation}\label{eq:likelyhoodfactored}
p(\mathcal{D}^m \mid \boldsymbol{\theta}^m)
=\prod_{j=1}^{J}\prod_{k=1}^{N}
\frac{1}{\sqrt{2\pi}\,\sigma^m_j}\,
\exp\!\left(-\frac{(\mathcal{D}^m_{j,k}-U^m_{j,k}(\boldsymbol{\theta}^m))^2}{2\,({\sigma^m_j})^2}\right).
\end{equation}
%=====================================================
%=====================================================
%=====================================================
%=====================================================
\subsection{Posterior evaluation and Markov-chain Monte-Carlo sampling}
\label{subsec:Posterior evaluation}
To avoid numerical stability (overflow/underflow) problems associated with the likelihood function, we evaluate the posterior on the log scale \cite{Rinkens2023}. %With two conditionally independent measurement variables, given the windowed data $\mathcal{D}=\{\tilde U_{\cair},\tilde U_T\}$ and the Gaussian measurement model from Section \ref{subsec:bayes_rule}, the log-posterior is
Taking the logarithm of (the proportional variant of) Bayes' rule \eqref{eq:BayesProp} yields
\begin{equation}
\log p(\boldsymbol{\theta}\mid \mathcal{D}) \propto\;
\log p(\mathcal{D}^{\cair}\mid \boldsymbol{\theta}^{\cair})\;+\;
\log p(\mathcal{D}^{T} \mid \boldsymbol{\theta}^T)\;+\;
\log p(\boldsymbol{\theta}),
\end{equation}
and the log-likelihood follows from \eqref{eq:likelyhoodfactored} as
\[
\log p(\mathcal{D}^m \mid \boldsymbol{\theta})
=-\frac12\sum_{j=1}^{J}\sum_{k=1}^{N}
\left\{\frac{(\mathcal{D}^m_{j,k}-U^m_{j,k}(\boldsymbol{\theta}^m))^2}{{\sigma^m_j}^2}
+\log\!\big(2\pi({\sigma^m_j})^2\big)\right\},
\quad m\in\{\cair,T\}.
\]
In this work, the measurement-noises $\sigma^m_j$ are treated as unknown parameters and inferred jointly with the physical parameters $\boldsymbol{\theta}$. This leads to one standard deviation per node $j$ and each $m\in\{\cair,T\}$ to be inferred. This choice considers the variety of sensor quality across rooms and measurement uncertainty into the posterior predictive bands. Hence, the Gaussian normalization term depends on $\sigma^m_j$. Including $\log(2\pi{(\sigma^m_j})^2)$ ensures a properly normalized Gaussian likelihood.

\begin{figure}[!b]
  \centering
  \includegraphics[width=0.78\linewidth]{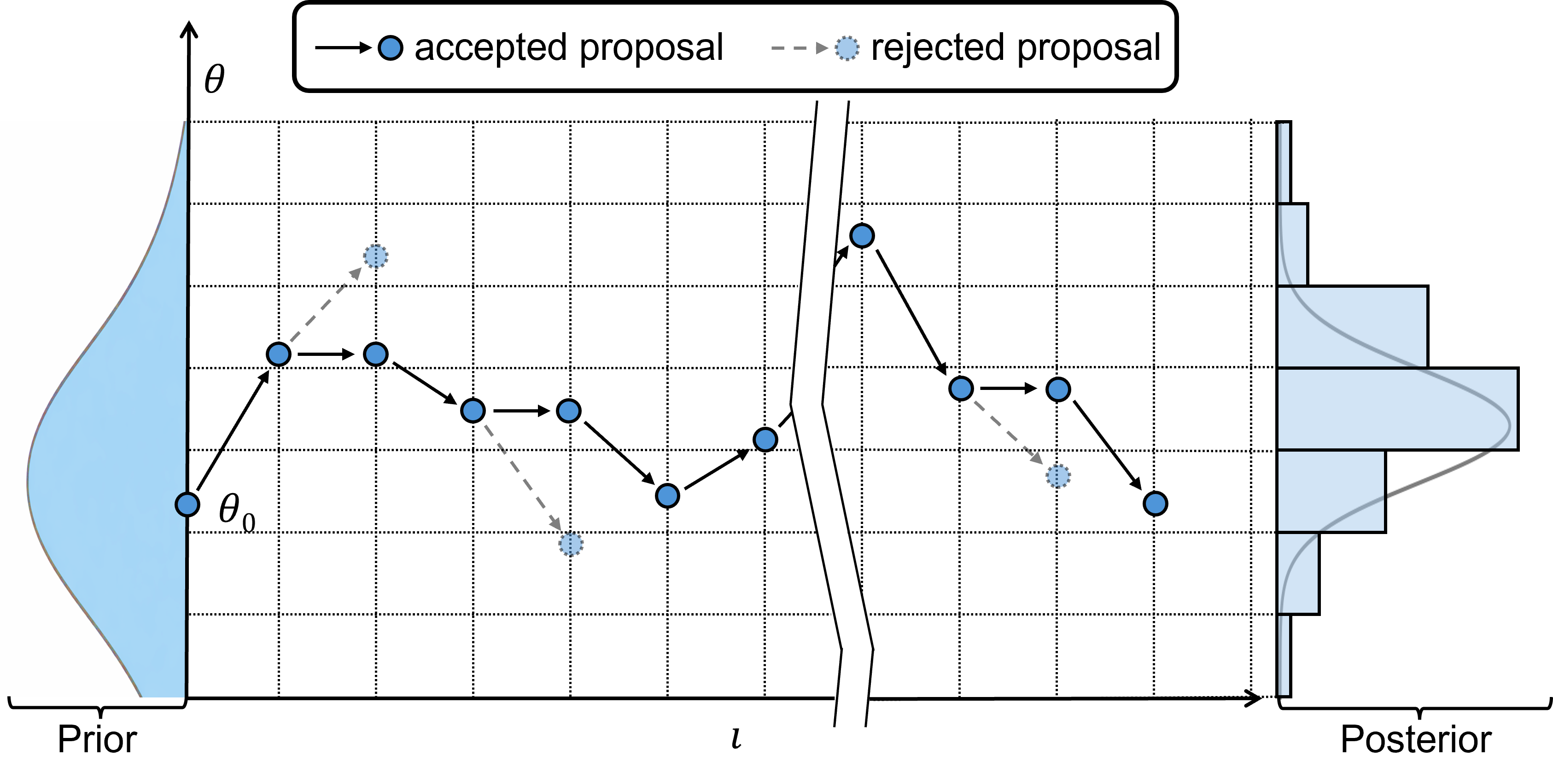}
  %\caption{\textbf{Schematic of a Metropolis–Hastings Markov-chain Monte-Carlo sampling from a one-dimensional posterior distribution.}}
  \caption{\textbf{Schematic of Metropolis--Hastings sampling in one dimension.} From the current state \(\theta_\imath\), a Gaussian random-walk proposal \(\theta^\star\) is accepted or rejected, producing samples whose histogram approximates the posterior; proposals outside the prior bounds are rejected.}
  \label{fig:MH}
\end{figure}

While the coupled \cotwo--temperature forward model is deterministic, the inverse problem is probabilistic because the measurements are noisy and some parameters might not be uniquely identifiable from the available data. We therefore approximate the posterior \(p(\boldsymbol{\theta}\mid\mathcal D)\) to quantify uncertainty and propagate it into the posterior predictive credible bands. We approximate this posterior using Markov-chain Monte-Carlo sampling, specifically the Robust Adaptive Metropolis (RAM) algorithm \cite{vihola2012robust}, an adaptive random-walk Metropolis--Hastings scheme that tunes both the global step size and the proposal covariance.

Figure~\ref{fig:MH} illustrates the Metropolis--Hastings sampling algorithm for a one-dimensional parameter in a Bayesian setting, where the starting point is sampled from the prior.
From the current state $\boldsymbol{\theta}_\imath$ (blue marker), a proposal is drawn using a Gaussian random-walk,
\[
\boldsymbol{\theta}^\star \sim \mathcal{N}\!\big(\boldsymbol{\theta}_\imath,\; s_\imath^{2}\Sigma_\imath\big),
\]
where $s_\imath>0$ is a scalar step-size factor at iteration $\imath$, and $\Sigma_\imath$ is the proposal covariance, which sets the direction-wise proposal variances and captures cross-parameter correlations in the Gaussian random walk. In RAM, $s_\imath$ and $\Sigma_\imath$ are not fixed user-tuned hyperparameters, but they are updated online from the evolving chain history to achieve efficient exploration, starting from an initial guess $(s_0,\Sigma_0)$. To decide whether a proposal is acceptable, the Metropolis--Hastings acceptance probability $\alpha$ is calculated as 
\[
\alpha(\boldsymbol{\theta}_\imath,\boldsymbol{\theta}^\star)
=
\min\!\left\{
1,
\frac{
p(\boldsymbol{\theta}^\star\mid\mathcal{D})
q(\boldsymbol{\theta}_\imath\mid\boldsymbol{\theta}^\star)
}{
p(\boldsymbol{\theta}_\imath\mid\mathcal{D})
q(\boldsymbol{\theta}^\star\mid\boldsymbol{\theta}_\imath)
}
\right\},
\]
where $\boldsymbol{\theta}^\star$ denotes the proposed parameter vector, $\boldsymbol{\theta}_\imath$ is the current parameter vector, and $q(\cdot\mid\cdot)$ denotes the proposal density. In the present case, the proposal is a Gaussian random walk centered at the current state. Conditional on the current proposal covariance $s_\imath^2\Sigma_\imath$, this proposal is symmetric, so that the proposal-density terms therefore cancel:
\[
q(\boldsymbol{\theta}^\star\mid\boldsymbol{\theta}_\imath)
=
q(\boldsymbol{\theta}_\imath\mid\boldsymbol{\theta}^\star).
\]

In the schematic, proposals outside the prior support are rejected ($\alpha=0$). If accepted, the chain moves to the proposed state $\boldsymbol{\theta}^\star$,
otherwise it remains at $\boldsymbol{\theta}_\imath$. Over iterations \(\imath\), the accepted samples concentrate in high-posterior regions, yielding the posterior histogram shown on the right.

After each iteration, RAM adapts the proposal parameters $(s_\imath,\Sigma_\imath)$ based on the chain history to keep the acceptance rate near a target value (we use $\approx 0.23$ for multi-dimensional random walks) \cite{RobertsGelmanGilks1997}. The RAM updates become progressively smaller as $\imath$ increases, so after an initial adaptation phase the proposal is effectively fixed. The sampler then behaves like a standard Metropolis--Hastings Markov-chain with a time-invariant transition kernel \cite{AndrieuThoms2008}. We initialize $\Sigma_0$ as a diagonal matrix whose entries are based on empirical parameter estimations, obtained from preliminary simulations and physical considerations. The initial step-size $s_0$ is chosen such that the proposal standard deviation of each parameter is approximately $1\%$ of its initial values. Further details on the RAM algorithm are available in \cite{Vihola2012}.
%=====================================================

\section{Synthetic benchmark problem}
\label{sec:Synthetic benchmark problem}

We first explore the performance of our modeling strategy in the context of a synthetic benchmark problem. This enables us to (i) validate the coupled \cotwo--temperature forward model against a controlled ground truth, (ii) examine how well occupancy- and airflow-driven effects can be identified from joint \cotwo and temperature measurements, and (iii) characterize how inference settings (e.g., window length, and noise level) influence reconstruction and predictive accuracy.

\subsection{Benchmark setup}
The benchmark we will consider is illustrated in Figure~\ref{fig:room_layout}. It represents a typical layout for commercial buildings, with four small offices (i.e., rooms A-D) and two larger meeting rooms (i.e., rooms E-F). All rooms have openings connected to the outside environment that allow air exchange. A centered hallway connects all the rooms. This hallway is divided into two zones (H1, H2) in which sensors are placed, allowing the model to capture the spatial inhomogeneity of the temperature and \cotwo states.

\subsubsection{Network structure}
\begin{figure}[H]
    \centering
    % (a)
    \subfloat[Room layout\label{fig:room_layout}]{
    \hspace*{0.04\textwidth} % left extra space
    \raisebox{0.55cm}{\includegraphics[height=0.213\textheight]{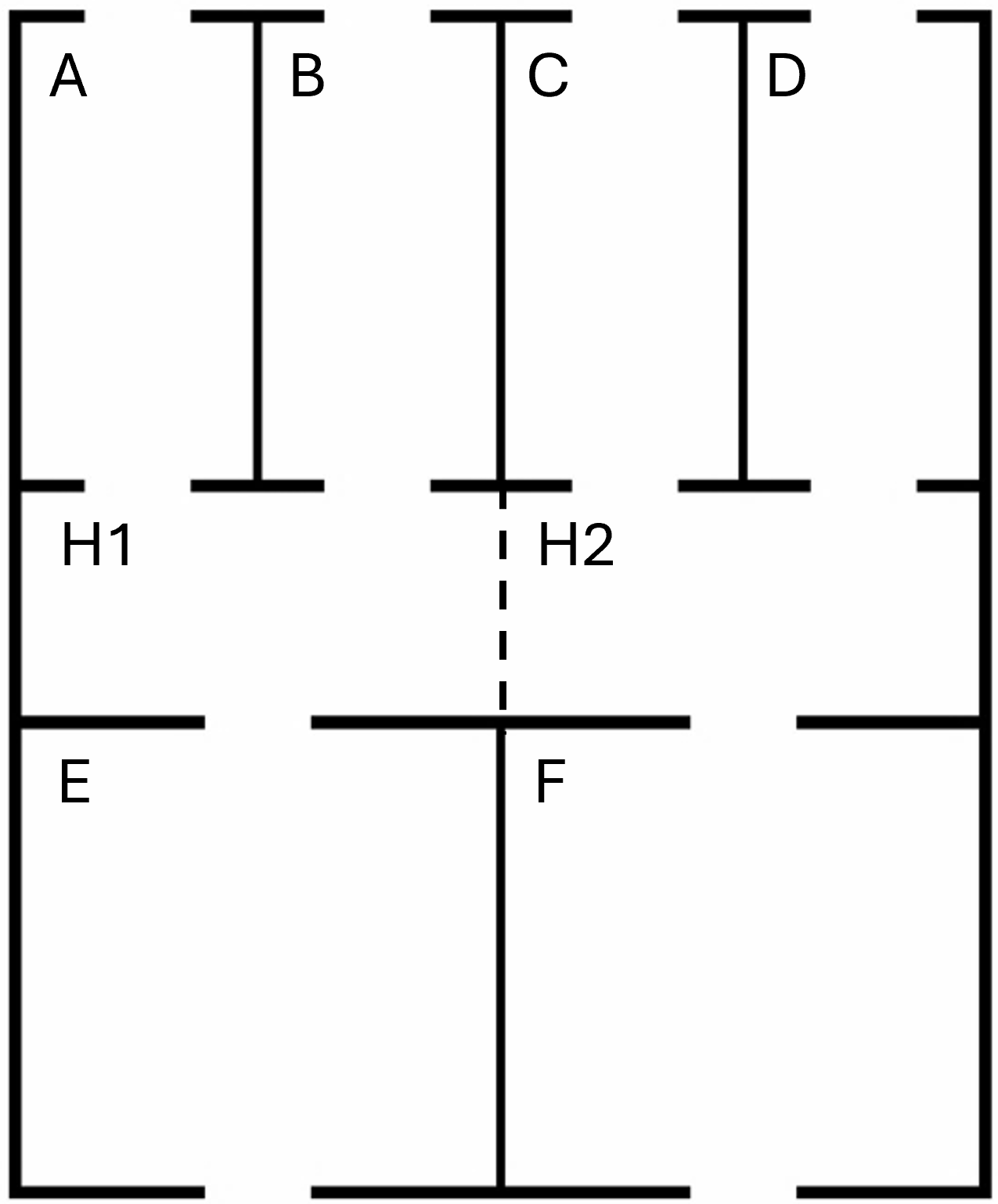}}
    \hspace*{0.04\textwidth} % right extra space
    }
    \hfill
    % (b)
    \subfloat[\cotwo network structure\label{fig:air_layout}]{
        \includegraphics[height=0.27\textheight]{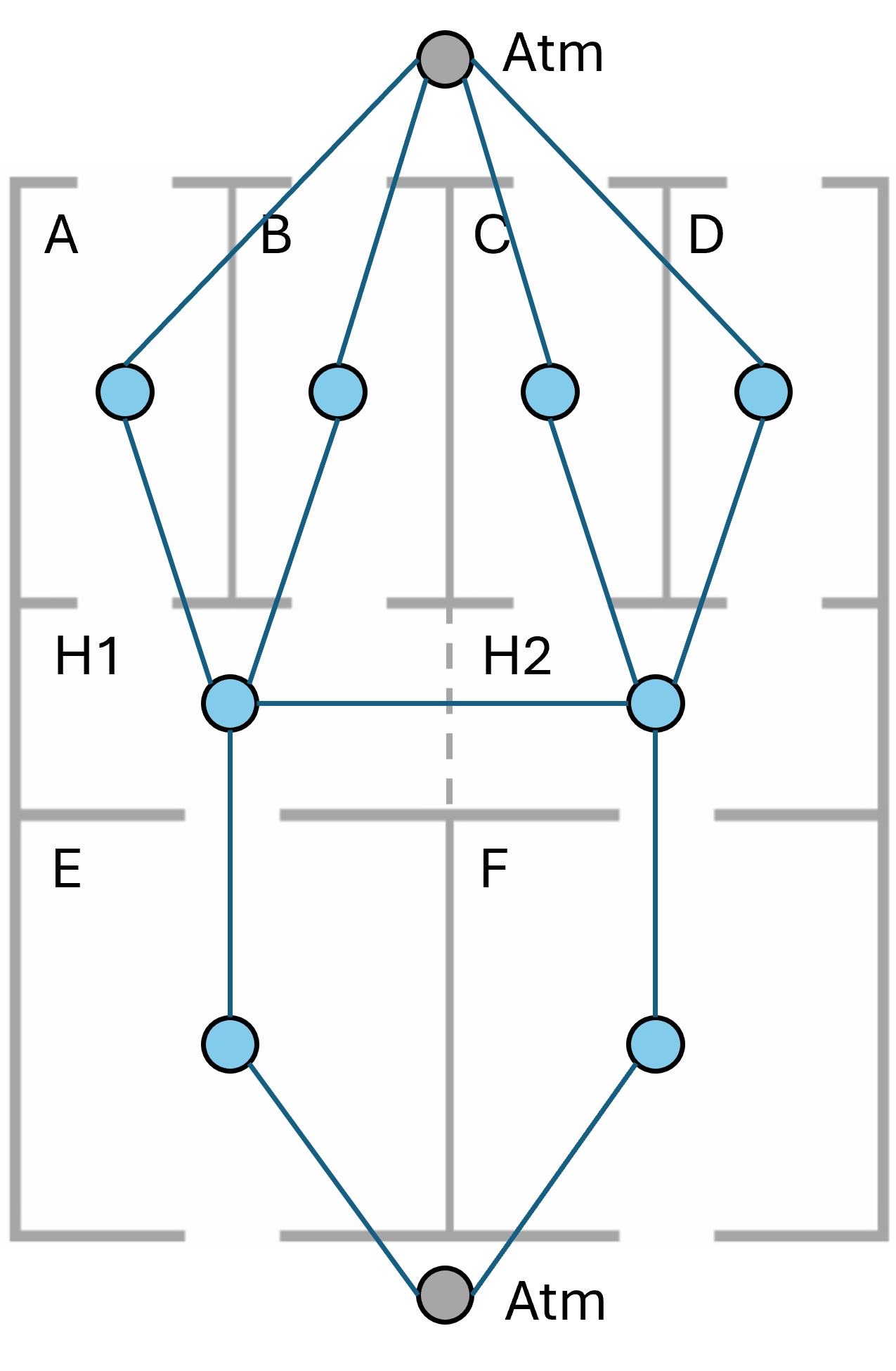}
    }
    \hfill
    % (c)
    \subfloat[Temperature network structure\label{fig:thermal_layout}]{
        \includegraphics[height=0.27\textheight]{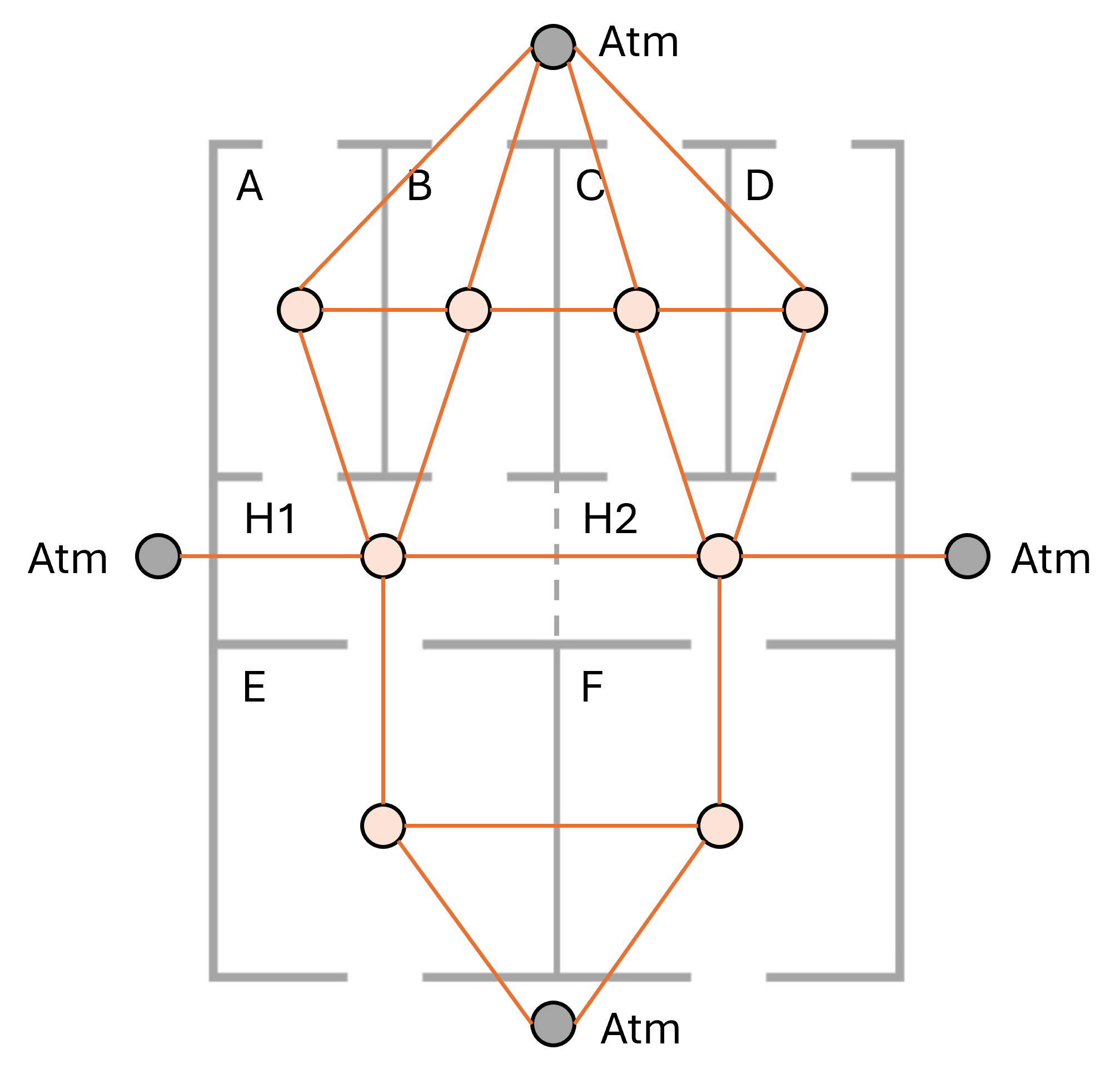}
    }
    
    \caption{Illustration of the benchmark problem room layout and network structure.}
    \label{fig:layouts}
\end{figure}

As discussed in Section~\ref{sec:Model}, the lumped temperature and \cotwo networks regard each room or zone as a node. Based on the layout of the building, Figures~\ref{fig:air_layout} and~\ref{fig:thermal_layout} represent the structure of the \cotwo and temperature networks. In both networks, the ``Atm'', atmosphere, nodes represent the same boundary conditions, which is the atmospheric temperature and the background \cotwo. %In the model, there is only one atmospheric node, which is shown multiple times in Figure \ref{fig:layouts} to ensure readability. 
Excluding the atmosphere nodes, the \cotwo and temperature networks share the same eight nodes but they differ in terms of connections (or edges). In the temperature network, all neighboring rooms are connected due to heat conduction through the walls, which are omitted in the \cotwo network since there is no active transport phenomenon between rooms that are not connected by open pathways (e.g., from room A to B).

\subsubsection{Parameter configuration and data generation}
\label{subsec:Parameter configuration and data generation}

Our benchmark problem considers two sets of parameters, parameters to be inferred ($\boldsymbol{\theta} = \boldsymbol{\theta}^{\cair} \cup \boldsymbol{\theta}^{ T }$), and constants or known parameters ($\boldsymbol{\kappa} = \boldsymbol{\kappa}^{\cair} \cup \boldsymbol{\kappa}^{ T }$). The parameters to be inferred and the constants are listed as:
\begin{subequations}
\begin{align}
    \boldsymbol{\theta} &= \{\boldsymbol n_{\mathrm{ppl}},~\boldsymbol{\qair},~\boldsymbol{R},~\boldsymbol{\capacitor}\}, \\
    \boldsymbol{\kappa} &= \{\text{BC}^{\cair},~\text{BC}^{T},~\{V_i\},~\qair_{\mathrm{exh}},~\cair_{\mathrm{exh}},~Q_{\mathrm{ppl}},~c_{p,\mathrm{air}},~\rho_{\mathrm{air}}\} \,.
\end{align}
\end{subequations}
Tables~\ref{tab:fixed parameter} and~\ref{tab:ground truth estimated} overview the parameter choices used to generate the synthetic data. Table~\ref{tab:fixed parameter} lists physical constants and boundary conditions that are treated as known in the Bayesian inference, whereas Table~\ref{tab:ground truth estimated} reports the ground-truth values of the parameters that we aim to infer. The thermal resistances and capacitances are chosen to reflect typical wall conduction and room-air thermal storage properties in office environments \cite{Mazzarella2020}.

Synthetic benchmark data are generated by forward simulation of the coupled \cotwo--temperature network model from Section~\ref{sec:Model}. We integrate the model over a 4--hour time window with a one-minute sampling interval. The occupancy schedule emulates an office-to-meeting transition: during the first two hours ($t<120$~min), four single-occupant offices (rooms A--D) are occupied while meeting rooms E and F remain empty; at $t=120$~min all occupants move to meeting room~F and the other rooms become unoccupied. %This transition is designed to induce \cotwo and temperature dynamics of sufficient magnitude to constrain the inference.

% Airflow boundary conditions are prescribed as a small set of independent façade exchange flows. This combination provides outdoor exchange while maintaining global consistency of the airflow field, thereby improving identifiability in the sensor signals.

The simulated \cotwo concentration and temperature trajectories constitute noiseless model outputs. To emulate realistic measurements, we add zero-mean Gaussian noise with standard deviations of $5$~ppm for \cotwo and $0.1^\circ$C for temperature, representative of low-cost office sensors. The resulting synthetic data set provides a controlled, yet realistic, test case for evaluating the performance and robustness of the Bayesian inference framework in Section~\ref{sec:Bayesian}.

\begin{table}[H]
\centering
\caption{Constants used for synthetic benchmark generation.}
\label{tab:fixed parameter}
\renewcommand{\arraystretch}{1.3} % Adjust spacing factor
\begin{tabular}{l p{0.65\textwidth}}
\hline\\[-10pt]
Ambient (outdoor) & $T_{\mathrm{amb}}=20~^\circ\mathrm{C}$, $\cair_{\mathrm{amb}}=400$ ppm \\[4pt]
Initial conditions & $T_{0,i}=20~^\circ\mathrm{C}$, $\cair_{0,i}=400$ ppm \\[4pt]
Room volume & $V_{\mathrm{A,B,C,D}}=30~\mathrm{m}^3$, $V_{\mathrm{E,F}}=60~\mathrm{m}^3$, $V_{\mathrm{H1,H2}}=15~\mathrm{m}^3$ \\[4pt]
\cotwo production per person & $q_{\mathrm{exh}}=1.0\times 10^{-5}~\mathrm{m}^3\!/\mathrm{s}$, \\&$\cair_\mathrm{exh}=50000~ppm$ (volume fraction) \\[4pt]
Heat production per person & $Q_{\mathrm{ppl}}=100~\mathrm{W}$ per person \\[4pt]
Air thermal properties 
& $c_{p,\mathrm{air}} = 1000~\mathrm{J}\!/\!(\mathrm{kg}\!\cdot\!\mathrm{K})$,\quad $\rho_\mathrm{air} = 1.2~\mathrm{kg}\!/\!\mathrm{m}^3$ \\[4pt]
\hline
\end{tabular}
\end{table} % TABLE CONSTANTS

\begin{table}[H]
\centering
\caption{Ground-truth values of the parameters inferred in the synthetic benchmark. The occupancy and boundary-airflow parameters are piecewise constant with a switch at $t_{\mathrm{sw}}=120$~min. Positive $q_{\mathrm{Atm},i}$ denotes inflow from outdoors into zone $i$.}
\label{tab:ground truth estimated}
\renewcommand{\arraystretch}{1.3}
\begin{tabular}{p{0.40\textwidth} p{0.55\textwidth}}
\hline\\[-10pt]
\multicolumn{2}{l}{\textbf{\cotwo network: occupancy and boundary airflow}} \\[4pt]
Occupancy & 0 to 120 min: 1 person in each room A, B, C, D,\\
          &120 to 240 min: 4 persons all move to room F. \\[6pt]
Independent boundary flows [$\mathrm{m}^3/\mathrm{s}$] &
$(q_{\mathrm{Atm},A},q_{\mathrm{Atm},B},q_{\mathrm{Atm},C},q_{\mathrm{Atm},D},q_{\mathrm{Atm},E})=$\\
&$(0.01,\,0.01,\,0.01,\,0.01,\,-0.01)\ $.\\[8pt]

\multicolumn{2}{l}{\textbf{Temperature network: thermal parameters}} \\[4pt]
Thermal resistances $\boldsymbol{R}$ {[K/W]} &
$R_{H1,H2}=1.0.$\\
& $R_{A,\mathrm{Atm}},R_{D,\mathrm{Atm}}=2,$\\
& $R_{B,\mathrm{Atm}},R_{C,\mathrm{Atm}}=R_{H1,\mathrm{Atm}}=R_{H2,\mathrm{Atm}}=3,$\\
& $R_{A,B}, R_{B,C},R_{C,D},R_{A,H1},R_{B,H1},R_{C,H2},$\\
& $R_{D,H2},R_{H1,E},R_{H2,F},R_{E,F},R_{E,\mathrm{Atm}},R_{F,\mathrm{Atm}}=1.5,$
\\[6pt]
Thermal capacitances $\boldsymbol{C}$ {[J/K]} &
$C_A,C_B,C_C,C_D=24000$; \ $C_E,C_F=48000$;\\
& $C_{H1},C_{H2}=12000$.\\[2pt]
\hline
\end{tabular}
\end{table}

\subsubsection{Inference configuration}
\label{subsec:Inference configuration}

The Bayesian inference framework is now adopted to infer the coupled \cotwo-temperature model parameter vector $\boldsymbol{\theta}$. To do so effectively, we introduce two more parameters to this set, namely the standard deviation of the sensor noise distributions, and the initial values for the \cotwo and temperature values within a time window. This leads to the following augmented set of parameters:
\[
\boldsymbol{\theta}^*=\Big\{
\underbrace{\boldsymbol{n}_{\text{ppl}}}_{\substack{\text{occupancy}\\\text{numbers (8)}}},\;
\underbrace{\boldsymbol{q}}_{\substack{\text{air}\\\text{flows (5)}}},\;
\underbrace{\boldsymbol{\cair}_0}_{\substack{\text{\cotwo}\\\text{initials (8)}}},\;
\underbrace{\boldsymbol{R}}_{\text{resistances (19)}},\;
\underbrace{\boldsymbol{C}}_{\text{capacitances (8)}},\;
\underbrace{\boldsymbol{T}_{0}}_{\substack{\text{temperature}\\\text{initials (8)}}},\;
\underbrace{\boldsymbol{\sigma}^{\cair},\quad \boldsymbol{\sigma}^T}_{\text{noise stds (8+8)}}
\Big\}.
\]

Table~\ref{tab:priors} lists the prior distributions chosen for these parameters. We use Gaussian priors $\mathcal{N}(\mu,\sigma^2)$ for unconstrained quantities and truncated Gaussian priors $\mathrm{TruncNorm}(\mu,\sigma^2;a,b)$ for parameters with physical bounds, where $\mu$ and $\sigma^2$ denote the mean and variance, and the support is restricted to $[a,b]$ with proper re-normalization. 

Most of the priors are truncated normals. The occupancy priors are truncated to enforce non-negativity and to reflect reasonable ranges given the respective room size. Boundary airflows are assigned zero-centered Gaussian priors. The initial states are specified per node: $\cair_{0,i}$ denotes the initial latent \cotwo concentration in zone $i$ at the start of a window, and we set $\cair_{0,i}\sim\mathcal{N}(\tilde{\cair}_i(t{=}0),\,\sigma_{\cair,i}^2)$, where $\tilde{\cair}_i(t{=}0)$ is the first \cotwo observation in that window and $\sigma^{\cair}_{i}$ is the corresponding node-wise measurement-noise standard deviation. Thermal resistances and capacitances are constrained to plausible building-physics ranges and grouped by room type. Finally, noise standard deviations are treated as hyperparameters, with priors centered on the synthetic noise levels ($\sigma_{\cair}=5$ ppm, $\sigma_T=0.1^\circ$C) that are also used to generate the benchmark data.The initial guesses adopted for the inference in each time window are set to reasonable physical values (occupancies $\approx 0.5$, boundary airflow $0.01~\mathrm{m}^3/\mathrm{s}$, $C$ around nominal capacitances, $T_0=20^\circ$C, $\cair_0=400$ ppm).

For each moving window, the Markov chain includes an initial burn-in phase, in which the samples are used to move the chain toward the high-probability region of the posterior distribution and to allow the adaptive proposal parameters to stabilize. These burn-in samples are discarded and are not used for posterior inference. In this study, we run $100{,}000$ iterations, discard the first $50{,}000$ as burn-in, and
retain the remaining $50{,}000$ draws for inference. Between windows we warm-start the sampler by setting the initial state
$\boldsymbol{\theta}_0$ of the next window to the posterior mean of the current window.

\begin{table}[H]
\centering
\caption{Prior distributions for all parameters that are to be inferred. % used for inference. TruncNorm$(\mu,\sigma,a,b)$ denotes a truncated normal with mean $\mu$, standard deviation $\sigma$, and support $[a,b]$.
}
\label{tab:priors}
\renewcommand{\arraystretch}{1.25}
\begin{tabular}{p{0.32\textwidth} p{0.6\textwidth}}
\hline\\[-10pt]
\textbf{Occupancies $n_{\mathrm{ppl},i}$} &
Rooms A--D: TruncNorm$(1,\,1;\,0,\,2)$; \\
&Rooms E--F: TruncNorm$(3,\,2;\,0,\,7)$; \\
&Hallway H1--H2: TruncNorm$(0,\,0.3;\,0,\,1)$. \\[7pt]

\textbf{Boundary airflows $q_{\mathrm{Atm},i}$ [$\mathrm{m}^3/\mathrm{s}$]} &
TruncNorm$(0,\,0.05;\,-0.1,\,0.1)$ for the five boundary nodes. \\[7pt]

\textbf{Initial \cotwo states $\cair_{0,i}$ [ppm]} & $\mathcal{N}\!\big(\tilde{\cair}_i(t{=}0),\,{\sigma^{\cair}_i}^2\big)$ per node (anchored to the first observation in each window). \\[7pt]

\textbf{Thermal resistances $R_{i,j}$ [K/W]} & TruncNorm$(1.0,\,0.5;\,0.5,\,2.0)$ for all 19 resistances. \\[7pt]

\textbf{Thermal capacitances $C_i$ [J/K]} &
Room A--D: TruncNorm$(20000,\,10000;\,10000,\,30000)$;\\
& Rooms E--F: TruncNorm$(40000,\,20000;\,20000,\,60000)$;\\
& Hallway H1--H2: TruncNorm$(10000,\,5000;\,5000,\,20000)$. \\[7pt]

\textbf{Initial temperatures $T_{0,i}$ [$^\circ$C]} & $\mathcal{N}\!\big(\tilde{T}_i(t{=}0),\,{\sigma^T_i}^2\big)$ per node (anchored to the first observation in each window). \\[7pt]

\textbf{CO$_2$ noise levels $\sigma_{\cair,i}$ [ppm]} &
$\sigma_{\cair,i} \sim \text{TruncNorm}(5.0,\,2.0;\,0,\,8.0)$, \quad $i\in\{\mathrm{A},\dots,\mathrm{H2}\}$. \\[7pt]

\textbf{Temperature noise levels $\sigma_{T,i}$ [$^\circ$C]} &
$\sigma_{T,i} \sim \text{TruncNorm}(0.1,\,0.1;\,0,\,0.5)$, \quad $i\in\{\mathrm{A},\dots,\mathrm{H2}\}$. \\[6pt]

\hline
\end{tabular}
\end{table}

\subsubsection{Error metrics}
\label{subsec:Error metrics and sweep over window size and noise}

To quantify the predictive performance of the moving-window Bayesian inference, we evaluate the discrepancy between the posterior predictive trajectories and the synthetic ground-truth data by the mean-normalized root-mean-square error (nRMSE), computed separately for \cotwo and temperature.

Let $\mathcal{D}^m \in \mathbb{R}^{8 \times N}$ denote the synthetic data points for $m \in \{\cair, T\}$ over time, and let $U^{m}(\boldsymbol{\theta})$ be the corresponding prediction generated from the posterior mean of the inferred parameters $\boldsymbol{\theta}$ for each inference window. % (propagated from the window start time up to $t = 240$~min). For a given verification interval in the prediction period $[t_{\text{start}}, t_{\text{end}}] \subset [0,240]$, 
Suppose the interval of verification contains $N_{\text{time}}$ time samples, indexed by $k = 1,\dots,N_{\text{time}}$. The root-mean-square error (RMSE) over this interval is then given as:
\begin{equation}
\mathrm{RMSE}
= \sqrt{\frac{1}{8\,N_{\text{time}}}
\sum_{j=1}^{8} \sum_{k=1}^{N_{\text{time}}}
\big(U^{m}_{j,k} - \mathcal{D}^m_{j,k}\big)^2},
\qquad m \in \{\cair, T\},
\end{equation}
where the factor $8\,N_{\text{time}}$ represents the total number of room–time pairs
(8 rooms and $N_{\text{time}}$ time samples).
The mean-normalized RMSE (nRMSE) now follows as
\begin{equation}
\mathrm{nRMSE}
= 100 \, \frac{\mathrm{RMSE}}{\displaystyle \frac{1}{8\,N_{\text{time}}}
\sum_{j=1}^{8} \sum_{k=1}^{N_{\text{time}}} \big|U^m_{j,k}\big|}.
\end{equation}

% In the synthetic benchmark, we consider the verification intervals to be $[80,120]$~min. Whenever the end of an inference window coincides with the beginning of a verification interval, the nRMSE is evaluated by comparing the corresponding prediction with the ground-truth benchmark over that verification interval.

%=====================================================
%=====================================================
\subsection{Synthetic benchmark results}
\label{sec:Synthetic benchmark problem results}
%This section evaluates the proposed moving-window Bayesian inference scheme on a fully controlled synthetic benchmark. 
Using the ground-truth parameters and prescribed occupancy/airflow drivers, we now assess (i) whether the proposed moving-window Bayesian inference scheme is able to reconstruct \cotwo and temperature trajectories within each window, (ii) how prediction accuracy and uncertainty evolve across the occupancy transition, and (iii) how inferred parameters and noise levels behave under different window sizes and measurement-noise settings.
\subsubsection{Synthetic experiment data}
\label{subsec:Synthetic experiment data}

\begin{figure}[H]
    \centering
    % (a)
    \subfloat[Synthetic \cotwo \label{fig:synthetic_co2}]{
        \includegraphics[width=1\linewidth]{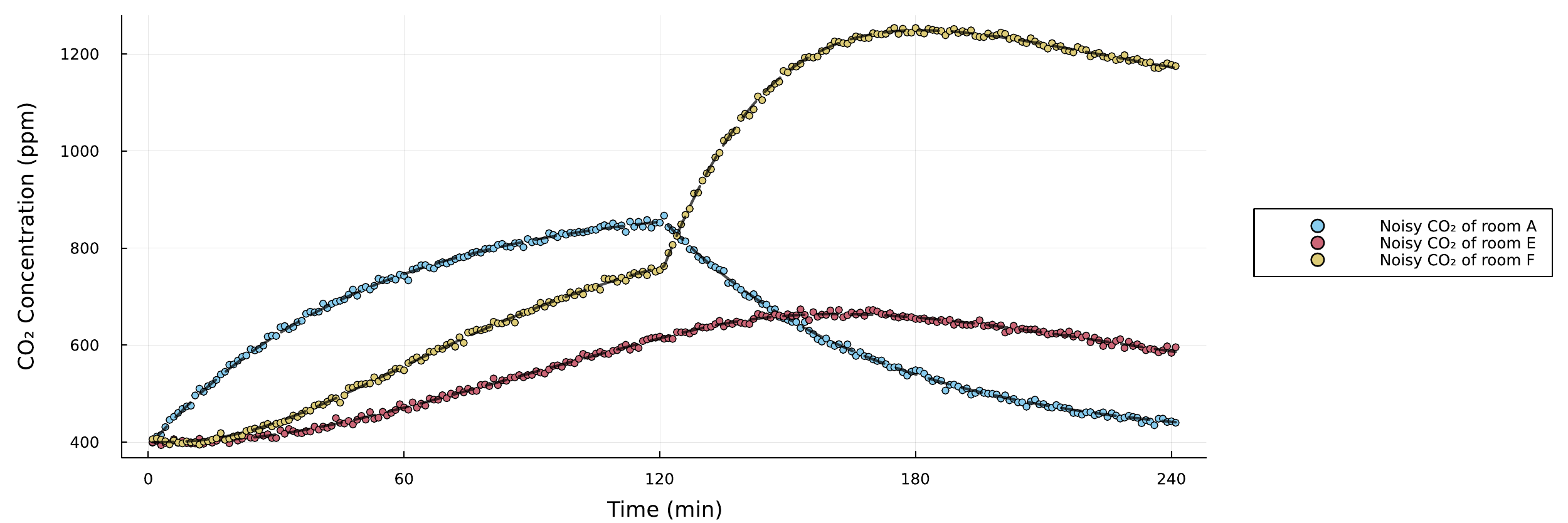}
    }\\[1em] % line break + small vertical gap

    % (b)
    \subfloat[Synthetic temperature\label{fig:synthetic_temp}]{
        \includegraphics[width=1\linewidth]{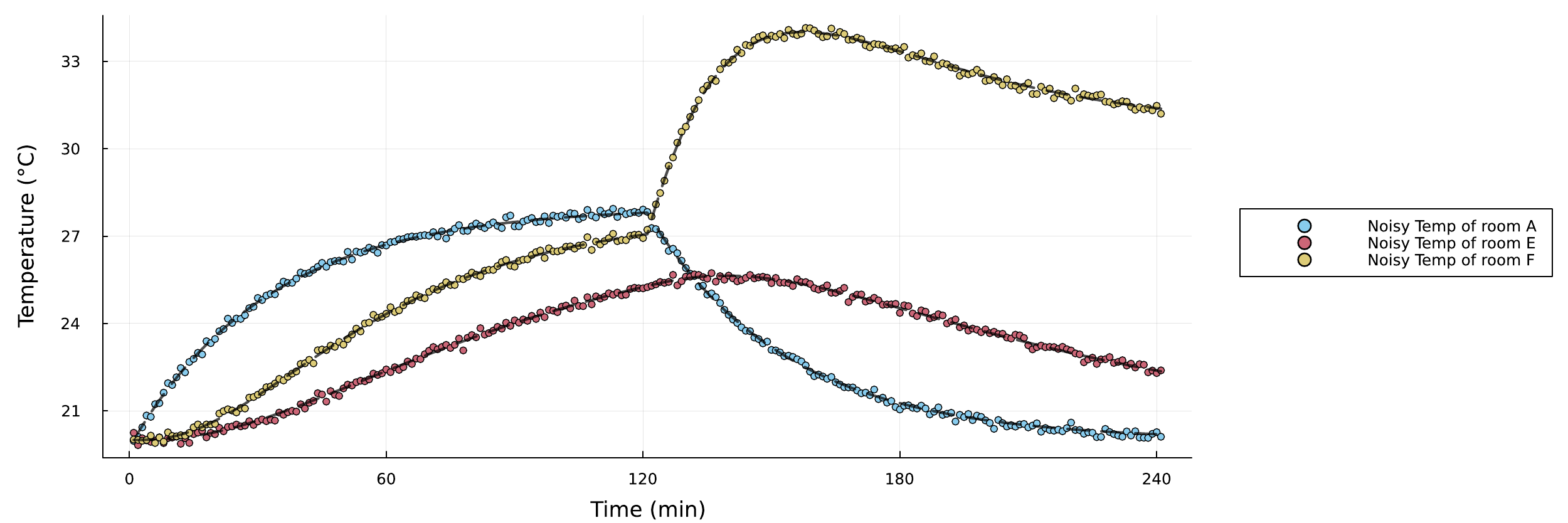}
    }

    \caption{Synthetic ground-truth benchmark results for rooms A, E, and F, with and without noise.}
    \label{fig:synthetic_results}
\end{figure}

Figure~\ref{fig:synthetic_results} shows the resulting \cotwo and temperature trajectories in three representative zones (offices A and E, and meeting room F). The dashed black curves correspond to the noise–free model output, while the colored markers represent synthetic measurements obtained by adding independent Gaussian noise with standard deviations of $5$ ppm for \cotwo and $0.3$~\textcelsius~for temperature. % to the simulation results of the ground-truth model. These standard deviations are consistent with low–cost indoor sensors. 
In room A, both \cotwo concentration and temperature initially rise from ambient levels as the single occupant generates \cotwo and heat, reaching peaks of approximately $850$ ppm and $28$~\textcelsius~before gradually decaying once the room is vacated. Room E, which remains unoccupied, exhibits only a modest increase driven by exchange with adjacent zones. In contrast, room F stays close to ambient conditions in the first half of the experiment and then shows a sharp increase in both \cotwo (up to about $1200$ ppm) and temperature (above $33$~\textcelsius) after the occupants move into that meeting room.

\subsubsection{Moving-window prediction and parameter tracking}
\label{subsec:Moving-window prediction and parameter tracking}

Figure~\ref{fig:moving-window-6} illustrates how the Bayesian framework operates with the moving-window scheme for rooms A, E and F. Each panel shows an inference window of 40 data points (indicated by the dashed black frame), together with the noisy synthetic measurements and the corresponding noiseless ground truth (colored markers and solid lines). The dashed curves denote the posterior predictive mean of the \cotwo concentration, computed from the mean of the inferred parameters in each window, while the gray ribbons indicate the associated 95\% credible bands for each room.

Within each window, the inferred trajectories closely follow the ground truth, whereas outside the window the short-term forecasts rely solely on the information available up to the right edge of the dashed box. Panels~(a) and~(b) correspond to windows located entirely before the occupancy transition at $t=120$~min. In this regime the sampler recovers the pre-transition occupancy and airflow pattern, which leads to accurate reconstructions of the rapidly rising \cotwo levels in the single-occupant offices and the much smaller increase in the meeting rooms, but, by construction, does not anticipate the occupancy change after $t=120$~min. When the window overlaps the transition point (panels~(c) and~(d)), the model is inherently unable to accurately predict the evolution behavior visible in the data. The strength of a Bayesian inference framework is that this unidentifiability of suitable model parameters is quantifiable, as reflected by
%the inference rapidly adapts: the posterior reallocates the occupants from the small offices to meeting room~F, so that the predicted \cotwo evolution closely matches the ground truth and 
the uncertainty bands temporarily widening around the change point. Moreover, the framework is able to adapt and stabilize after a mismatch. In the later window (panels~(e)), after the transition has been fully observed, the calibrated model again produces predictions with narrow uncertainty that remain close to the ground truth in the new regime.
%%%%%%%%%%%%%%%%%%%%%%%%%%%%%%%%%%%%%%%%%%%%%%%%%%%%%%%%%%%%%%%%%%

\begin{figure}[H]
  \centering
    % Row 1
    \subfloat[Window 11-51]{\includegraphics[width=0.5\textwidth,trim={0.7cm 1cm 11cm 0.3cm},clip] {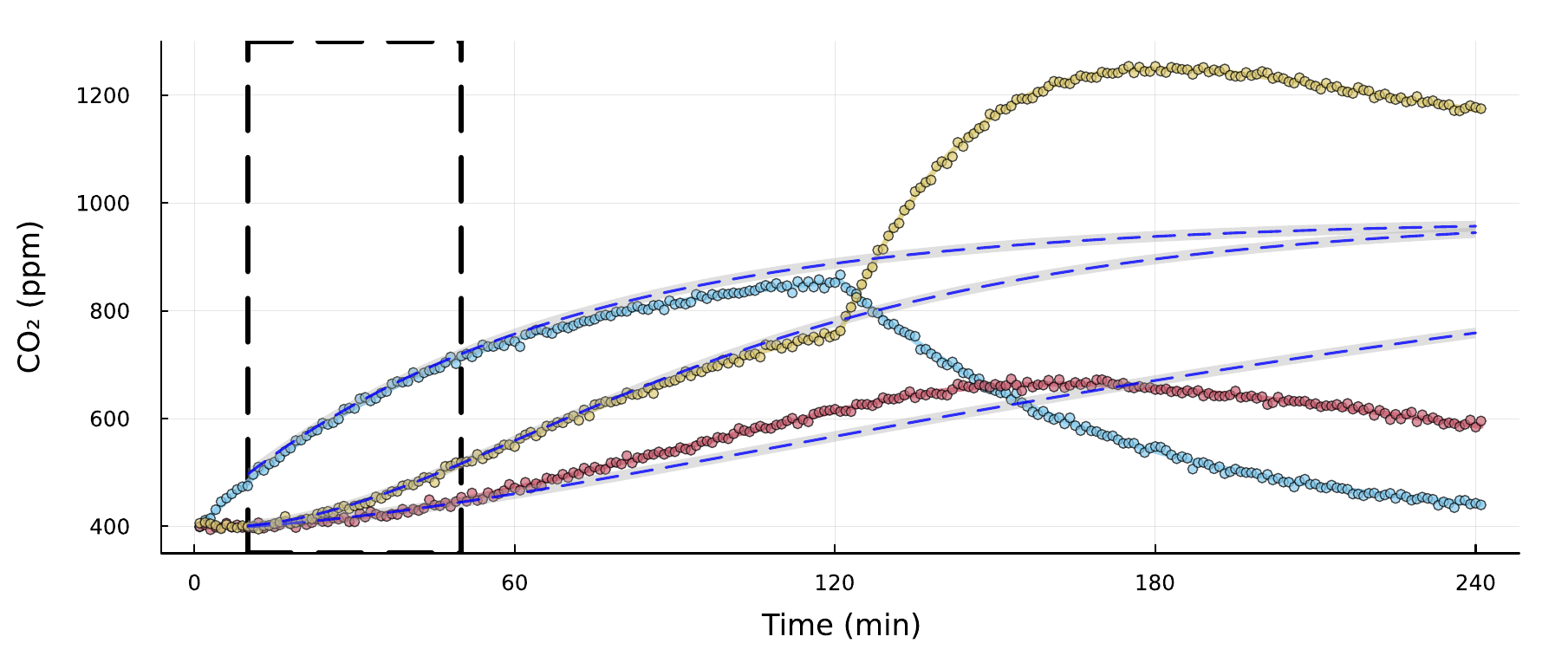}\label{fig:co2window11}}
    \subfloat[Window 51-91]{\includegraphics[width=0.5\textwidth,trim={0.7cm 1cm 11cm 0.3cm},clip]{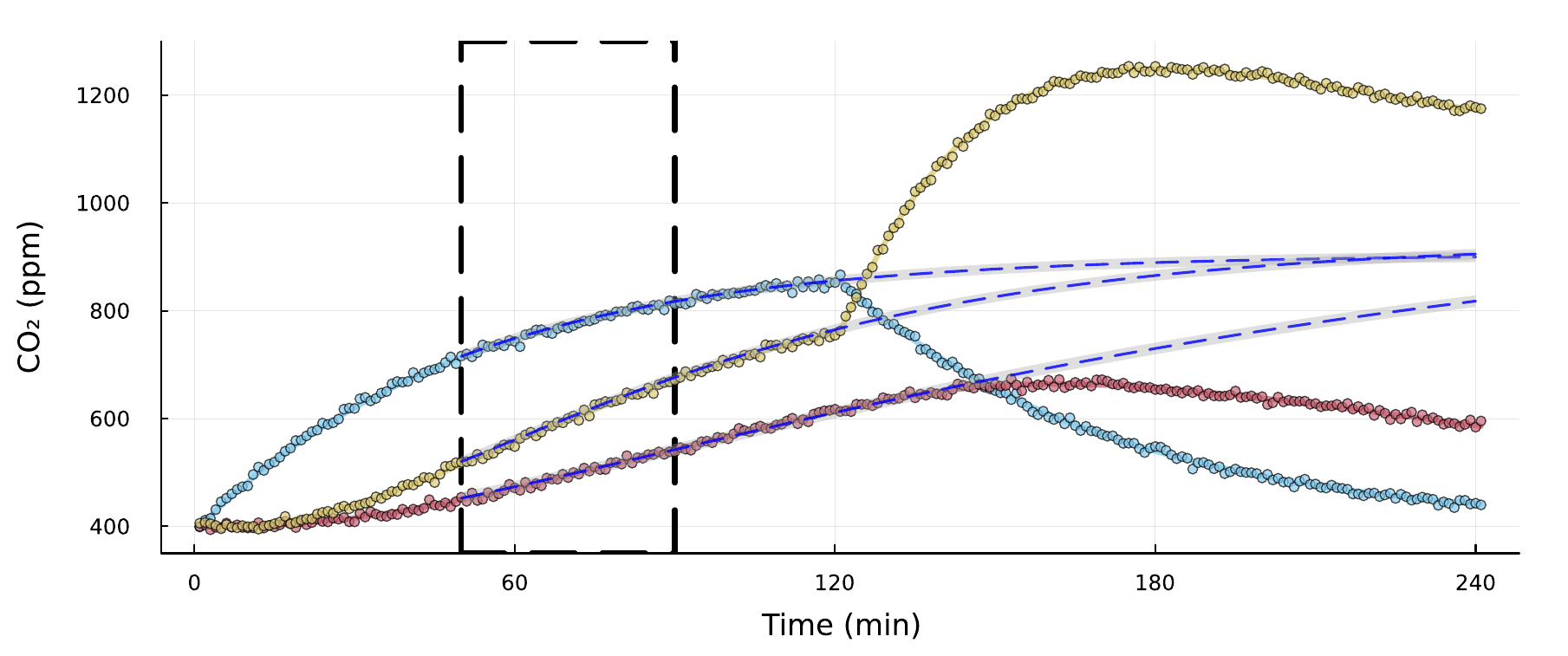}\label{fig:co2window51}}
    
    % Row 2
    \subfloat[Window 91-131]{\includegraphics[width=0.5\textwidth,trim={0.7cm 1cm 11cm 0.3cm},clip] {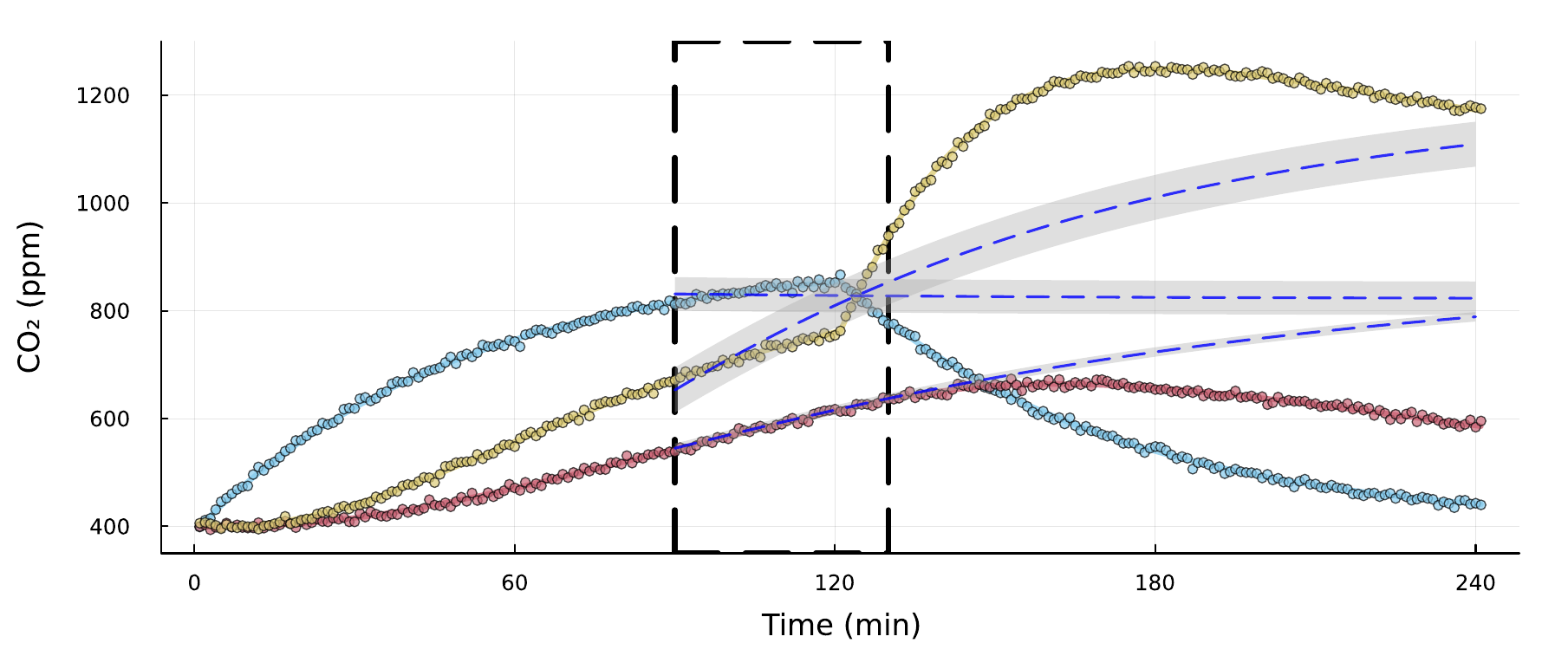}\label{fig:co2window91}}
    \subfloat[Window 101-141]{\includegraphics[width=0.5\textwidth,trim={0.7cm 1cm 11cm 0.3cm},clip]{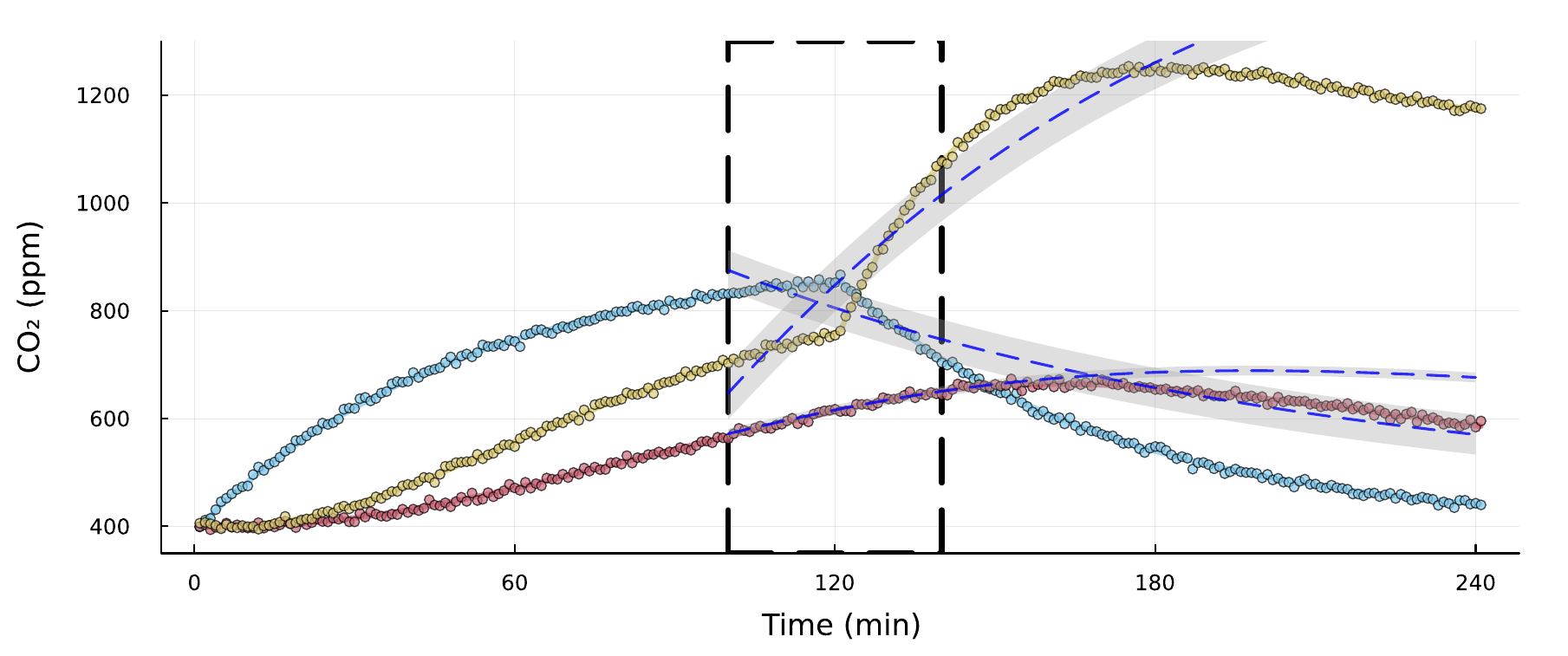}\label{fig:co2window101}}

    % row 3
    \subfloat[Window 121-161]{\includegraphics[width=0.5\textwidth,trim={0.7cm 1cm 11cm 0.3cm},clip]{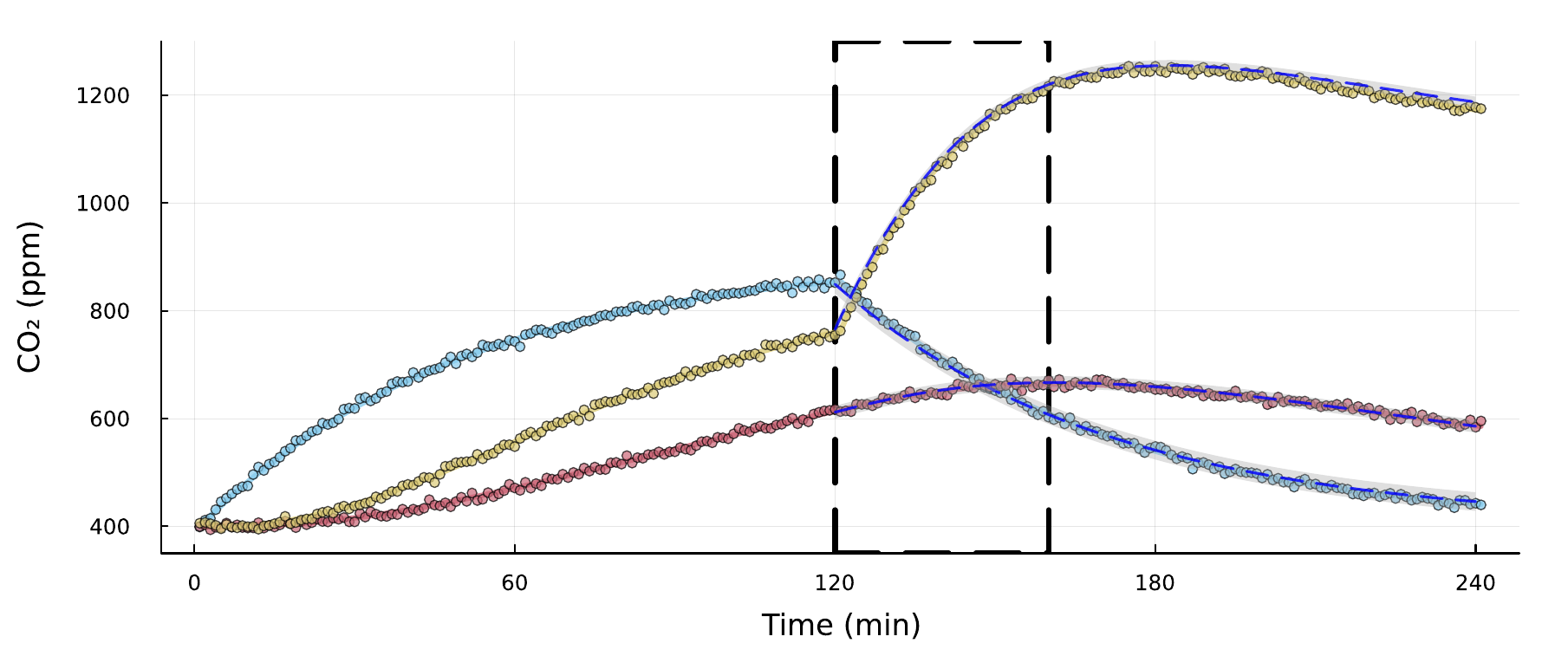}\label{fig:co2window121}}\hfill
    \subfloat{\includegraphics[width=0.35\textwidth,trim={30.5cm 5cm 2.5cm 4.2cm},clip]{Figures/Results/moving_window/121c.pdf}}\hspace{0.08\textwidth}
    
    \caption{Posterior predictive \cotwo trajectories in rooms A, E and F %obtained with the moving-window Bayesian scheme 
    for the synthetic-data experiment. %Panels (a)–(e) show five overlapping 40-min inference windows (dashed rectangles), with noisy synthetic measurements and ground truth, posterior predictive mean, and shared legend.
    }
  \label{fig:moving-window-6}
\end{figure}

\begin{figure}[H]
  \centering
    % Row 1
    \subfloat[Window 11-51]{\includegraphics[width=0.5\textwidth,trim={0.7cm 1cm 11cm 0.3cm},clip] {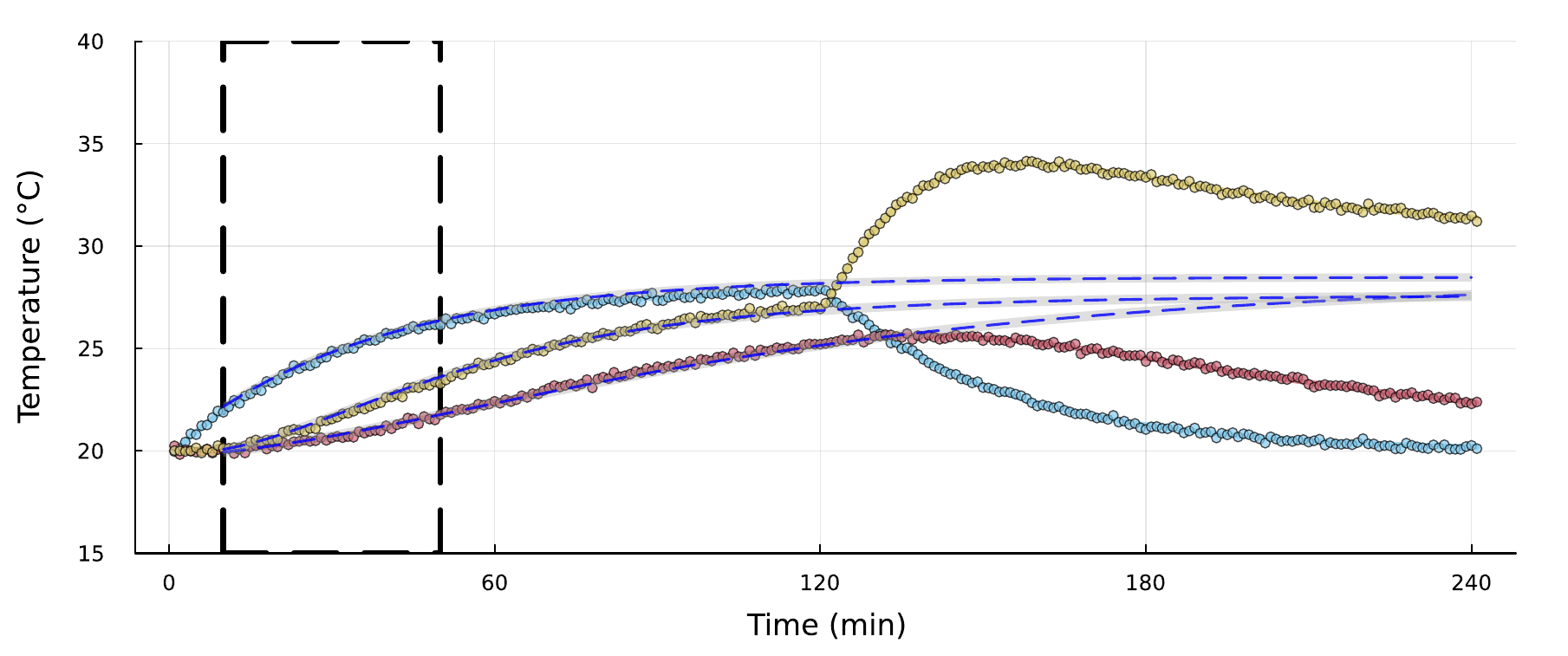}\label{fig:tempwindow11}}
    \subfloat[Window 51-91]{\includegraphics[width=0.5\textwidth,trim={0.7cm 1cm 11cm 0.3cm},clip]{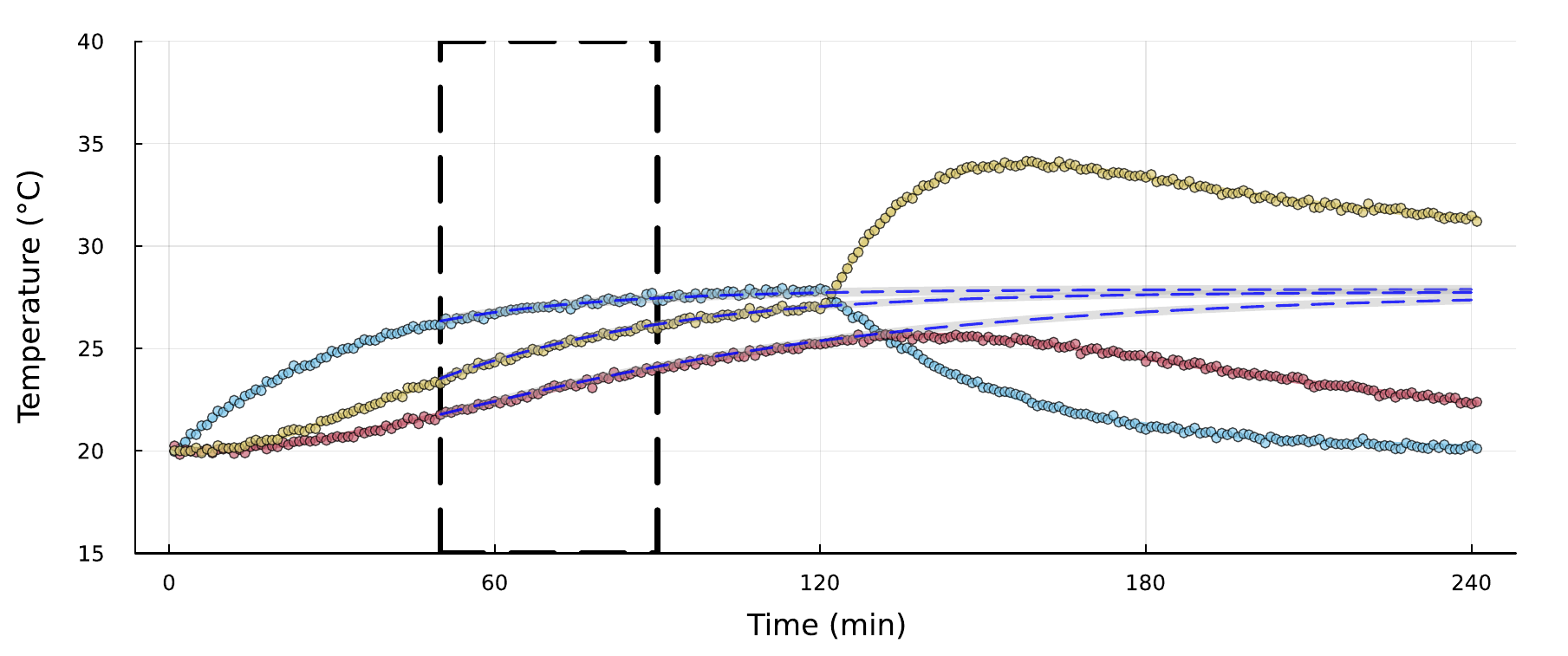}\label{fig:tempwindow51}}
    
    % Row 2
    \subfloat[Window 91-131]{\includegraphics[width=0.5\textwidth,trim={0.7cm 1cm 11cm 0.3cm},clip] {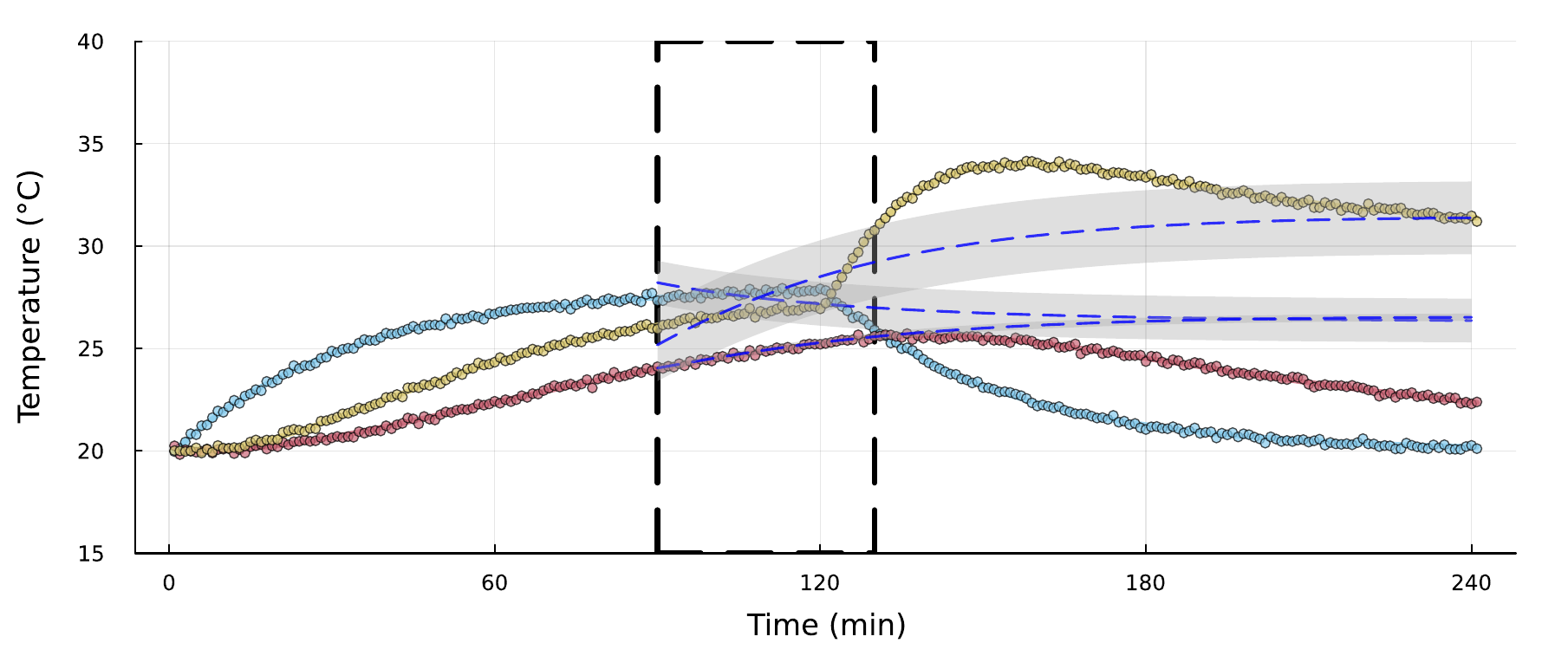}\label{fig:tempwindow91}}
    \subfloat[Window 101-141]{\includegraphics[width=0.5\textwidth,trim={0.7cm 1cm 11cm 0.3cm},clip]{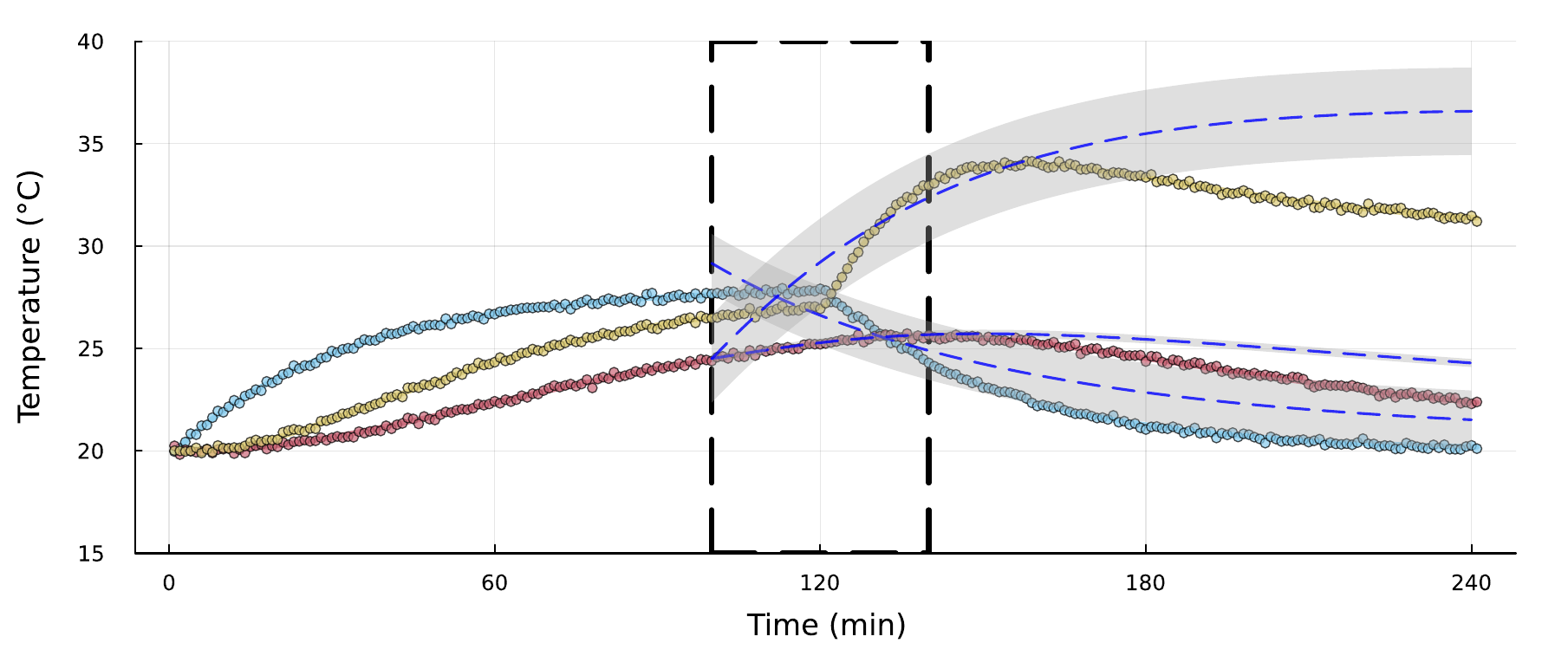}\label{fig:tempwindow101}}

    % row 3
    \subfloat[Window 121-161]{\includegraphics[width=0.5\textwidth,trim={0.7cm 1cm 11cm 0.3cm},clip]{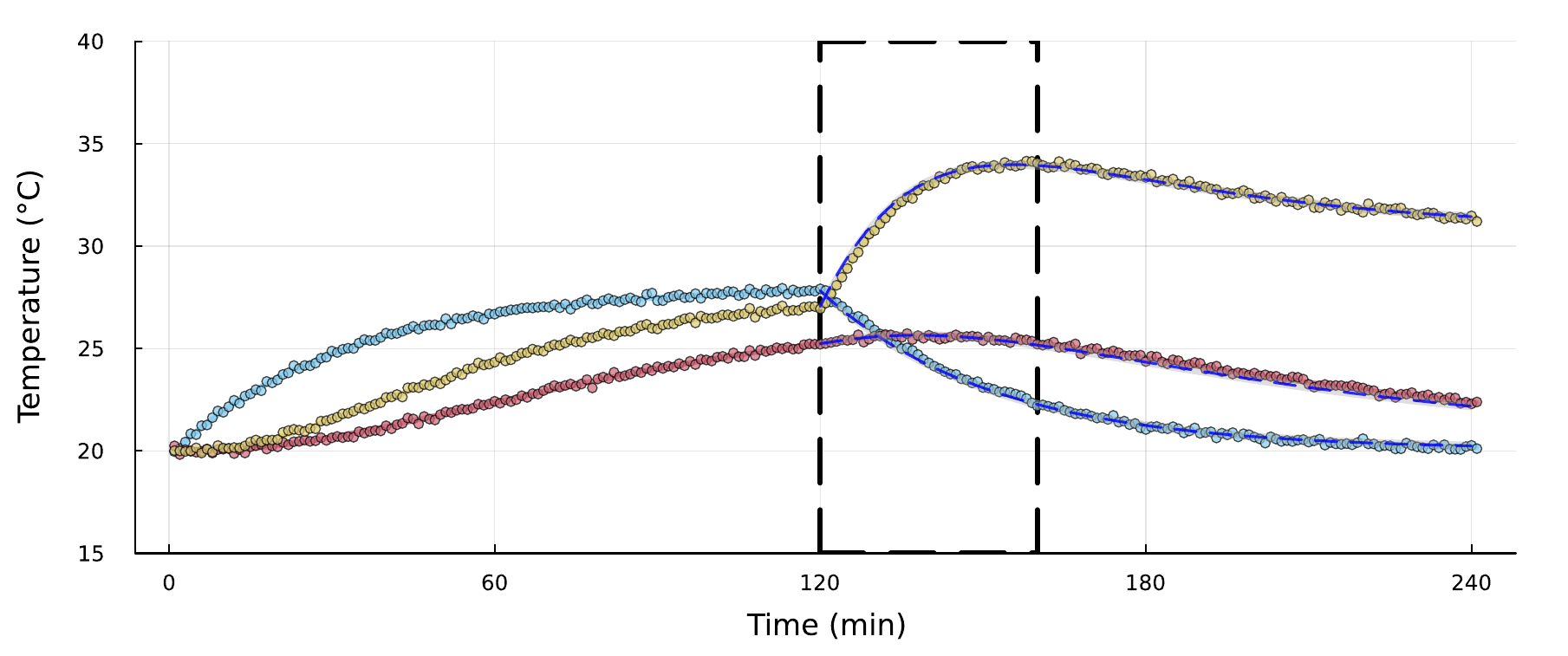}\label{fig:tempwindow121}}\hfill
    \subfloat{\includegraphics[width=0.35\textwidth,trim={30.5cm 5cm 2.5cm 4.2cm},clip]{Figures/Results/moving_window/121t.pdf}}\hspace{0.08\textwidth}
    
    \caption{Posterior predictive temperature trajectories in rooms A, E and F %obtained with the moving-window Bayesian scheme 
    for the synthetic-data experiment. %Panels (a)–(e) show five overlapping 40-min inference windows (dashed rectangles), with noisy synthetic measurements and ground truth, posterior predictive mean, and shared legend.
    }
  \label{fig:moving-window-6}
\end{figure}

% %%%%%%%%%%%%%%%%%%%adjusting figure%%%%%%%%%%%%%%%%%%

\begin{figure}[!t]
  \centering
    % Row 1
    \subfloat[Occupancy parameter tracking]{%
    \includegraphics[width=0.8\textwidth,trim={2cm 14cm 0.5cm 27.1cm},clip]
    {Figures/Results/track_parameter/all_params_4panels}
    \label{fig:occupancy_tracking}}
    
    \subfloat[Airflow parameter tracking]{%
    \includegraphics[width=0.8\textwidth,trim={2cm 1cm 0.5cm 40.5cm},clip]
    {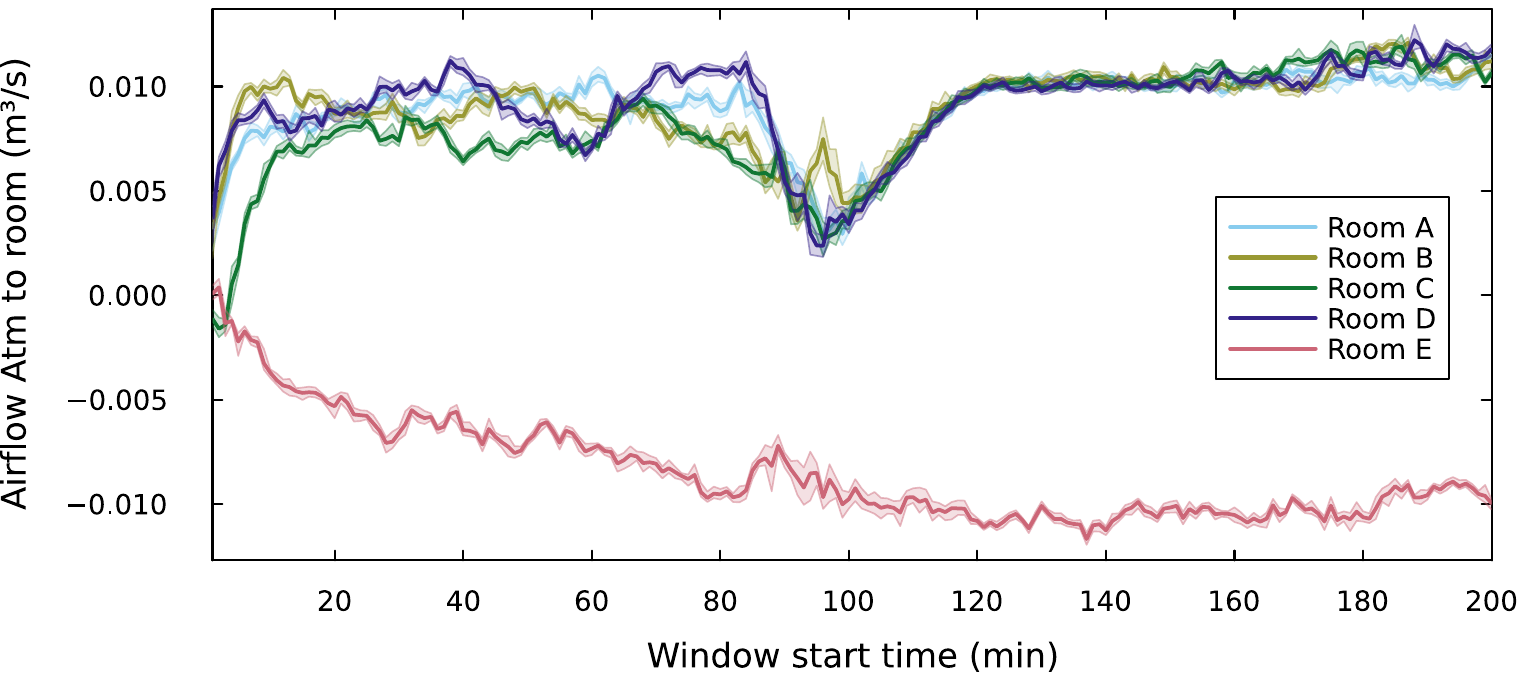}
    \label{fig:airflow_tracking}}
    
  \subfloat[Thermal capacitance tracking]{%
    \includegraphics[width=0.8\textwidth,trim={2cm 40.7cm 0.5cm 0.7cm},clip]
    {Figures/Results/track_parameter/all_params_4panels.pdf}%
    \label{fig:C_tracking}%
  }
\end{figure}

\begin{figure}[!t]
  \ContinuedFloat 
  \subfloat[Thermal resistance tracking]{%
    \includegraphics[width=0.8\textwidth,trim={2cm 27.5cm 0.5cm 14cm},clip]
    {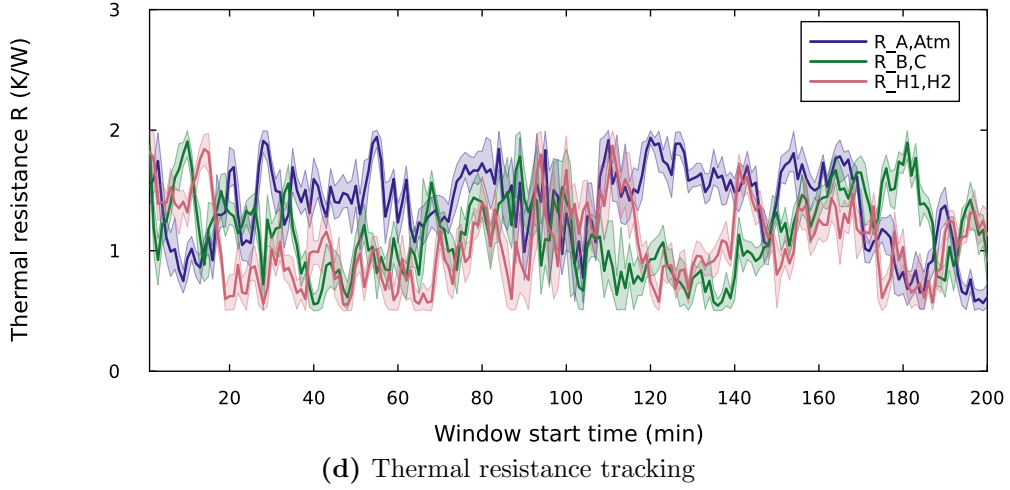}%
    \label{fig:R_tracking}%
  }
    \caption{Tracking of model parameters under moving-window Bayesian inference, showing the posterior mean (solid lines) and 95\% credible intervals (shaded bands).}
    \label{fig:params-tracking}
\end{figure}

% %%%%%%%%%%%%%%%%%%%adjusting figure%%%%%%%%%%%%%%%%%%

The corresponding evolution of the inferred parameters of the coupled \cotwo and temperature model %airflow-model 
are shown in Figure~\ref{fig:params-tracking}. The top panel reports the posterior mean and 95\% credible intervals of the occupancies in rooms A, E and F as a function of the window start time. Before $t \approx 120$~min, the chain assigns approximately one person to room~A and almost zero occupants to rooms~E and~F. As the windows cross the transition, the mean occupancy in room F increases sharply to four people, while the occupants in rooms A and E drop to zero, which is in agreement with the synthetic ground truth. Panel (b) displays the posterior mean and 95\% credible intervals of the boundary airflow from the ambient to the rooms A–E. Inferred boundary flows remain close to their true values and evolve smoothly across windows, with only a moderate increase in uncertainty around the transition period. %Together, Figures~\ref{fig:moving-window-6} and~\ref{fig:air-params-tracking} demonstrate that the moving-window Bayesian scheme not only produces accurate short-term \cotwo forecasts, but also tracks the underlying occupancy and airflow parameters physically consistent. 
Panel (c) reports the inferred thermal capacitances \(C_i\) for all zones. Across the experiment, the inferred capacities remain close to the ground-truth magnitudes and preserve the expected grouping between the smaller offices (A--D), larger rooms (E--F), and hallways (H1--H2). The bottom panel shows the inference results of representative thermal resistances. Compared with the capacities, the resistance tracks show noticeably larger window-to-window variability, which is expected because multiple resistance configurations can yield similar temperature trajectories over a finite window. Overall, it indicates that the moving-window Bayesian inference tracks RC parameters in a temporally coherent and physically plausible range under the controlled synthetic ground truth.

% The inferred temperature-network parameters is shown in Figure~\ref{fig:thermal-params-tracking}. The top panel 

% \begin{figure}[H]
%   \centering
%   % Row 1
%   \subfloat[Thermal capacity tracking]{%
%     \includegraphics[width=0.8\textwidth,trim={1cm 13cm 0cm 0cm},clip]{Figures/Results/track_parameter/thermal_params_tracking.pdf}%
%     \label{fig:C_tracking}%
%   }

%   \subfloat[Thermal resistance tracking]{%
%     \includegraphics[width=0.8\textwidth,trim={1cm 1cm 0cm 13cm},clip]{Figures/Results/track_parameter/thermal_params_tracking.pdf}%
%     \label{fig:R_tracking}%
%   }

%   \caption{Tracking of temperature-network (RC) parameters under moving-window Bayesian inference.
%   (a) Posterior mean and 95\% credible intervals of the inferred thermal capacities \(C_i\).
%   (b) Posterior mean and 95\% credible intervals of a representative subset of thermal resistances \(R_{i,j}\).}
%   \label{fig:thermal-params-tracking}
% \end{figure}

\subsubsection{Noise level}
\label{subsec:Noise level}

As discussed in Section~\ref{subsec:bayes_rule}, the nodewise noise scales $\sigma_{m,j}$ are inferred jointly with the other parameters in each moving window. While the synthetic benchmark data are generated with a spatially and temporally uniform \cotwo data noise level of $5$~ppm and temperature data noise level of $0.1$~\textcelsius , the inferred noise scales may vary across windows. The inferred noise level is a measure of the precision with which the the Bayesian framework is able to reproduce the data at each local time window. %\zw{I am not sure if return is the correct word.}

\begin{figure}[htbp]
  \centering
  \includegraphics[width=0.8\linewidth]{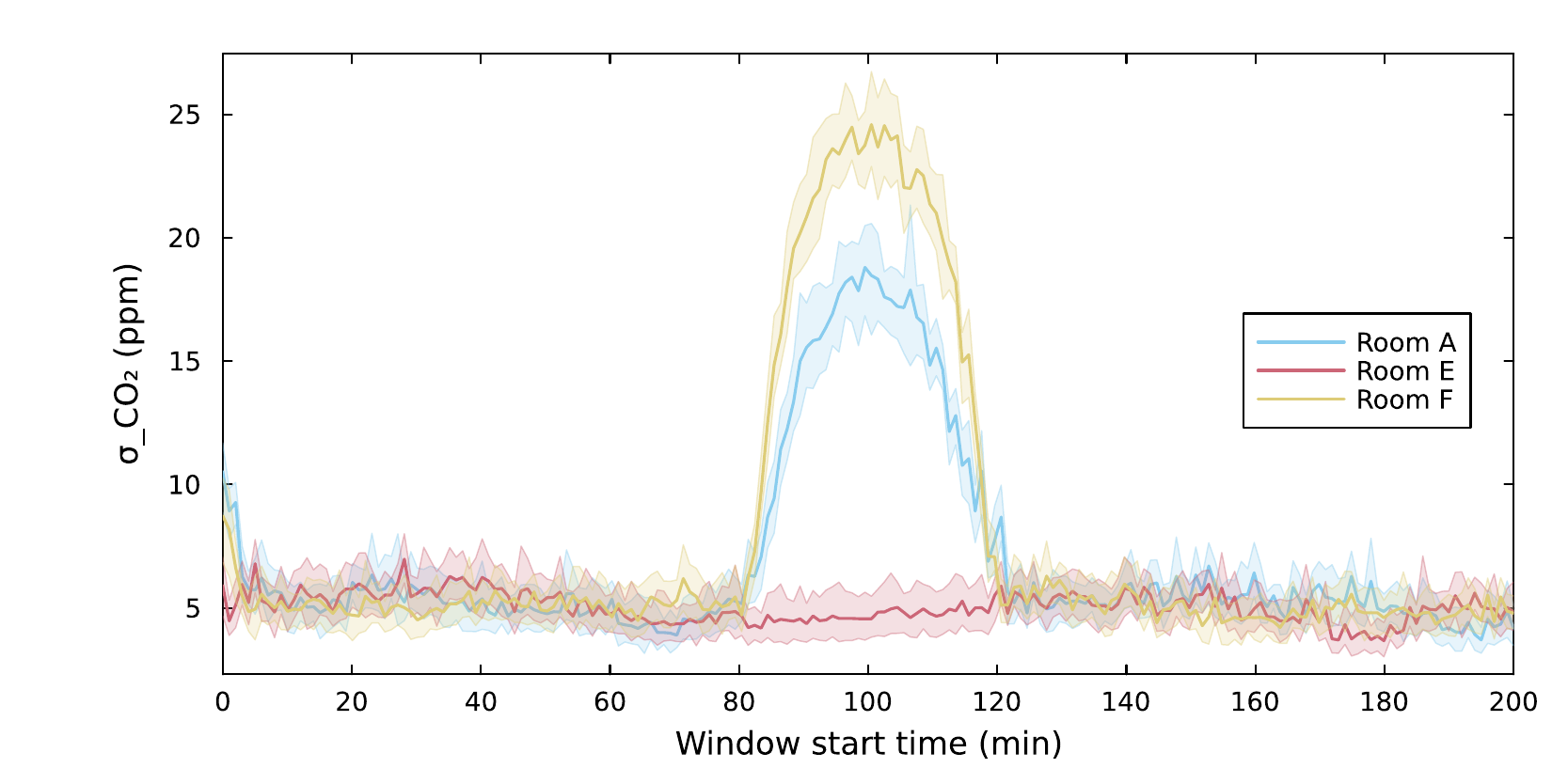}
  \caption{Inferred \cotwo noise levels in rooms A, E and F. Solid lines denote posterior means and shaded bands indicate 95\% credible intervals as a function of the window start time.}
  \label{fig:Noise CO2}
\end{figure}

Figure~\ref{fig:Noise CO2} shows the posterior mean and 95\% credible intervals of the inferred \cotwo noise levels $\sigma_{\mathrm{CO_2},j}$ in rooms A, E and F of each window start time with window size of 40 data points. For most windows, the posterior mean values remain close to the synthetic noise level. It shows that Bayesian inference is able to recover the measurement dispersion when the model is an accurate reflection of reality. %there is no change of occupancy or ventilation patterns. 

Around the occupancy transition at $t=120$~min, the behavior becomes more varied. For room F, where the occupants gathered while the \cotwo dynamics changed rapidly, the inferred $\sigma_{\mathrm{CO_2}}$ shows a significant increase: the posterior mean increases to values around 25 ppm until the window moves past the transition. Room A, where one person moves out, shows a less rapid change while room E almost shows no change in noise level. In that case, the noise level signal can be used for diagnosing if system behavior within the window size cannot be appropriately captured by the model.

\subsubsection{Change window size and noise level}
\label{subsec:Change window size and noise level}
To evaluate the sensitivity of the Bayesian inference to the moving-window size and measurement noise, we vary the inference window size from 10 to 60 data points under noise level of synthetic \cotwo and temperature data for four pairs of configuration $(\sigma_{\cotwo},\sigma_T)\in\{(5,0.1),(10,0.2),(15,0.3),(20,0.4)\}$. Performance is quantified by the nRMSE of prediction and ground-truth data in the verification interval $[80,120]$~min.

\begin{figure}[H]
  \centering
    % Row 1
    \subfloat[\cotwo]{\includegraphics[width=0.5\textwidth] {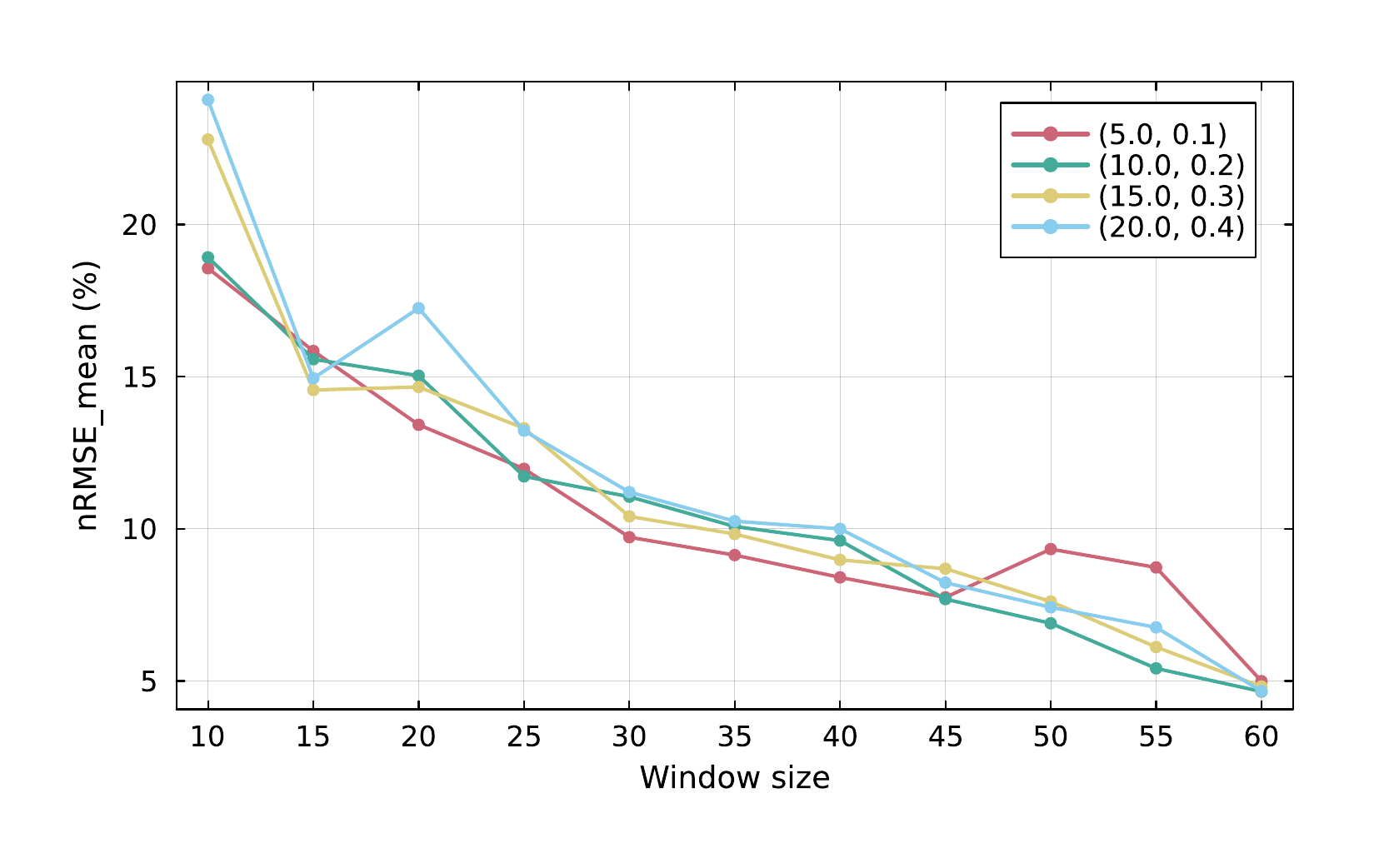}\label{fig:co2period1}}
    \subfloat[Temperature]{\includegraphics[width=0.5\textwidth]{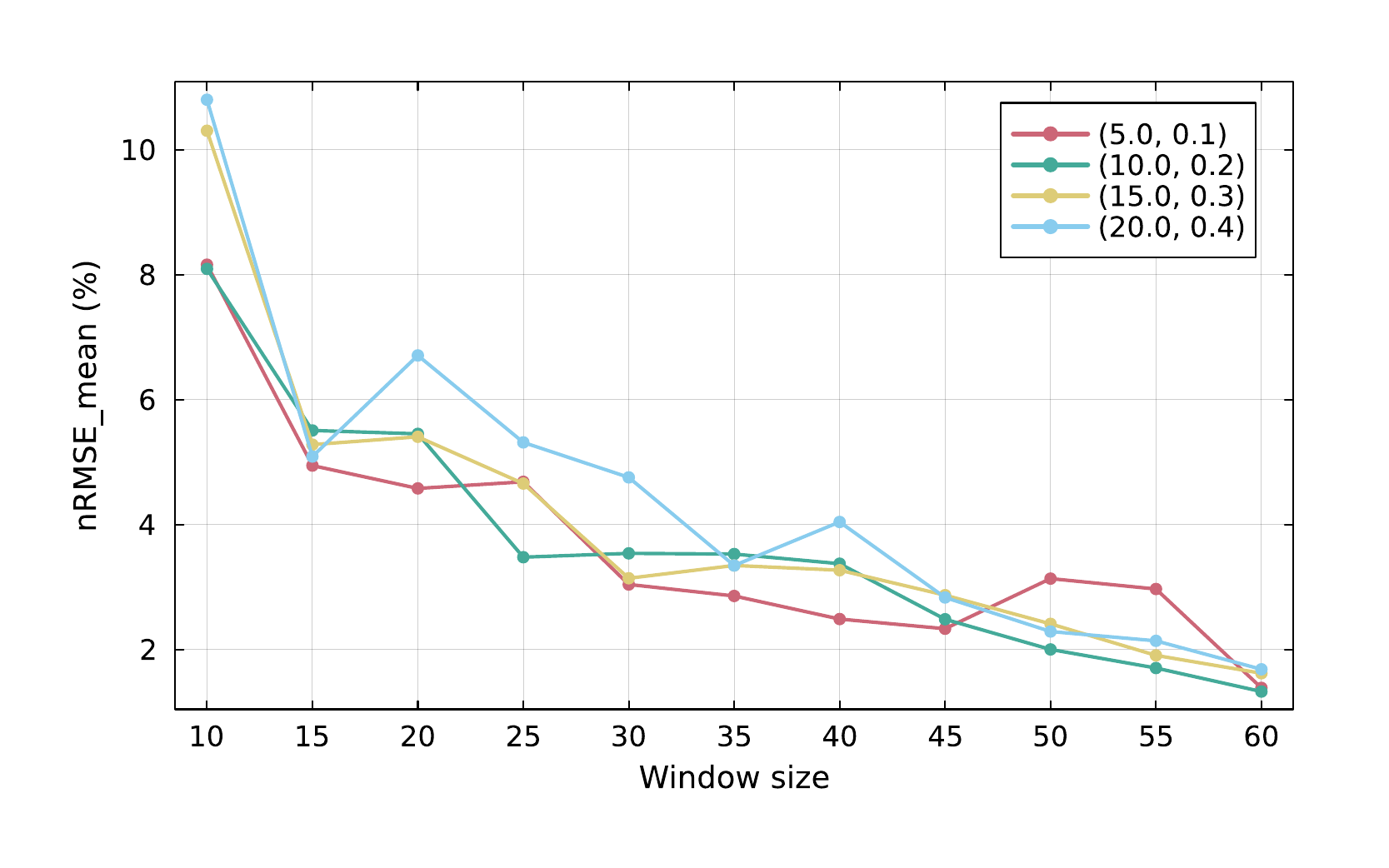}\label{fig:tempperiod1}}

    \caption{nRMSE of prediction and ground-truth data, evaluated under different inference window size 10--60 and paired measurement-noise settings $(\sigma_{\cotwo},\sigma_T)$. (a) shows \cotwo and (b) shows temperature.}
    \label{fig:noise_window_size}
\end{figure}

Figure~\ref{fig:noise_window_size} shows that increasing the size of the inference window consistently reduces nRMSE for both \cotwo and temperature across all noise levels. 
The impact of measurement noise is most pronounced for shorter window sizes. At 10--25 size, higher noise levels lead to higher nRMSE, particularly for \cotwo. At larger window sizes, the performance curves converge, and the difference between noise settings becomes small. This indicates that larger window sizes partly compensate for the increased noise by providing more samples per inference step, and therefore the inference can better estimate the noise level of the data.

The appropriate choice of window size depends on the required model accuracy. For \cotwo (Fig.~\ref{fig:co2period1}), nRMSE decreases from 19--24\% at window size 10 to around 5\% at a window size 60, with most of the improvement occurring between 10 and 35. For temperature (Fig.~\ref{fig:tempperiod1}), errors are substantially lower overall and decrease from  8--11\% at 10 to  $1.3\sim1.7\%$ at 60.

%=====================================================
%=====================================================
\section{Experimental validation}
\label{sec:Experiment validation}
To assess the predictive performance of our modeling framework on real measurement data that correspond to time-varying operating conditions, we now implement the moving-window Bayesian inference framework for a (scaled) physical experiment of an office building.

\subsection{Scaled experiment}
\label{subsec:Scaled setup}
First, we describe the setup of the physical experiment and the data acquisition procedure.

\subsubsection{The physical setup}
\label{subsubsec:The setup}
\begin{figure}[htbp]
  \centering
  \includegraphics[width=0.8\linewidth]{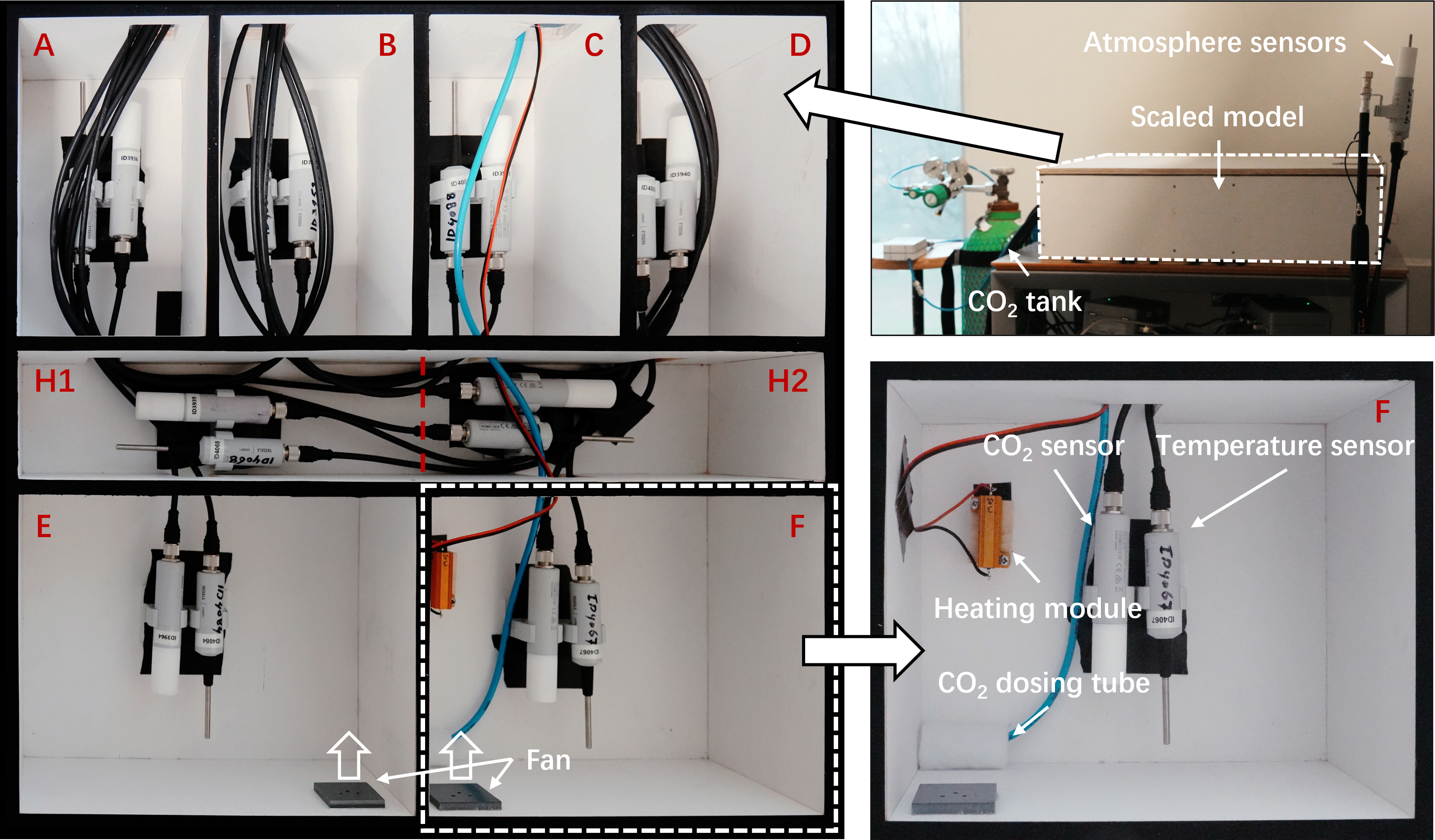}
  \caption{Experiment setup}
  \label{fig:experiment setup whole}
\end{figure}

%To validate the developed framework, a physical scaled experiment was considered. (
Figure \ref{fig:experiment setup whole} shows the experimental setup. It corresponds to the spatial layout of the synthetic building from Section~\ref{sec:Synthetic benchmark problem}. Four smaller rooms (A--D) and two larger rooms (E and F) are connected by a hallway (H1 and H2). There are openings in the walls that represent windows and doors, allowing air exchange between the rooms and the environment. 

To impose controlled airflow conditions, two $30\times 30$~mm square fans are installed in rooms E and F. The fans are operated at fixed electrical settings (powered by 1.8~V) to create a controlled flow. The fans are mounted on the wall and oriented to drive flow from the ambient to the room as indicated in Figure~\ref{fig:experiment setup whole}. The fans are partially blocked to reduce the airflow to a reasonably level, so that the resulting air-exchange rate in the scaled setup produces \cotwo concentration changes that are traceable given the sensors’ 10-second sampling resolution. The fan operating state is kept constant throughout the experimental run.

The rooms and the ambient are equipped with a total of nine pairs of \cotwo concentration and temperature sensors. The installed sensors are the Vaisala Carbon Dioxide Probe GMP252~\cite{vaisala_gmp252} and the Vaisala Humidity and Temperature Probe HMP1~\cite{vaisala_hmp1}, for which the accuracy is documented as ±40 ppm for the \cotwo sensor and $\pm0.2^\circ\mathrm{C}$ for the temperature sensor in the range of 0–3000 ppm.
The sensors are placed at the room center positions. Atmospheric reference sensors are installed outside the enclosure to record background \cotwo and temperature for boundary conditions.

To simulate human-body CO$_2$ and heat production in room F, a premixed CO$_2$ gas cylinder and a heating module were used. The premixed gas had a nominal CO$_2$ concentration of 5000~ppm, as supplied by Nippon Gases \cite{nippon_gases_certified_gas}. The dosing flow rate is regulated to $25$~mL/min using a flow meter (GASTEC Model GSP-300FT-2 \cite{gastec_gsp300ft2}). %Gas was delivered via a dosing tube connected between the flow meter and the room F. 
The end of the dosing tube is covered by foam to promote mixing with the air in the room. %Meanwhile, the \zw{check this in the lab: resistive heating module} was installed in room F to generate controlled heating power. 
The heating element is operated with a constant 25 Watt power.

\subsubsection{Experimental conditions}
\label{subsec:Experimental conditions}

The scaled experiment considers a single controlled run in which the fan-driven airflow is kept constant and the described combined \cotwo and heating source is applied for 30 minutes and then turned off to simulate a clear point of change.
%The heater power was set to $25$~W, and the \cotwo source was a premixed cylinder at nominal $5000$~ppm delivered to room F through the dosing tube.
Table \ref{tab:scaled_experiment_schedule} shows detailed timings of the experiment. We initialize the experiment with a short pre-run calibration period (0--10 min) to compute sensor offsets for both sensor types. Then, the \cotwo dosing and heating module are switched on (10--40 min). During the experiment, the two fans in rooms E and F remain at fixed electrical settings, and the \cotwo dosing and heating in room F are applied constantly. At 40 min, the \cotwo dosing and heating module are closed while the fans are kept running for another 30 minutes (40--70 min). All signals are logged at 10 seconds intervals.

\begin{table}[htbp]
\centering
\caption{Scaled experiment run plan used for Bayesian inference.}
\label{tab:scaled_experiment_schedule}
\begin{tabular}{lll}
\hline
Phase & Time (minute) & Condition \\
\hline
Offset calibration & 0--10 & No intervention; used to align sensor baselines \\
Dosing& 10--40 & Fans on; \cotwo dosing into room F; heater on in room F \\
Decay& 40--70 & \cotwo dosing off and heater off; fans kept constant \\
\hline
\end{tabular}
\end{table}

The raw sensor data are first converted to time-series for the eight zones (A--F, H1--H2). To reduce inter-sensor bias, we perform an offset calibration using the 0--10 minute baseline period: for each room and sensor type (\cotwo and temperature), we subtract the room-specific mean deviation from the grand mean over that calibration interval. After calibration, the inference interval 10--70 min is extracted and used as input to the moving-window algorithm.

\subsubsection{Inference configuration}
\label{subsec:Scaled experiment inference configuration}

The inferred parameter vector follows the coupled layout used throughout Section~\ref{sec:Bayesian} and includes the (i) \cotwo and temperature evolution parameters (8 effective source strengths and 5 boundary airflow controls for rooms A--E, 19 resistances and 8 capacitances), (ii) initial states for \cotwo and temperature in each room, and (iii) per-room measurement noise standard deviations for \cotwo and temperature. The likelihood assumes independent Gaussian measurement noise per room, constant over each inference window. Table \ref{tab:scaled_priors_scaled} presents the full prior configuration for the Bayesian inference in the validation experiment. We use truncated normal priors consistent with the synthetic benchmark: thermal resistances $R$ and capacitances $C$ are constrained to physically plausible ranges, boundary airflow parameters are constrained to a small range (in L/s) to represent low-speed fan-driven or leakage-scale exchanges, and per-room noise parameters are constrained to positive ranges.

The Bayesian estimation uses a moving-window size of $40$ samples (corresponding to 400 seconds), with a step size of $6$ samples (corresponding to 60 seconds), producing overlapping windows across the full 0--70 min recording until the end of the window reaches $t=70$ min. For each window, we sample the log-posterior using the robust adaptive Metropolis algorithm as for the synthetic experiment. To encourage temporal continuity and reduce the burn-in cost, the posterior mean of the previous window is used to initialize the next window's chain. Posterior predictive trajectories are generated from the posterior mean for visualization and parameter tracking.

\begin{table}[H]
\centering
\caption{Prior distributions for all parameters that are inferred in the validation experiment.}
\label{tab:scaled_priors_scaled}
\renewcommand{\arraystretch}{1.25}
\begin{tabular}{p{0.32\textwidth} p{0.6\textwidth}}
\hline\\[-10pt]
\textbf{Occupancies $n_{\mathrm{ppl},i}$} &
Rooms A--D: TruncNorm$(1,\,1;\,0,\,1)$; \\
&Rooms E--F: TruncNorm$(1,\,1;\,0,\,2)$; \\
&Hallway H1--H2: TruncNorm$(0,\,0.3;\,0,\,0.3)$. \\[7pt]

\textbf{Boundary airflows $q_{\mathrm{Atm},i}$ [L/s]} &
TruncNorm$(0,\,0.05;\,-0.1,\,0.1)$ (for the 5 boundary nodes). \\[7pt]

\textbf{Initial \cotwo states $\cair_{0,i}$ [ppm]} &
$\mathcal{N}\!\big(\tilde{\cair}_i(t{=}0),\,\sigma_{\cair,i}^2\big)$ per node (anchored to the first observation in each window). \\[7pt]

\textbf{Thermal resistances $R_{i,j}$ [K/W]} &
TruncNorm$(0.5,\,2.0;\,0.05,\,8.0)$ for all 19 resistances. \\[7pt]

\textbf{Thermal capacitances $C_i$ [J/K]} &
TruncNorm$(3000,\,2000;\,500,\,15000)$ for all 8 capacitances. \\[7pt]

\textbf{Initial temperatures $T_{0,i}$ [$^\circ$C]} &
$\mathcal{N}\!\big(\tilde{T}_i(t{=}0),\,\sigma_{T,i}^2\big)$ per node (anchored to the first observation in each window). \\[7pt]

\textbf{\cotwo noise levels $\sigma_{\cair,i}$ [ppm]} &
$\sigma_{\cair,i} \sim \text{TruncNorm}(5.0,\,2.0;\,0,\,8.0)$, 
\quad $i\in\{\mathrm{A},\dots,\mathrm{H2}\}$. \\[7pt]

\textbf{Temperature noise levels $\sigma_{T,i}$ [$^\circ$C]} &
$\sigma_{T,i} \sim \text{TruncNorm}(0.1,\,0.1;\,0,\,0.5)$, 
\quad $i\in\{\mathrm{A},\dots,\mathrm{H2}\}$. \\[6pt]
\hline
\end{tabular}
\end{table}

%=====================================================
%=====================================================
%=====================================================
%=====================================================
\subsection{Validation results}
\label{subsec:Experiment validation result}
We validate the moving-window Bayesian inference scheme on the scaled experiment by examining (i) posterior predictive trajectories in both \cotwo and temperature, (ii) the corresponding time-varying parameter tracks for occupancy and boundary exchange, and (iii) post-window predictive accuracy quantified via nRMSE over a fixed prediction horizon.

\subsubsection{Predictive posterior results }

\begin{figure}[H]
  \centering
    % Row 0
    \subfloat[Window start at 3 min, \cotwo.]{\includegraphics[width=0.5\textwidth,trim={0.7cm 1cm 13cm 0.3cm},clip] {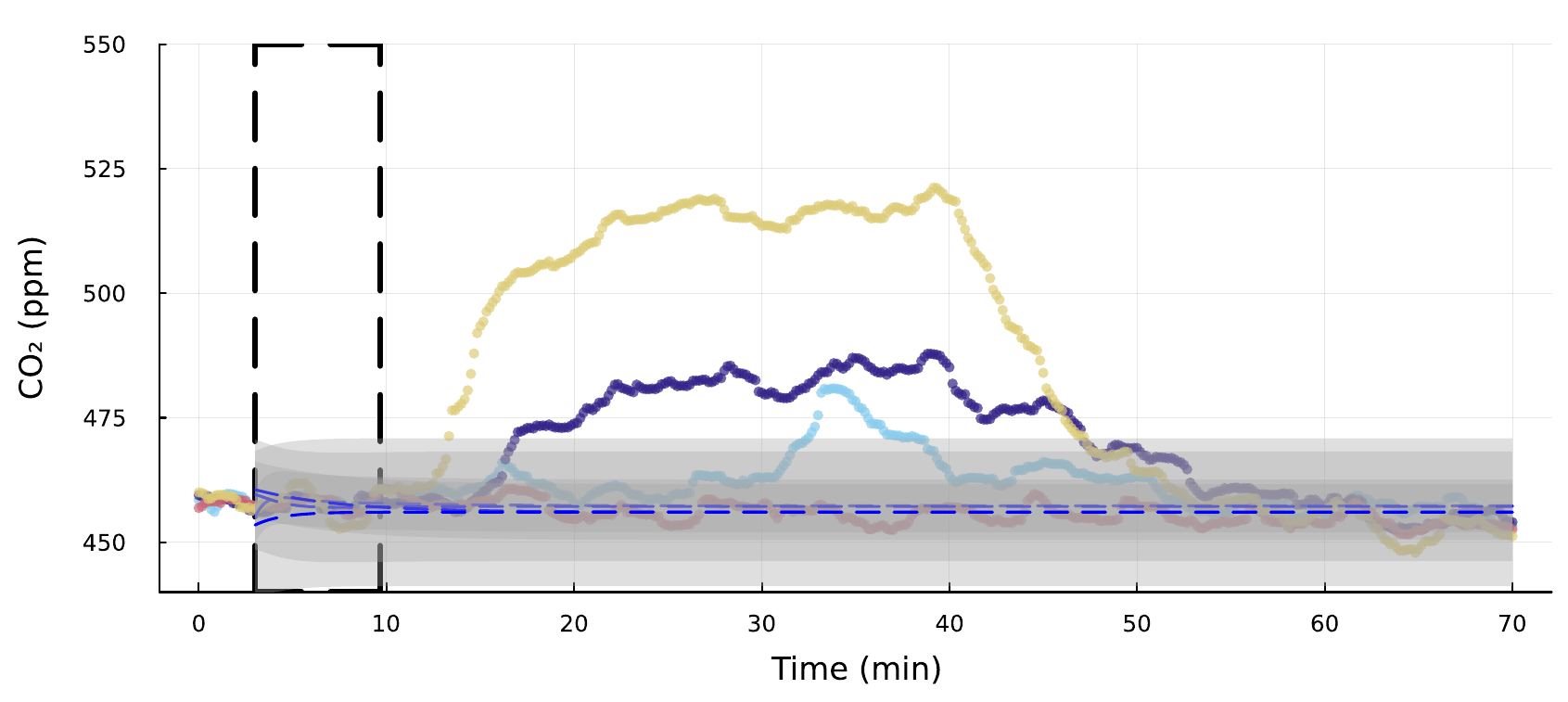}\label{fig:scaledCO2Window19}}
    \subfloat[Window start at 3 min, temperature.]{\includegraphics[width=0.5\textwidth,trim={0.7cm 1cm 13cm 0.3cm},clip]{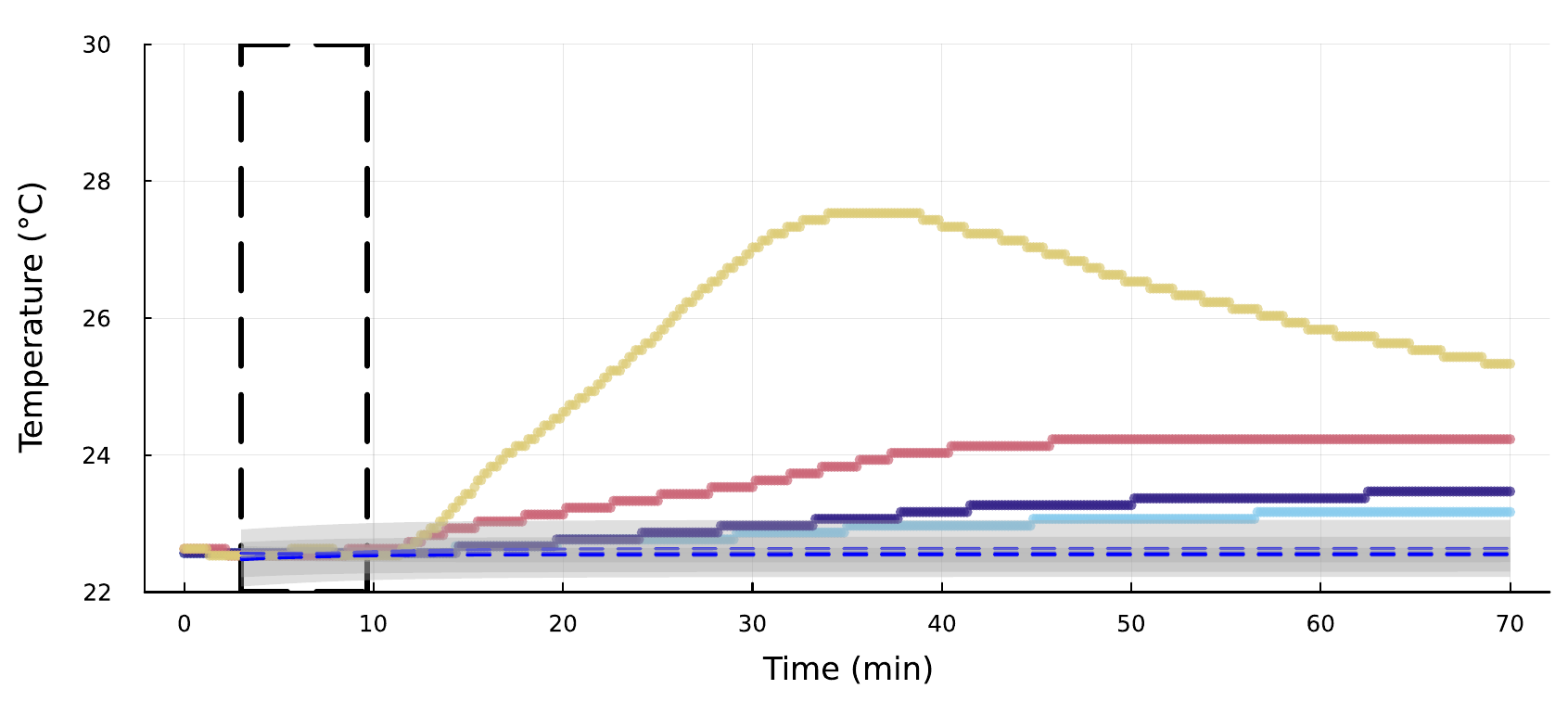}\label{fig:scaledTempWindow19}}
  
    % Row 1
    \subfloat[Window start at 20 min, \cotwo.]{\includegraphics[width=0.5\textwidth,trim={0.7cm 1cm 13cm 0.3cm},clip] {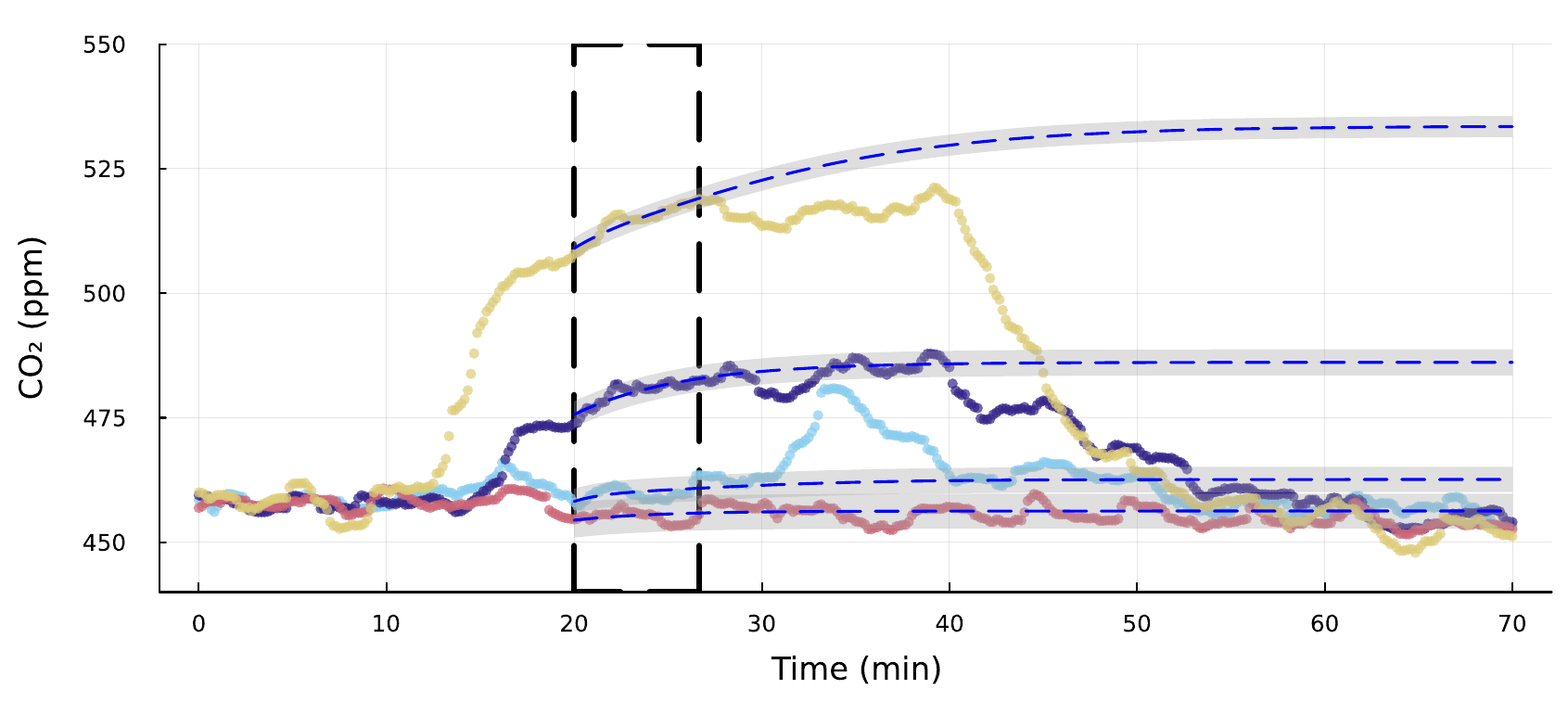}\label{fig:scaledCO2Window121}}
    \subfloat[Window start at 20 min, temperature.]{\includegraphics[width=0.5\textwidth,trim={0.7cm 1cm 13cm 0.3cm},clip]{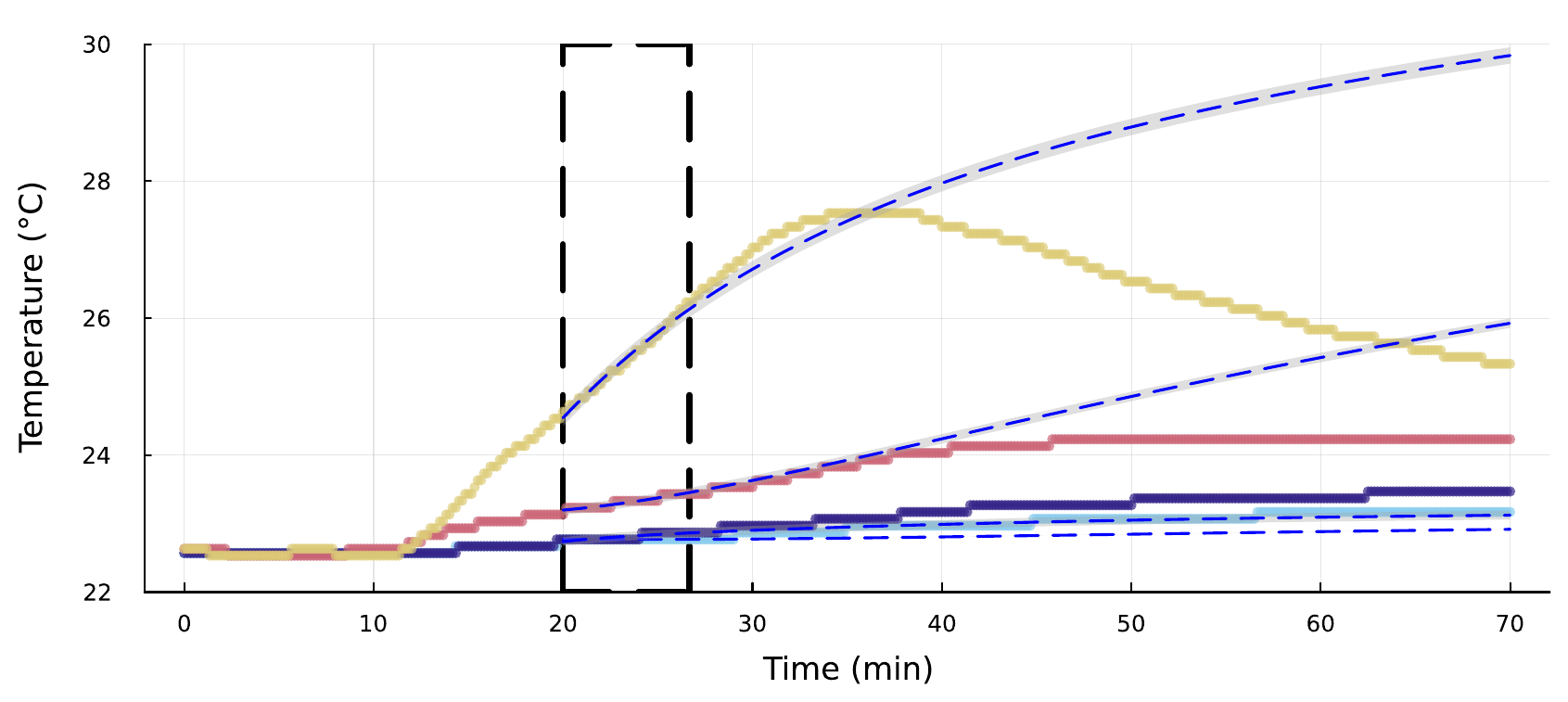}\label{fig:scaledTempWindow121}}
    
    % Row 2
    \subfloat[Window start at 31 min, \cotwo.]{\includegraphics[width=0.5\textwidth,trim={0.7cm 1cm 13cm 0.3cm},clip] {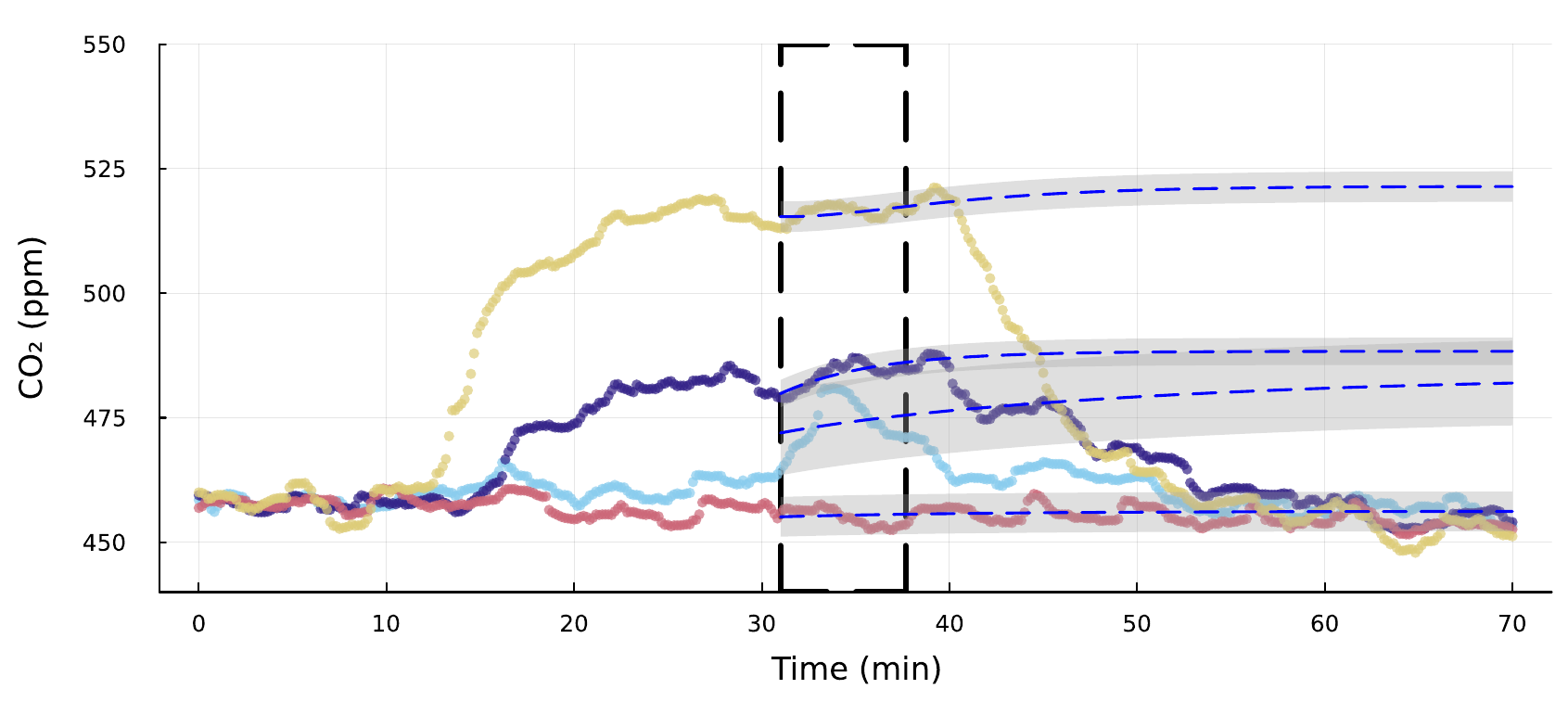}\label{fig:scaledCO2Window187}}
    \subfloat[Window start at 31 min, temperature.]{\includegraphics[width=0.5\textwidth,trim={0.7cm 1cm 13cm 0.3cm},clip]{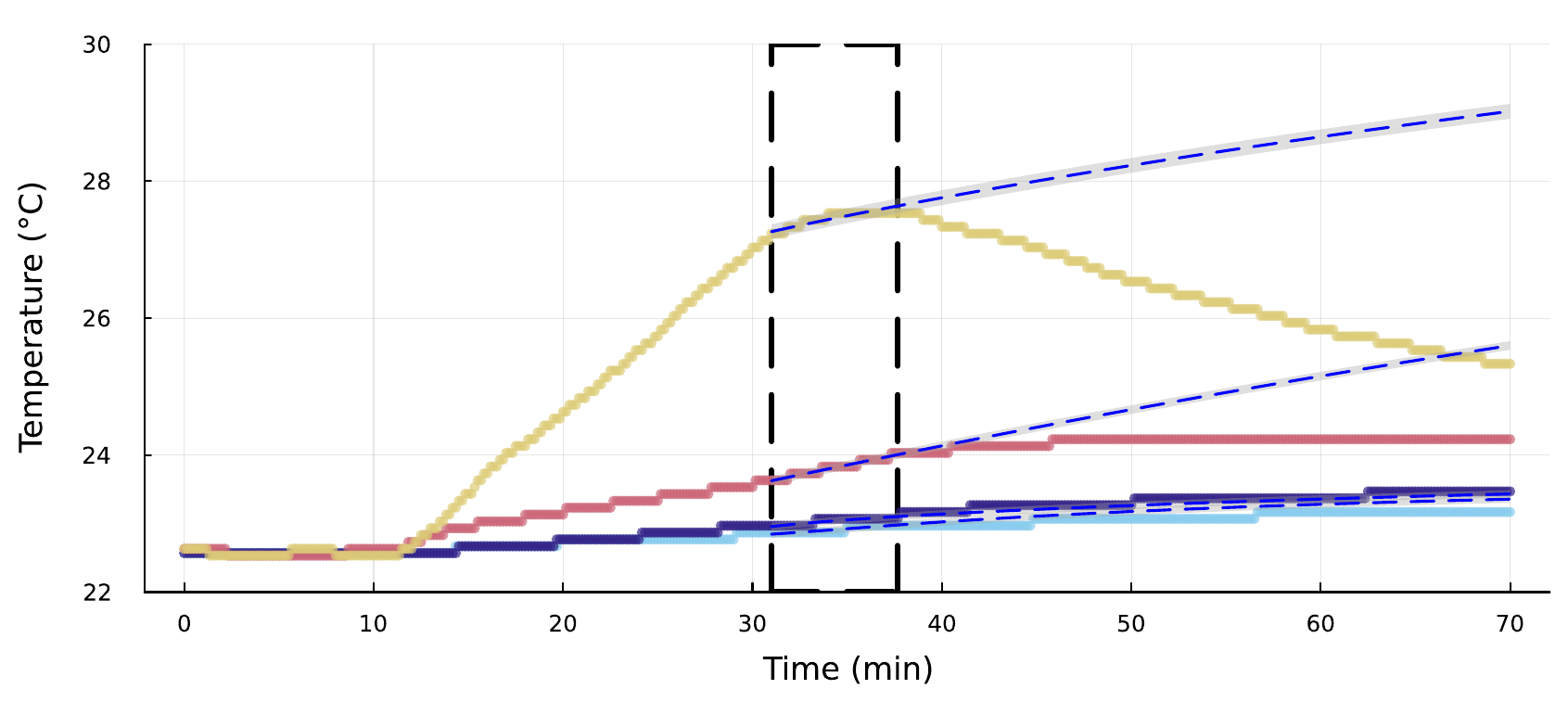}\label{fig:scaledTempWindow187}}

    % row 3
    \subfloat[Window start at 44 min, \cotwo.]{\includegraphics[width=0.5\textwidth,trim={0.7cm 1cm 13cm 0.3cm},clip] {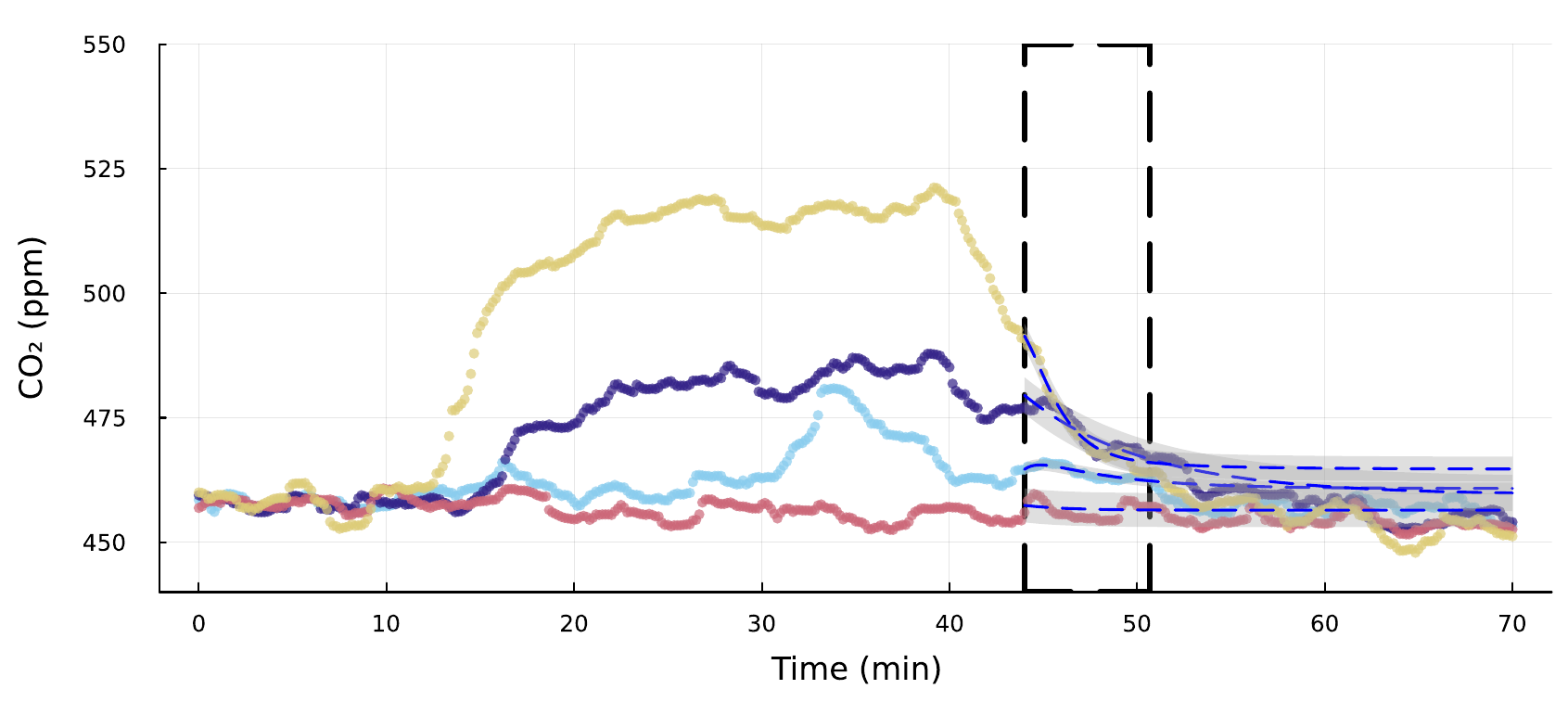}\label{fig:scaledCO2Window265}}
    \subfloat[Window start at 44 min, temperature.]{\includegraphics[width=0.5\textwidth,trim={0.7cm 1cm 13cm 0.3cm},clip]{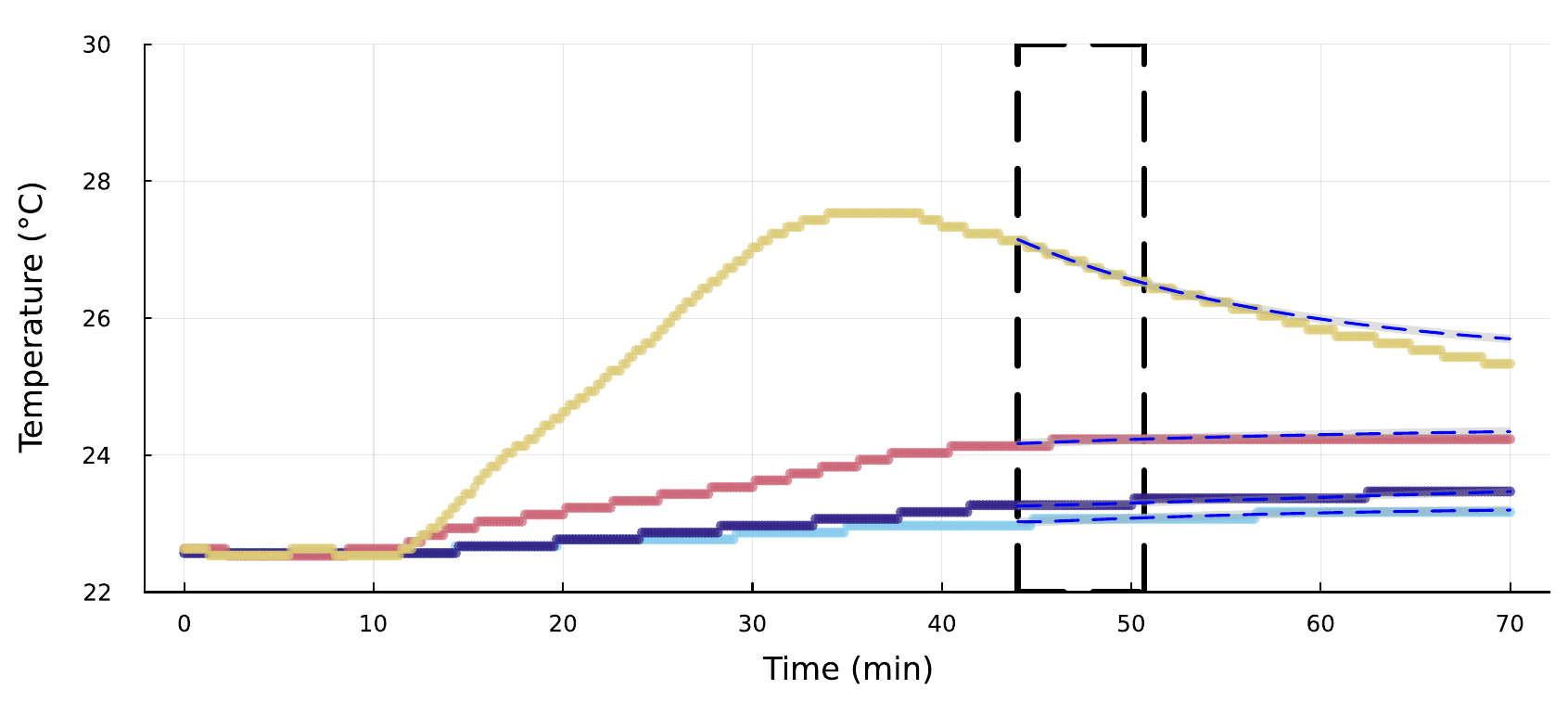}\label{fig:scaledTempWindow265}}

    \subfloat{\includegraphics[width=0.5\textwidth,trim={20.5cm 6.8cm 2.5cm 5cm},clip]{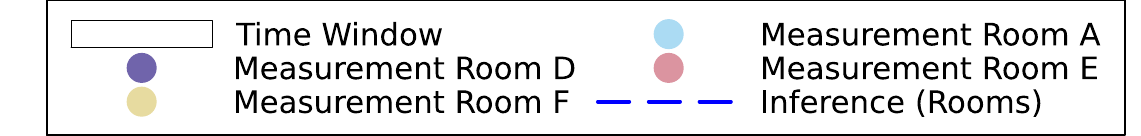}\label{fig:scaled_legend}}
    
    \caption{Posterior predictive \cotwo concentration (left column) and temperature (right column) trajectories in rooms A, D, E, and F, for the scaled physical experiment. %obtained with the inference windows.
    }
  \label{fig:scaled-moving-window}
\end{figure}

Figure~\ref{fig:scaled-moving-window} summarizes the behavior of the Bayesian moving-window algorithm on the experiment by showing posterior predictive trajectories for three representative windows for both \cotwo and temperature. In each panel, the dashed rectangle indicates the data interval used for inference in that step, the colored markers show the measured room trajectories (rooms A, D, E, and F), the dashed curves report the posterior predictive mean obtained from the posterior mean parameter set, and the gray ribbons indicate the associated 95\% credible bands. %Because each window only conditions on the measurements inside the dashed window, agreement within the window reflects local calibration performance, while deviations immediately outside the window highlight model prediction behavior and potential mismatch when operating conditions change.

The early window in the first row (panels (a)–(b), start time 3 min) corresponds to the early inference steps (0--10 min) and is taken before the main system dynamics initiate. Within this window, both \cotwo and temperature remain close to their background levels in all rooms, with only small fluctuations around the baseline, which leads to a nearly flat posterior predictive mean. % Comment on thickness of 95% confidence interval, if this is not a bug

Windows in the second row (panels (c)–(d), start time 20 min) lie in the occupancy-on regime, where \cotwo dosing and heating are applied in room F. In this rising phase, the posterior predictive mean closely tracks the within-window increase in both \cotwo and temperature in room F, while the other rooms exhibit smaller responses consistent with heat exchange through the connected network layout. Immediately after the right edge of the window, the predictions continue the prevailing trends, which is expected because the model has not yet observed the subsequent change point. In the third row window (panels (c)–(d), start time 31 min), the inference operates later in the excitation period, where the room-F trajectories approach a near-peak and the differences between rooms become more pronounced. Here, the moving-window fit remains locally consistent, and the credible bands tend to broaden when the dynamics are steep or when the window goes across a transition in slope, indicating elevated parameter and state uncertainty under rapidly changing conditions.

The window in the fourth row (panels (g)–(h), start time 44 min) is positioned in the occupancy-off regime, where both \cotwo dosing and heating are disabled while the fan-power remains constant. In this decay phase, the model captures the main decay in the driven room~F for both \cotwo and temperature, and smaller changes elsewhere.

Overall, Figure~\ref{fig:scaled-moving-window} conveys that the moving-window posterior can provide physically interpretable, locally accurate reconstructions of the measured \cotwo and temperature trajectories. However, unmodeled behavior, such as non-uniform mixing, boundary leakage, or fan-induced jet patterns, might lead to systematic mismatch. The modeling framework is able to identify such mismatches by increased width of the credibility bands and recovers when the mismatches disappear. %These observations motivate the parameter-tracking results in Section \ref{subsec:Parameter inference}, where inferred occupancy and airflow parameters are examined.

\subsubsection{Parameter inference}
\label{subsec:Parameter inference}

\begin{figure}[H]
  \centering
    % Row 1
    \subfloat[Occupancy.]{\includegraphics[width=0.8\textwidth,trim={1cm 12cm 0cm 0cm},clip] {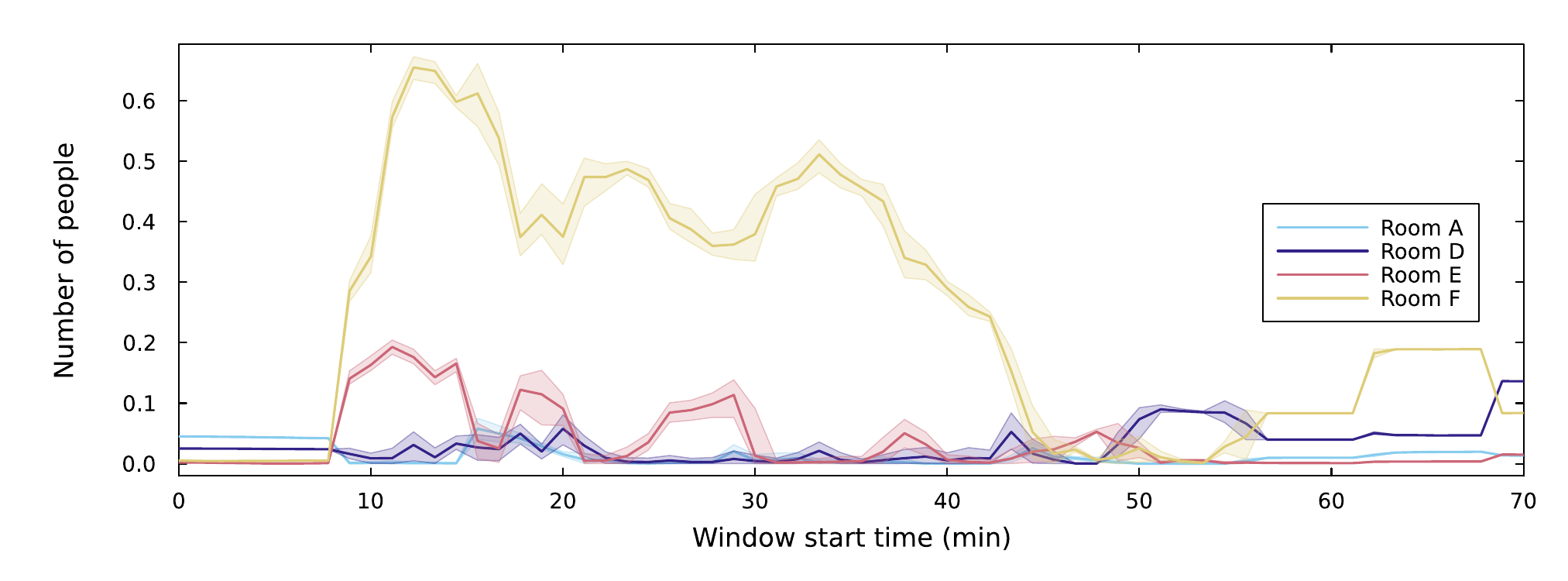}\label{fig:scale_occupancy_tracking}}
    
    \subfloat[Airflow.]{\includegraphics[width=0.8\textwidth,trim={1cm 0cm 0cm 12cm},clip]{Figures/Scaled_results/scale_air_params_tracking.pdf}\label{fig:scale_airflow_tracking}}

    \caption{ Posterior mean (solid lines) and 95\% credible intervals (shaded bands) of the inferred occupancies and façade airflow, for the scaled physical experiment.}
    \label{fig:scale-air-params-tracking}
\end{figure}

Figure~\ref{fig:scale-air-params-tracking} shows the moving-window tracking of the parameter groups that drive both balance models: the equivalent occupancy/source-strength parameters (Figure~\ref{fig:scale_occupancy_tracking}) and the boundary airflow control parameters for rooms~A--E (Figure~\ref{fig:scale_airflow_tracking}). In both panels, solid curves denote posterior means and the shaded ribbons show 95\% credible intervals, plotted against the window start time. Smooth trends indicate persistent evidence across overlapping windows, while widened bands or abrupt shifts indicate locally weak identifiability or a regime transition. The occupancy parameters enter the \cotwo and temperature balance equations as source terms (fixed exhalation concentration and fixed per-person sensible heat), so they should be interpreted as equivalent occupancies rather than literal head counts.

During the initial calibration interval (0--10~min), where the measurements show no appreciable \cotwo excitation, both panels remain close to their baseline levels and are therefore largely influenced by the sampler initialization (weak data constraints). Later, Figure~\ref{fig:scale_occupancy_tracking} assigns the dominant source to room~F during the excitation period (10--40~min), and this source strength drops rapidly after the switch-off (40~min) and remains near zero thereafter, showing that the windowed posterior responds sharply once the windows begin to include post-switch-off data. Meanwhile, Figure~\ref{fig:scale_airflow_tracking} shows that the inferred ambient exchange controls for rooms~A--E exhibit time-varying adjustment around the switch-off, indicating that the estimator reallocates effective boundary exchange to explain the observed decay and redistribution patterns in \cotwo and temperature once the imposed source in room~F is removed, while still satisfying the constraints of network mass-balance through the hallway nodes and room~F.

As a supplementary explanation, intermittent small non-zero equivalent-source estimates in non-actuated zones (including the hallways) should be interpreted as model--experiment mismatch compensation -- e.g., imperfect mixing, leakage paths, or fan-induced jet structure not represented by the well-mixed, low-order airflow parameterization -- rather than evidence of truly distributed sources. Correspondingly, variations in the boundary-control parameters should be read primarily as effective exchanges that support data-consistent transport, not as direct measurements of physical fan flow rates.

Overall, Figure~\ref{fig:scale-air-params-tracking} indicates that the moving-window Bayesian scheme can track regime changes in real measurements, while also motivating improved experiment-to-model alignment (e.g., explicit fan-actuation inputs and a decoupled parameterization of \cotwo vs.\ heat injection) to reduce compensation effects and improve physical interpretability.

\subsubsection{Prediction evaluation}
\label{subusbsec:Prediction evaluation}

To assess how well each posterior generalizes beyond the inference interval, we compute normalized root-mean-square error (nRMSE) between the model prediction and experimental measurements over a horizon of $80$ samples after the end of each inference window, and calculate the mean value over eight rooms, for both \cotwo and temperature. The evaluation is performed for every window end-time.

\begin{figure}[t]
  \centering
  \includegraphics[width=0.8\linewidth]{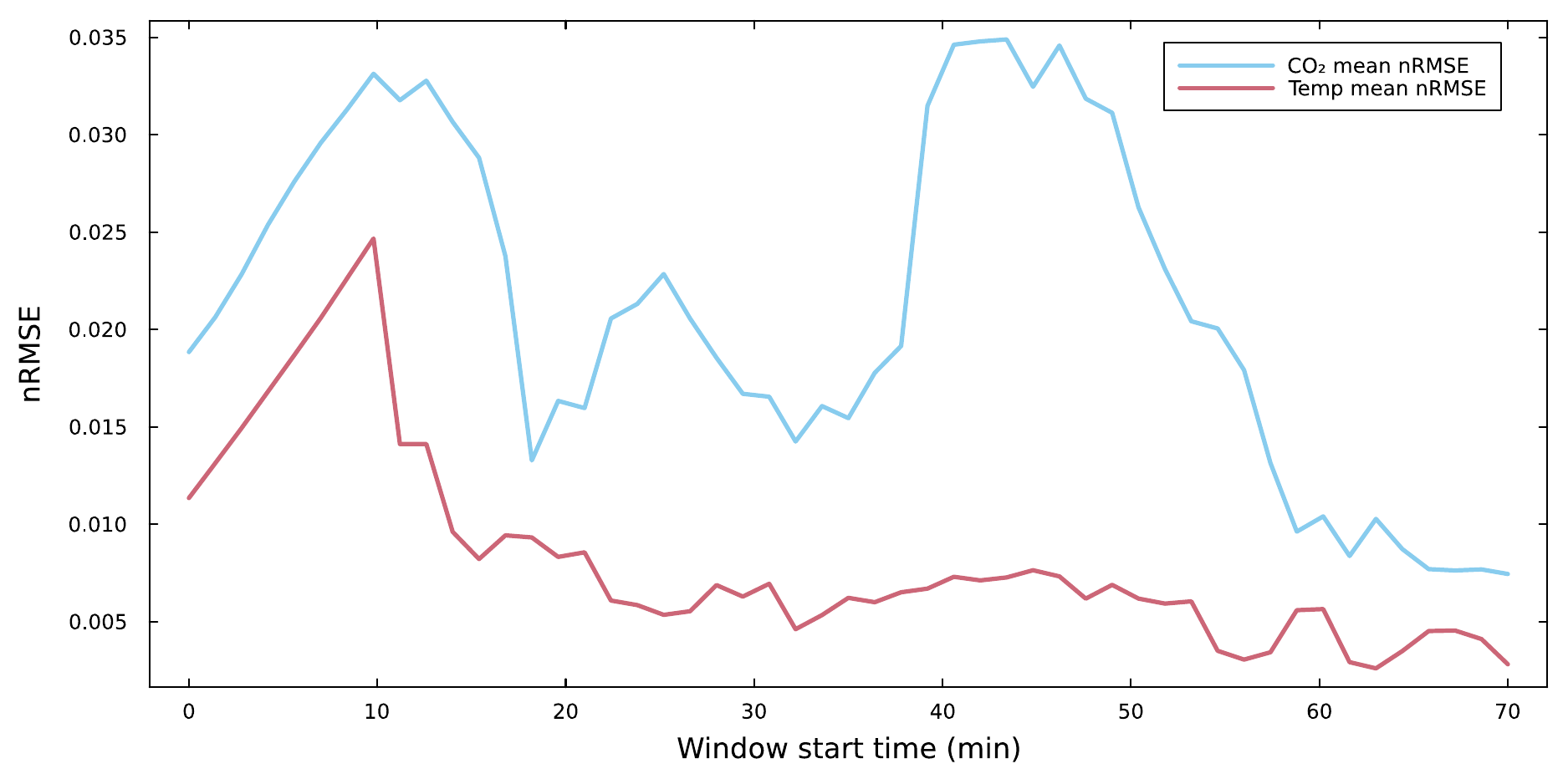}
  \caption{Mean nRMSE over an $80$ samples post-window prediction horizon, averaged over 8 rooms, for \cotwo and temperature.}
  \label{fig:nrmse_mean_h80}
\end{figure}

The results are presented in Figure~\ref{fig:nrmse_mean_h80}. It is observed that when the forecast interval remains within a relatively stationary operating regime, both \cotwo and temperature errors are low and comparatively stable, most clearly in the later part of the experiment where the curves decrease and flatten. Elevated errors occur for windows whose 80-step forecast overlaps a change point in the physical inputs (e.g., \cotwo dosing and heater on/off) or a rapid transient (steep rise/decay). In these cases, a window-local posterior can match \emph{local} dynamics, but model prediction cannot anticipate imminent input switching. %, so extrapolation accumulates bias over the prediction. 
The comparatively larger \cotwo nRMSE excursions (relative to temperature) match the earlier observation from the synthetic benchmark case, per Figure~\ref{fig:noise_window_size}. For this physical experiment, it indicates that the prediction of the short-horizon is more sensitive to unmodeled or abruptly changing source terms and ventilation-driven transport during the excitation transitions. A room-wise inspection (not shown in the figure) indicates that the largest forecast errors concentrate in the directly actuated zone (room~F), which is expected because it experiences the strongest source terms and the sharpest regime change. The remaining zones show smaller errors driven primarily by transport. This asymmetry suggests that forecast accuracy is limited less by the network propagation itself and more by the inability of a window-local parameterization to represent abrupt source/input transitions over the forecast horizon.

%=====================================================
%=====================================================
%=====================================================
%=====================================================

\section{Conclusion}
\label{sec:conclusion}

We proposed a coupled \cotwo--temperature network model that links multi-zone \cotwo transport and temperature dynamics through shared latent drivers, namely inter-zonal airflow and occupancy. We developed a moving-window Bayesian inference method to infer the occupancy, airflow, and temperature-network parameters from windowed \cotwo and temperature measurements, while simultaneously estimating per-room sensor noise levels and producing posterior predictive trajectories with uncertainty bands. 

In a synthetic benchmark experiment, we showed that the moving-window posterior accurately reconstructs the within-window trajectories and provides accurate model-based predictions, as long as the measurement data can be described by our coupled \cotwo--temperature network model. If the system evolution does not adhere to the constraints of the model (e.g., when the inference window overlaps an occupancy-regime transition), then this mismatch can be observed as a widening of the % Consistent with the expected behavior of a window-local estimator, 
posterior uncertainty and an increase in the predicted noise levels. The framework is able to recover quickly after such a phase. %around the change point and contracts once the new regime is fully observed. Moreover,
Furthermore, we studied the relation between window size and prediction error for different sensor noise levels. As expected, increasing the window size reduces the prediction error for both \cotwo and temperature. Interestingly, these trends and the error values are not strongly affected by sensor noise levels. %, indicating that longer windows provide stronger constraints for identifying the coupled dynamics. 

The validation experiment reinforced the conclusions drawn from the synthetic benchmark experiment. Our framework produced locally consistent reconstructions and low prediction errors %in both equivalent occupancy engage and leave phases, and parameter tracking also reflects the scenario. Forecast % errors remain low 
in relatively well-behaved regimes, but errors increase when the prediction window includes abrupt input changes or steep transients. In those cases, the posterior uncertainty bands widen and the predicted noise levels increase. The largest errors concentrated in the directly actuated room (where heat and \cotwo were released), indicating that short-horizon predictive accuracy is limited primarily by unmodeled behavior.

From these experiments, several practice-related limitations and suggestions should be noted. Firstly, the lumping procedure in our network model assumes well-mixed room states. This assumption is likely violated in practice, especially in scaled setups, leading to fundamental modeling flaws. Secondly, in a real life application, the inferred occupancy parameters should be interpreted as equivalent source strengths consistent with fixed per-person \cotwo and sensible heat constants. Parameter inference might fail when the experimental actuation does not perfectly match the parameterization of the model. Thirdly, identifiability depends on excitation. During periods with little to no heat or \cotwo excitation, parameters can only be weakly informed and uncertainty bands would be broad. Finally, the moving-window Markov-chain Monte-Carlo configuration, while robust, can be computationally demanding for real-time deployment, especially when scaling to larger buildings or shorter update intervals.

Based on the above, future research could focus on: (i) incorporating regime-change handling directly in the estimator, for example via change-point models~\cite{AdamsMacKay2007}; (ii) exploring more computationally efficient inference alternatives to support the use of our model in real-time, building scale, digital-twin applications; (iii) integrating richer physics into the model, to reduce structural mismatches while still maintaining computational efficiency.

%=====================================================

\section*{Acknowledgments}
The authors acknowledge the ACD (Academic Committee for Design) of  the Department of Mechanical Engineering, Eindhoven University of Technology (TU/e). The authors gratefully acknowledge the support of colleagues in the TU/e BPS Laboratory during the setup of the experiments. The authors also thank Dr. Marcel Loomans for his constructive and valuable suggestions.

\bibliographystyle{ieeetr}
%\bibliography{cm_references}
\bibliography{MyBib}

\end{document}